\RequirePackage{fix-cm}
\documentclass[smallextended]{svjour3}       % onecolumn (second format)
\smartqed  % flush right qed marks, e.g. at end of proof
\usepackage{graphicx}
\usepackage{graphicx} % Required for inserting images

\usepackage{multirow}%
\usepackage{amsmath,amssymb,amsfonts}%
\usepackage{mathrsfs}%
\usepackage[title]{appendix}%
\usepackage{xcolor}%
\usepackage{textcomp}%
\usepackage{manyfoot}%
\usepackage{booktabs}%
\usepackage{algorithm}%
\usepackage{algorithmicx}%
\usepackage{algpseudocode}%
\usepackage{listings}%

\usepackage[utf8]{inputenc}

\usepackage{multirow}

\usepackage{amssymb}
\usepackage{latexsym}
\usepackage{amsmath,amstext,amssymb,bbm,amsbsy}
\usepackage{mathtools}
\usepackage{tabularx}
\usepackage{url}

\usepackage{subcaption}

\usepackage{tikz}
\usetikzlibrary{shapes.geometric, arrows}
\usepackage{adjustbox}
\usepackage{siunitx}

\newtheorem{rem}{Remark}[section]

\DeclareMathOperator*{\argmin}{arg\,min}
\DeclareMathOperator*{\argmax}{arg\,max}

\def\R {\mathbb{R}}

\def\E {\mathbb{E}}
\def\dom {\mathrm{dom~}}

\def\Rn {\R^n}

\def\ind {\chi}  % indicator function (values 0 or +infty)
\begin{document}

\title{A Primal-dual hybrid gradient method for solving optimal control problems and the corresponding Hamilton-Jacobi PDEs\thanks{T. Meng, S. Liu and S. Osher are partially funded by Air Force Office of Scientific Research (AFOSR) MURI FA9550-18-502 and Office of Naval Research (ONR) N00014-20-1-2787. W. Li is supported by Air Force Office of Scientific Research (AFOSR) MURI FP 9550-18-502, Air Force Office of Scientific Research (AFOSR) YIP award No. FA9550-23-1-0087, National Science Foundation (NSF) DMS-2245097, and National Science Foundation (NSF) RTG: 2038080.}
}
% \subtitle{Do you have a subtitle?\\ If so, write it here}


\author{Tingwei Meng         \and
        Siting Liu           \and 
        Wuchen Li            \and
        Stanley Osher
}


\institute{Tingwei Meng \at
              Mathematics Department, University of California, Los Angeles \\
              \email{tingwei@math.ucla.edu}      
           \and
           Siting Liu \at
              Mathematics Department, University of California, Los Angeles \\
              \email{siting6@math.ucla.edu}
           \and
           Stanley Osher \at
              Mathematics Department, University of California, Los Angeles \\
              \email{sjo@math.ucla.edu}
           \and
           Wuchen Li \at
              Mathematics Department, University of South Carolina \\
              \email{wuchen@mailbox.sc.edu}
}

% \affil[1]{\orgdiv{Mathematics Department}, \orgname{University of California, Los Angeles}, \orgaddress{\street{520 Portola Plaza}, \city{Los Angeles}, \postcode{90095}, \state{CA}, \country{USA}}}

% \affil[2]{\orgdiv{Mathematics Department}, \orgname{University of South Carolina}, \orgaddress{\street{1523 Greene Street}, \city{Columbia}, \postcode{29208}, \state{SC}, \country{USA}}}

\date{Received: date / Accepted: date}

\maketitle

\begin{abstract}
Optimal control problems are crucial in various domains, including path planning, robotics, and humanoid control, demonstrating their broad applicability. The connection between optimal control and Hamilton-Jacobi (HJ) partial differential equations (PDEs) underscores the need for solving HJ PDEs to address these control problems effectively. While numerous numerical methods exist for tackling HJ PDEs across different dimensions, this paper introduces an innovative optimization-based approach that reformulates optimal control problems and HJ PDEs into a saddle point problem using a Lagrange multiplier. Our method, based on the preconditioned primal-dual hybrid gradient (PDHG) method, offers a solution to HJ PDEs with first-order accuracy and numerical unconditional stability, enabling larger time steps and avoiding the limitations of explicit time discretization methods. 
Our approach has ability to handle a wide variety of Hamiltonian functions, including those that are non-smooth and dependent on time and space, through a simplified saddle point formulation that facilitates easy and parallelizable updates. Furthermore, our framework extends to viscous HJ PDEs and stochastic optimal control problems, showcasing its versatility. Through a series of numerical examples, we demonstrate the method's effectiveness in managing diverse Hamiltonians and achieving efficient parallel computation, highlighting its potential for wide-ranging applications in optimal control and beyond.
\end{abstract}

\section{Introduction}
Optimal control problems are pivotal in numerous practical applications, such as trajectory planning~\cite{Coupechoux2019Optimal,Rucco2018Optimal,Hofer2016Application,Delahaye2014Mathematical,Parzani2017HJB,Lee2021Hopf}, robot manipulator control~\cite{lewis2004robot,Jin2018Robot,Kim2000intelligent,Lin1998optimal,Chen2017Reachability}, and humanoid robot control~\cite{Khoury2013Optimal,Feng2014Optimization,kuindersma2016optimization,Fujiwara2007optimal,fallon2015architecture,denk2001synthesis}, highlighting their significance across various fields. The relationship between optimal control and Hamilton-Jacobi (HJ) partial differential equations (PDEs) is well-established (see~\cite{Bardi1997Optimal,Yong1999Stochastic}), underscoring the importance of solving HJ PDEs for addressing optimal control challenges. Various numerical strategies have been developed to tackle HJ PDEs and optimal control problems, ranging from high-order grid-based methods in lower-dimensional settings to innovative approaches designed to overcome the curse of dimensionality in higher dimensions. 
For lower-dimensional problems, advanced grid-based methods like essentially nonoscillatory (ENO) schemes~\cite{Osher1991High}, weighted ENO schemes~\cite{Jiang2000Weighted}, and the discontinuous Galerkin method~\cite{Hu1999Discontinuous} are frequently utilized. To tackle the complexity of higher dimensions, several sophisticated strategies have been introduced.
These approaches encompass a diverse range of methodologies, including max-plus algebra methods~\cite{mceneaney2006max,akian2006max,akian2008max,dower2015max,Fleming2000Max,gaubert2011curse,McEneaney2007COD,mceneaney2008curse,mceneaney2009convergence}, dynamic programming and reinforcement learning~\cite{alla2019efficient,bertsekas2019reinforcement}, tensor decomposition techniques~\cite{dolgov2019tensor,horowitz2014linear,todorov2009efficient}, sparse grids~\cite{bokanowski2013adaptive,garcke2017suboptimal,kang2017mitigating}, model order reduction~\cite{alla2017error,kunisch2004hjb}, polynomial approximation~\cite{kalise2019robust,kalise2018polynomial}, optimization methods~\cite{darbon2015convex,darbon2016splitting,kirchner2017time,kirchner2018primal,darbon2019decomposition,Darbon2016Algorithms,yegorov2017perspectives,chow2017algorithm,chow2019algorithm,lin2018splitting,chen2021hopf,chen2021lax}, and neural network-based solutions~\cite{darbon2019overcoming,bachouch2018deep,Djeridane2006Neural,jiang2016using,Han2018Solving,hure2018deep,hure2019some,lambrianides2019new,Niarchos2006Neural,reisinger2019rectified,royo2016recursive,Sirignano2018DGM,darbon2021some,darbon2023neural,zhou2021actor,zhou2023policy}.

This paper introduces a novel optimization-based framework with time-implicit discretization to solve optimal control problems and HJ PDEs by reformulating them into a saddle point problem using a Lagrange multiplier. We propose an algorithm based on the preconditioned primal-dual hybrid gradient (PDHG) method~\cite{chambolle2011first} for finding the saddle point, effectively solving the HJ PDEs. Unlike other grid-based methods, our approach achieves first-order accuracy and numerical unconditional stability through implicit time discretization, allowing for larger time steps and avoiding the Courant–Friedrichs–Lewy (CFL) condition limitation of explicit time discretization methods.

Our method distinguishes itself within optimization-based approaches by its capability to handle a broad set of Hamiltonians, notably including those that exhibit non-smooth behavior and dependencies on both time and space variables. The saddle point formulation simplifies iteration updates, eliminating the need for complex nonlinear inversion. The most intricate update step can be executed efficiently and in parallel through the proximal point operator. For scenarios frequently encountered that involve quadratic or $1$-homogeneous Hamiltonians, possibly dependent on $(x,t)$, the proximal point operator can be straightforwardly implemented using explicit expressions.
We provide a convergence analysis of the PDHG algorithm under certain conditions and establish connections with traditional numerical schemes for HJ PDEs. Furthermore, our approach is applicable to viscous HJ PDEs and corresponding stochastic optimal control problems, demonstrating its versatility and broad applicability.

The technique of converting PDE problems into saddle point problems and solving them with the PDHG method has been effectively applied to reaction-diffusion equations~\cite{liu2023firstorder,fu2023high,carrillo2023structure}, conservation laws~\cite{liu2023primal}, and HJ PDEs~\cite{meng2023primal}. While the reaction-diffusion and conservation law studies employ PDHG alongside integration by parts for solving initial value problems, this method is not readily applicable to HJ PDEs without modifications. The approach for HJ PDEs in~\cite{meng2023primal} employs the Fenchel-Legendre transform within its saddle point formulation, leading to an objective function that is linear with respect to the HJ PDE solution $\varphi$, thereby avoiding nonlinear updates on $\varphi$. This method also facilitates the use of larger time steps thanks to implicit temporal Euler discretization. However, extracting optimal controls directly from HJ PDE solutions necessitates additional steps. Expanding upon this groundwork, our research specifically focuses on optimal control problems by proposing a dedicated saddle point formulation and algorithm. Unlike the previous work, our saddle point formulation is derived from the intrinsic link between optimal control problems and HJ PDEs, not the Fenchel-Legendre transform. This adjustment retains the advantages of the earlier method while making it more appropriate for optimal control issues. Our method's simplicity in the saddle point formulation allows for updates that are either explicitly defined or amenable to parallel computation, thus broadening its utility.

In addition, our proposed saddle point formulation establishes a connection with the primal-dual formulation of potential mean-field games (MFGs), a concept introduced in \cite{lasry2007mean,huang2006large} for modeling the equilibrium behavior of large populations of agents engaged in strategic interactions. This framework has been applied across a variety of fields~\cite{lachapelle2011mean,gomes2014mean,Gomes2015,carmona2015mean,cardaliaguet2018mean,kizilkale2019integral,han2019mean,lee2021controlling}. An MFG system is characterized by coupled PDEs, comprising a forward-evolving Fokker-Planck equation and a backward-evolving HJ equation.
In our research, we utilize this interconnection to approach the solution of HJ PDEs. By drawing on this link, we can apply methods initially devised for MFGs to the task of solving HJ PDEs that have specific initial conditions. 
For example, the PDHG method, which has been used to solve potential MFGs through a saddle point formulation~\cite{papadakis2014optimal,briceno2018proximal,briceno2019implementation}, is adapted in this work to solve optimal control problems and HJ PDEs.

To evaluate the effectiveness of our proposed methodology, we present a series of numerical examples in both one-dimensional and two-dimensional settings\footnote{Codes are available at \url{https://github.com/TingweiMeng/PDHG-optimal-control}.}. 
These examples demonstrate the method's capability in handling diverse Hamiltonians that depend on the spatial variable.
A notable feature of our method is that, during each iteration, the updates at individual points are decoupled, facilitating the application of parallel computing techniques to expedite the computation process. Additionally, we showcase an example where the Lagrangian is defined by an indicator function, imposing control constraints without introducing running costs. This scenario yields a Hamiltonian $H(x,t,p)$ that is $1$-homogeneous with respect to $p$ for any given spatial variable $x$ and temporal variable $t$. The numerical outcomes reveal the emergence of discontinuous feedback control functions and bang-bang control strategies. Such phenomena are also observed in the stochastic optimal control problem variants, where, despite the smoothness of the solution to the viscous HJ PDE, the control function remains non-smooth, with evident discontinuities in the sampled open-loop control trajectories.

This paper is structured as follows: Section~\ref{sec:oc_problem} lays out the saddle point problem and introduces our algorithm tailored for optimal control problems and their related HJ PDEs. We delve into the formulation and algorithm within function spaces in Section~\ref{sec:saddle_point_cont} and proceed to discuss spatial and temporal discretization strategies for one-dimensional cases in Section~\ref{sec:discretization_det_1d} and two-dimensional cases in Section~\ref{sec:discretization_det_2d}. Several numerical experiments are shown in Section~\ref{sec:numerics_det}. Section~\ref{sec:stochastic_oc} proposes similar saddle point formulation and algorithm for stochastic optimal control problems and the corresponding viscous HJ PDEs, with further numerical experiments presented in Section~\ref{sec:numerics_sto}. The paper wraps up with a summary and potential directions for future work in Section~\ref{sec:summary}. The appendices provide in-depth technical insights: Appendix~\ref{appendix:conv} examines our algorithm's convergence analysis for specific scenarios, and Appendix~\ref{sec:appendix_consistency} connects our methodology to traditional numerical solvers for HJ PDEs. Discretization details for the saddle point formulation and algorithm applied to viscous HJ PDEs and stochastic control problems are outlined in Appendix~\ref{sec:discretization_soc}.

\section{Optimal control problems}\label{sec:oc_problem}

The optimal control problems are in the following form
\begin{equation}\label{eqt:oc_problem}
{\small
\begin{split}
    \min_{\alpha(\cdot)}\left\{\int_t^T
    L(\gamma(s),s, \alpha(s)) ds + g(\gamma(T)) \colon \gamma(t)=x, \dot{\gamma}(s) = f(\gamma(s),s,\alpha(s)) \forall s\in(t,T) \right\},
\end{split}
}
\end{equation}
where $\alpha\colon [t,T]\to \R^m$ is called control, $\gamma\colon [t,T]\to \Omega$ is called state ($\Omega\subseteq \Rn$ is the state space), $f\colon \Omega\times [0,T]\times \R^m\to \Rn$ is the source term of ODE constraint, $L\colon \Omega\times [0,T]\times \R^m\to \R$ is called Lagrangian, and $g\colon \Omega \to \R$ is a terminal cost. Assume these functions are all continuous and satisfy certain conditions such that the optimal control problem has a unique minimizer. 

Denote the optimal value in~\eqref{eqt:oc_problem} by $\phi(x,t)$. Then the function $\phi$ solves the following Hamilton-Jacobi (HJ) PDE
\begin{equation*}
\begin{dcases}
\frac{\partial \phi(x,t)}{\partial t} + \inf_{\alpha\in \R^m} \{\langle f(x, t, \alpha), \nabla_x \phi(x,t)\rangle + L(x,t,\alpha)\} = 0, & x\in \Omega, t\in [0,T],\\
\phi(x,T) = g(x), & x\in \Omega,
\end{dcases}
\end{equation*}
where the terminal condition is given by the terminal cost $g$ in~\eqref{eqt:oc_problem}, and $\Omega$ is the state space in the optimal control problem. 
Throughout this paper, we predominantly consider $\Omega$ to be a rectangular domain subject to periodic boundary conditions. Nonetheless, our methodology is adaptable to alternative boundary conditions, as illustrated in Section~\ref{sec:numerical_newton_det}.
Utilizing a time reversal strategy and setting $\varphi(x,t) = \phi(x,T-t)$, we derive the initial value HJ PDE as follows:
\begin{equation}\label{eqt:cont_HJ_initial}
\begin{dcases}
\frac{\partial \varphi(x,t)}{\partial t} + H(x,t,\nabla_x \varphi(x,t)) = 0, & x\in \Omega, t\in [0,T],\\
\varphi(x,0) = g(x), & x\in \Omega,
\end{dcases}
\end{equation}
where the Hamiltonian $H$ is defined by
\begin{equation}\label{eqt:def_H}
H(x,t,p) = \sup_{\alpha\in \R^m} \{-\langle f(x, T-t, \alpha), p\rangle - L(x,T-t,\alpha)\}.
\end{equation}
This Hamiltonian $H$ is convex with respect to $p$ for any $x\in\Omega$ and $t\in [0,T]$.
From the HJ PDE~\eqref{eqt:cont_HJ_initial}, we obtain the function $\alpha\colon \Omega\times [0,T] \to \R^m$ by
\begin{equation}\label{eqt:feedback_initial}
\alpha(x,t) = \argmax_{\alpha\in \R^m} \{-\langle f(x, T-t, \alpha), \nabla_x \varphi(x,t)\rangle - L(x,T-t,\alpha)\}.
\end{equation}
Following the application of time reversal, this function yields the feedback control function represented by $\alpha(x,T-t)$.
The optimal trajectory $\gamma^*$ and optimal open-loop control $\alpha^*$ in the optimal control problem~\eqref{eqt:oc_problem} are then computed by
\begin{equation}\label{eqt:ode_openloop}
\begin{dcases}
\dot \gamma^*(s) = f(\gamma^*(s), s, \alpha(\gamma^*(s), T-s)), & s\in (t,T),\\
\gamma^*(t) = x, \\
\alpha^*(s) = \alpha(\gamma^*(s), T-s), & s\in [t,T].
\end{dcases}
\end{equation}
For additional insights into the mathematical relationship between optimal control problems and HJ PDEs, we refer readers to~\cite{Bardi1997Optimal}.

\subsection{Saddle point formulation}\label{sec:saddle_point_cont}
We introduce a saddle point formulation to solve the optimal control problem~\eqref{eqt:oc_problem} and its associated initial value HJ PDE~\eqref{eqt:cont_HJ_initial}. Following this, we propose a numerical algorithm derived from the PDHG method~\cite{chambolle2011first,valkonen2014primal,darbon2021accelerated} to effectively solve the saddle point problem.
To simplify the notation, we will use $f_{x,t}$ and $L_{x,t}$ to denote the functions $f(x,T-t,\cdot)$ and $L(x,T-t,\cdot)$, respectively.

We treat the HJ PDE given in~\eqref{eqt:cont_HJ_initial} as a constraint within an optimization problem, where the objective function is defined as $-c\int_{\Omega} \varphi(x,T) dx$. This specific objective function is adopted based on insights from numerical experiments detailed in~\cite{meng2023primal}. By employing the Lagrange multiplier $\rho$, we derive:
\begin{equation*}
\resizebox{\hsize}{!}{ $
\begin{aligned}
&\min_{\varphi \text{ satisfying~\eqref{eqt:cont_HJ_initial}}}  -c\int_{\Omega} \varphi(x,T) dx\\
=& \min_{\substack{\varphi\\ \varphi(x,0)=g(x)}}\max_{\rho}  \int_0^T \int_\Omega \rho(x,t)\left(\frac{\partial \varphi(x,t)}{\partial t} + \sup_{\alpha\in \R^m} \{-\langle f_{x,t}(\alpha), \nabla_x \varphi(x,t)\rangle - L_{x,t}(\alpha)\}\right)dxdt - c\int_\Omega \varphi(x,T)dx\\
\geq & \min_{\substack{\varphi\\ \varphi(x,0)=g(x)}}\max_{\rho\geq 0, \alpha}  \int_0^T \int_\Omega \rho(x,t)\left(\frac{\partial \varphi(x,t)}{\partial t} -\langle f_{x,t}(\alpha(x,t)), \nabla_x \varphi(x,t)\rangle - L_{x,t}(\alpha(x,t)) \right)dxdt - c\int_\Omega \varphi(x,T)dx\\
=& \min_{\substack{\varphi\\ \varphi(x,0)=g(x)}}\max_{\rho\geq 0, \alpha}  \int_0^T \int_\Omega -\frac{\partial \rho(x,t)}{\partial t}\varphi(x,t) + \nabla_x\cdot(f_{x,t}(\alpha(x,t))\rho(x,t))\varphi(x,t) - \rho(x,t)L_{x,t}\left(\alpha(x,t)\right) dxdt \\
&\quad\quad\quad\quad - c\int_\Omega \varphi(x,T)dx + \int_\Omega (\rho(x,T)\varphi(x,T) - \rho(x,0)g(x))dx,
\end{aligned} 
$}
\end{equation*}
where the boundary term over $\partial \Omega$ vanishes during integration by parts, owing to the periodic boundary condition imposed on $\Omega$. While this approach can be adapted to different boundary conditions, it necessitates a thorough evaluation of the boundary terms specific to each case.
In this study, we introduce algorithms designed to solve HJ PDEs and their related optimal control problems through solving the subsequent saddle point problem
\begin{equation} \label{eqt:saddle_oc}
\begin{adjustbox}{width=0.99\textwidth}$
\begin{split}
\min_{\substack{\varphi\\ \varphi(x,0)=g(x)}}\max_{\rho\geq 0, \alpha}  \int_0^T \int_\Omega \rho(x,t)\left(\frac{\partial \varphi(x,t)}{\partial t} -\langle f_{x,t}(\alpha(x,t)), \nabla_x \varphi(x,t)\rangle - L_{x,t}(\alpha(x,t)) \right)dxdt - c\int_\Omega \varphi(x,T)dx.
\end{split}
$\end{adjustbox}
\end{equation}

Assuming that the stationary point $(\varphi, \rho, \alpha)$ satisfies $\rho > 0$, the first-order optimality condition is
\begin{equation}\label{eqt:first_order_optimality_cont}
\begin{dcases}
\partial_t \varphi(x,t) - \langle f_{x,t}(\alpha(x,t)), \nabla_x \varphi(x,t)\rangle - L_{x,t}(\alpha(x,t)) = 0,\\
\alpha(x,t) = \argmin_{\alpha\in\R^m}\{\langle f_{x,t}(\alpha), \nabla_x \varphi(x,t)\rangle + L_{x,t}(\alpha) \},\\
\partial_t \rho(x,t) - \nabla_x \cdot (f_{x,t}(\alpha(x,t))\rho(x,t)) = 0,\\
\varphi(x,0) = g(x),\quad  \rho(x,T) = c.
\end{dcases}
\end{equation}
Upon applying the Hamiltonian $H$, as defined in~\eqref{eqt:def_H}, the resulting simplification leads to:
\begin{equation}\label{eqt:coupled_pdes_det}
\begin{dcases}
\partial_t \varphi(x,t) + H(x,t, \nabla_x \varphi(x,t)) = 0,\\
\partial_t \rho(x,t) + \nabla_x \cdot (\nabla_p H(x,t,\nabla_x \varphi(x,t))\rho(x,t)) = 0,\\
\varphi(x,0) = g(x),\quad  \rho(x,T) = c,
\end{dcases}
\end{equation}
yielding a system containing the HJ PDE~\eqref{eqt:cont_HJ_initial} with a continuity equation. Furthermore, the stationary point $\alpha$ determines the function specified in~\eqref{eqt:feedback_initial}, which facilitates the calculation of the optimal control. Thus, addressing the saddle point problem~\eqref{eqt:saddle_oc} enables the derivation of the HJ PDE solution $\varphi$ along with the feedback optimal control function via $\alpha$. While theoretically, a discrepancy exists between the saddle point problem's stationary point $(\varphi, \rho, \alpha)$ and the HJ PDE solution if $\rho$ takes zero values within certain regions of $\Omega\times [0,T]$, such instances were not encountered during our numerical experiments.

To solve the saddle point problem~\eqref{eqt:saddle_oc}, we propose an iterative method, with the update at the $\ell$-th iteration described as follows:
\begin{equation} \label{eqt:pdhg_det_cont}
\begin{adjustbox}{width=0.99\textwidth}$
\begin{dcases}
\rho^{\ell+1}(x,t) = \left(\rho^{\ell}(x,t) + \tau_\rho\left(\partial_t \tilde\varphi^\ell(x,t) - \langle f_{x,t}(\alpha^\ell(x,t)), \nabla_x \tilde\varphi^\ell(x,t)\rangle - L_{x,t}(\alpha^\ell(x,t))\right)\right)_+,\\
\alpha^{\ell+1}(x,t) = \argmin_{\alpha\in\R^m}\left\{\langle f_{x,t}(\alpha), \nabla_x \tilde\varphi^\ell(x,t)\rangle + L_{x,t}(\alpha) + \frac{\rho^{\ell+1}(x,t)}{2\tau_\alpha} \|\alpha - \alpha^\ell(x,t)\|^2\right\},\\
\varphi^{\ell+1}(x,t) = \varphi^{\ell}(x,t) + \tau_{\varphi}(I-\Delta)^{-1}\left(\partial_t \rho^{\ell+1}(x,t) - \nabla_x \cdot (f_{x,t}(\alpha^{\ell+1}(x,t))\rho^{\ell+1}(x,t))\right),\\
\tilde\varphi^{\ell+1} = 2\varphi^{\ell+1} - \varphi^\ell.
\end{dcases}
$\end{adjustbox}
\end{equation}
It is pertinent to highlight that both the saddle point problem and our algorithm exclusively engage with affine functions of $\nabla_x \varphi$, from both an analytical and numerical standpoint. This approach strategically circumvents the nonlinear aspects inherent in HJ PDEs.
The methodology we propose utilizes a preconditioned PDHG approach (for an extensive discussion on PDHG, the reader is referred to~\cite{chambolle2011first,valkonen2014primal,darbon2021accelerated}). Further elaboration on this method is provided in the subsequent remark and detailed in Appendix~\ref{appendix:conv}.

\begin{remark}[Convergence of the Proposed Algorithm] \label{rem:conv_pdhg_det}
In Appendix~\ref{appendix:conv}, we delve into the specifics of the convergence of the PDHG algorithm when applied to~\eqref{eqt:saddle_oc} after a change of variable (see~\eqref{eqt:pdhg_det_cont_appendix}). Specifically, our discussion is limited to scenarios where $f_{x,t}$ is an affine function and \(L_{x,t}\) is non-negative, proper, convex, lower semi-continuous, and \(1\)-coercive with $L_{x,t}(0) = 0$. While this argument can extend to broader cases, we refrain from addressing the most general scenarios as our current assumptions encompass several significant instances, including quadratic Lagrangians or indicator functions of a set that includes \(0\). Such instances are relevant to quadratic Hamiltonians and \(1\)-homogeneous Hamiltonians. For illustrative examples, refer to Section~\ref{sec:numerics_det}.

In Remark~\ref{rem:A1_conv_pdhg_rho_alp}, we clarify the connection between the update methods detailed in~\eqref{eqt:pdhg_det_cont_appendix} and the updates we propose in~\eqref{eqt:pdhg_det_cont}. Our approach revises~\eqref{eqt:pdhg_det_cont_appendix} by opting for a sequential update strategy for \(\rho\) and \(\alpha\), as opposed to their concurrent modification. This strategy leverages the convergence principles outlined in Appendix~\ref{appendix:conv}, which pertain to the joint updates in~\eqref{eqt:pdhg_det_cont_appendix}. Based on that, the sequential update scheme in~\eqref{eqt:pdhg_det_cont} facilitates a more computationally streamlined and versatile algorithm. To address the errors introduced by the staggered updates for \(\rho\) and \(\alpha\), our implementation includes multiple iterative updates for \(\rho\) and \(\alpha\) followed by a single update for \(\varphi\).
\end{remark}

\begin{remark}[Details on Implementation]
For the variable $\alpha$, we employ the penalty term $\frac{\rho^{k+1}(x,t)}{2\tau_\alpha} \|\alpha - \alpha^k(x,t)\|^2$, while for $\varphi$, the penalty is based on the $H^1$-norm, leading to the inclusion of the preconditioning operator $(I-\Delta)^{-1}$ in $\varphi$'s update equation. These penalty terms are specifically selected to enhance the speed of the algorithm. The use of preconditioning methods is a well-established practice in the field, with references such as~\cite{jacobs2019solving,liu2023primal,meng2023primal} providing further context. For an in-depth discussion on the practical aspects of computing $(I-\Delta)^{-1}$, we direct the reader to~\cite{meng2023primal}.

In cases where the grid is extensively partitioned, leading to a high-dimensional optimization challenge, a reduction in dimensionality can be crucial for feasible numerical resolution. To this end, we implemented a strategy that divides the time interval, applying our algorithm within each segment sequentially. Additional information on this approach can be found in~\cite[Algorithm 3]{liu2023primal} and~\cite{meng2023primal}.
\end{remark}

\begin{remark}[Relation to MFGs]\label{rem:connection_mfg}
Upon establishing the stationary point $(\varphi, \rho, \alpha)$ where $\rho > 0$, we derive the coupled PDE system as delineated in~\eqref{eqt:coupled_pdes_det}. Implementing time reversal transforms the equations for $\varphi(x,T-t)$ and $\tilde \rho(x,t) = \rho(x,T-t)$, aligning them with the first-order optimality conditions of MFGs, which are expressed as follows:
\begin{equation*}
\begin{split}
\min_{\tilde\alpha}\Big\{\int_0^T\int_{\Omega}
L(x,s, \tilde\alpha(x,s))\tilde \rho(x,s) dxds + \int_\Omega g(x)\tilde \rho(x,T) dx \colon \quad\quad\quad\quad\\
\partial_t\tilde \rho(x,s) + \nabla_x \cdot (f(x,s,\tilde\alpha(x,s))\tilde\rho(x,s)) = 0 \ 
 \forall x\in\Omega, s\in [0,T],\\
 \tilde \rho(x,0) = c \ \forall x\in \Omega\Big\}.
\end{split}
\end{equation*}
This linkage emerges as inherently logical, given that an optimal control problem is fundamentally intertwined with an MFG scenario characterized by an initial density condition resembling a Dirac mass. Through this association, numerical methods devised for MFGs can also be utilized in addressing HJ PDEs and associated optimal control problems.
\end{remark}

\begin{remark}[Comparison with the Approach in~\cite{meng2023primal}]
While both the current work and~\cite{meng2023primal} leverage saddle point formulations to tackle HJ PDEs, distinctions arise particularly in the context of HJ PDEs derived from optimal control problems. The following analysis addresses scenarios where $f_{x,t}$ is a bijection for any $(x,t)$ within $\Omega\times [0,T]$.
Through implementing a change of variable $v = -f(x,T-t,\alpha)$, we obtain
\begin{equation*}
\begin{split}
H(x,t,p) &= \sup_{\alpha\in \R^m} \{-\langle f(x, T-t, \alpha), p\rangle - L(x,T-t,\alpha)\}\\
&= \sup_{v\in \R^n} \{\langle v, p\rangle - L(x,T-t,f^{-1}(x,T-t,-v))\},
\end{split}
\end{equation*}
with the inversion $f^{-1}$ applied solely concerning $\alpha$. Consequently, $H(x,t,\cdot)$'s Fenchel-Legendre transform emerges as the convex lower semi-continuous hull of $v\mapsto L(x,T-t,f^{-1}(x,T-t,-v))$.

Assuming the function $v\mapsto L(x,T-t,f^{-1}(x,T-t,-v))$ is convex and lower semi-continuous, the saddle point formulation presented in~\cite{meng2023primal} is:
\begin{equation}\label{eqt:saddle_point_duality}
\begin{adjustbox}{width=0.99\textwidth}$
\begin{split}
\min_{\substack{\varphi\\ \varphi(x,0)=g(x)}}\max_{\rho\geq 0, v} \int_0^T \int_\Omega \rho(x,t)\left(\frac{\partial \varphi(x,t)}{\partial t} + \langle v(x,t), \nabla_x\varphi(x,t)\rangle  - L(x,T-t, f^{-1}(x,T-t, -v(x,t)) \right)dxdt\\ - c\int_\Omega \varphi(x,T)dx.
\end{split}
$\end{adjustbox}
\end{equation}
Although both approaches employ saddle point methodologies, the specific formulations outlined in~\eqref{eqt:saddle_oc} and~\eqref{eqt:saddle_point_duality} lead to distinct algorithms. The strategy developed in this work is expressly designed for addressing optimal control problems, effectively providing solutions to HJ PDEs as well as the feedback optimal control functions crucial for determining optimal controls.

A key component of the above computation is the change of variable $-v = f(x,T-t, \alpha)$, a method also used in~\cite{lee2022convexifying}. Unlike~\cite{lee2022convexifying}, which concentrates on analyzing a single trajectory in the context of optimal control, our proposed approach provides solutions to HJ PDEs and feedback controls over the entire domain. These solutions facilitate the computation of optimal controls for various initial conditions.
\end{remark}

\subsection{Discretization}\label{sec:discretization_det}
This section is dedicated to discretizing the saddle point problem presented in~\eqref{eqt:saddle_oc} and outlining the associated algorithms for one-dimensional scenarios in Section~\ref{sec:discretization_det_1d} and for two-dimensional situations in Section~\ref{sec:discretization_det_2d}.
For each scenario, we start with upwind spatial discretization, an essential step for correctly approximating the viscous solutions of the HJ PDEs. Following this, implicit Euler temporal discretization is applied, enabling us to circumvent the restrictive CFL condition often encountered with explicit schemes.
\subsubsection{One-dimensional problems}\label{sec:discretization_det_1d}
\textbf{Spatial discretization.}
The first step in discretizing the saddle point problem involves partitioning the spatial domain $\Omega$. We define $\{x_i\}_{i=1}^{n_x}$ as the evenly spaced grid points within $\Omega$, with $\Delta x$ representing the distance between adjacent points. For a domain $\Omega = [a,b]$ (where $-\infty < a < b < \infty$) subject to periodic boundary conditions, we calculate the grid size as $\Delta x = \frac{b-a}{n_x}$ and assign grid points such that $x_i = a + (i-1)\Delta x$. Our objective is to determine the values of the functions at these grid points, specifically $\varphi_i(t) = \varphi(x_i,t)$ and $\rho_i(t) = \rho(x_i,t)$.
In addressing the HJ PDE, the treatment of first-order spatial derivatives necessitates consideration of both positive and negative velocities, corresponding to scenarios where the function $f$ assumes positive and negative values, respectively. To manage these cases, we introduce two variables, $\alpha_1$ and $\alpha_2$, and aim to calculate their values at the grid points, denoted as $\alpha_{1,i}(t)$ and $\alpha_{2,i}(t)$. 
Additionally, the numerical Lagrangian $\hat L\colon \Omega\times [0,T]\times \R^m\times \R^m\to \R$ is utilized as an approximation of $L$, and its determination will vary across different scenarios. The selection of $\hat L$ for particular cases will be discussed subsequently.

For simplicity, we represent the functions $f(x_i, T-t,\cdot)$ and $\hat L(x_i, T-t,\cdot)$ as $f_{i,t}$ and $\hat L_{i,t}$ respectively. Upon applying first-order spatial discretization, the saddle point problem expressed in~\eqref{eqt:saddle_oc} (after division by $\Delta x$) is reformulated as:
\begin{equation}\label{eqt:saddle_point_det_1d_semi}
{\scriptsize
\begin{split}
\min_{\substack{\varphi\\ \varphi_i(0)=g(x_i)}}\max_{\rho\geq 0, \alpha}  \int_0^T \sum_{i=1}^{n_x} \rho_i(t)\Bigg(\dot\varphi_i(t) - f_{i,t}(\alpha_{1,i}(t))_+ (D_x^+\varphi)_i(t) - f_{i,t}(\alpha_{2,i}(t))_- (D_x^-\varphi)_i(t) \\
- \hat L_{i,t}\left(\alpha_{1,i}(t), \alpha_{2,i}(t)\right) \Bigg)dt - c\sum_{i=1}^{n_x} \varphi_i(T).
\end{split}
}
\end{equation}
The variables of this saddle point problem contain all semi-discretized functions, including $\varphi_i(t)$, $\rho_i(t)$, $\alpha_{1,i}(t)$, and $\alpha_{2,i}(t)$, across the time interval $t\in[0,T]$ and for grid indices $i=1,\dots, n_x$.
Here, $f_+ = \max\{f,0\}$ and $f_- = \min\{f,0\}$ denote the positive and negative components of $f$, respectively. Meanwhile, $(D_x^+ \varphi)_i = \frac{\varphi_{i+1} - \varphi_i}{\Delta x}$ and $(D_x^- \varphi)_i = \frac{\varphi_{i} - \varphi_{i-1}}{\Delta x}$ represent the right and left finite differences, respectively. In this paper, for the sake of simplicity, we interchangeably use the notations $(D \eta)_i$ and $D(\eta_i)$ for any finite difference operator $D$, function $\eta$, and index $i$, regardless of whether they pertain to one-dimensional or two-dimensional scenarios.

Examining a stationary point $(\varphi, \rho, \alpha)$ within this saddle point problem, and under the assumption that $\rho_i(t) > 0$ for any $t\in [0,T]$ and $i=1,\dots, n_x$, we arrive at the first-order optimality conditions as follows: 
\begin{equation}\label{eqt:first_order_optimality_semi1d}
\begin{adjustbox}{width=0.99\textwidth}$
\begin{dcases}
\dot{\varphi}_i(t) - f_{i,t}(\alpha_{1,i}(t))_+ (D_x^+ \varphi)_i(t) - f_{i,t}(\alpha_{2,i}(t))_- (D_x^- \varphi)_i(t) - \hat L_{i,t}(\alpha_{1,i}(t), \alpha_{2,i}(t)) = 0,\\
(\alpha_{1,i}(t), \alpha_{2,i}(t)) = \argmin_{\alpha_1, \alpha_2\in\R^m}\{f_{i,t}(\alpha_1)_+ (D_x^+\varphi)_i(t) + f_{i,t}(\alpha_{2})_- (D_x^- \varphi)_i(t) + \hat L_{i,t}(\alpha_1, \alpha_2)\}, \\
\dot{\rho}_i(t) - D_x^-(f_{i,t}(\alpha_{1,i}(t))_+\rho_i(t)) - D_x^+(f_{i,t}(\alpha_{2,i}(t))_-\rho_i(t)) = 0.
\end{dcases}
$\end{adjustbox}
\end{equation}
By integrating the first two equations, we obtain
\begin{equation*}
\begin{adjustbox}{width=0.99\textwidth}$
\dot{\varphi}_i(t) +\sup_{\alpha_1,\alpha_2\in \R^m} \{-f_{i,t}(\alpha_{1})_+ (D_x^+ \varphi)_i(t) - f_{i,t}(\alpha_{2})_- (D_x^- \varphi)_i(t) - \hat L_{i,t}(\alpha_{1}, \alpha_{2})\} = 0,
$\end{adjustbox}
\end{equation*}
which provides a semi-discrete formulation for the HJ PDE~\eqref{eqt:cont_HJ_initial}. The corresponding numerical Hamiltonian $\hat H$ is given by
\begin{equation}\label{eqt:numericalH1d}
\begin{adjustbox}{width=0.99\textwidth}$
\hat H(x,t, D_x^+ \varphi, D_x^- \varphi) = \sup_{\alpha_1, \alpha_2\in \R^m} \{-f_{x,t}(\alpha_{1})_+ D_x^+ \varphi - f_{x,t}(\alpha_{2})_- D_x^- \varphi - \hat L_{x,t}(\alpha_{1}, \alpha_{2})\}.
$\end{adjustbox}
\end{equation}
This yields a semi-discrete approach to solving HJ PDEs and providing feedback optimal controls.

To obtain a valid numerical Hamiltonian, the following conditions must be met:
\begin{enumerate}
\item Monotonicity: Specifically, $\hat H$ should be non-increasing with respect to $D_x^+ \varphi$, and is non-decreasing with respect to $D_x^- \varphi$. This condition is ensured by the formulation given in~\eqref{eqt:numericalH1d}.
\item Consistency: This requires $\hat H(x,t, p, p) = H(x,t,p)$ for any $x\in \Omega$, $t\in [0,T]$, and $p\in\R$. The numerical Lagrangian $\hat L$ must be chosen to ensure that the corresponding numerical Hamiltonian is consistent. In scenarios where $f$ linearly depends on $\alpha$ and $L$ is a non-negative, convex, $1$-coercive function of $\alpha$ with $L_{x,t}(0) = 0$ for every $x$ and $t$, we can choose $\hat L_{x,t}(\alpha_1, \alpha_2) = L_{x,t}(\alpha_1) + L_{x,t}(\alpha_2)$. For more technical details, refer to Appendix~\ref{sec:appendix_consistency_1d}.
\end{enumerate}

With this spatial discretization, the update for the $\ell$-th iteration is described as follows:
\begin{equation}\label{eqt:pdhg_det_semi}
\begin{adjustbox}{width=0.99\textwidth}$
\begin{dcases}
\rho_i^{\ell+1}(t) = \Big(\rho_i^{\ell}(t) + \tau_\rho\big(\dot{\tilde\varphi}^\ell_i(t) - f_{i,t}(\alpha^{\ell}_{1,i}(t))_+ (D_x^+ \tilde\varphi^\ell)_i(t) - f_{i,t}(\alpha^\ell_{2,i}(t))_- (D_x^- \tilde\varphi^\ell)_i(t) - \hat L_{i,t}(\alpha^\ell_{1,i}(t), \alpha^\ell_{2,i}(t))\big)\Big)_+. \\
(\alpha^{\ell+1}_{1,i}(t), \alpha^{\ell+1}_{2,i}(t)) = \argmin_{\alpha_1, \alpha_2\in\R^m}\Big\{f_{i,t}(\alpha_1)_+ (D_x^+\tilde\varphi^\ell)_i(t) + f_{i,t}(\alpha_2)_- (D_x^-\tilde\varphi^\ell)_i(t) + \hat L_{i,t}(\alpha_1, \alpha_2)\\
\quad\quad\quad\quad\quad\quad\quad\quad
\quad\quad\quad\quad\quad\quad\quad\quad+ \frac{\rho_i^{\ell+1}(t)}{2\tau_\alpha} \left(\|\alpha_1 - \alpha^\ell_{1,i}(t)\|^2 + \|\alpha_2 - \alpha^\ell_{2,i}(t)\|^2\right)\Big\}.\\
\varphi^{\ell+1}_i(t) = \varphi^{\ell}_i(t) + \tau_{\varphi}(I-\partial_{tt}-D_{xx})^{-1}\left(\dot{\rho}^{\ell+1}_i(t) - D_x^-(f_{i,t}(\alpha^{\ell+1}_{1,i}(t))_+\rho^{\ell+1}_i(t)) - D_x^+(f_{i,t}(\alpha^{\ell+1}_{2,i}(t))_-\rho^{\ell+1}_i(t))\right).\\
\tilde\varphi^{\ell+1} = 2\varphi^{\ell+1} - \varphi^\ell,
\end{dcases}
$\end{adjustbox}
\end{equation}
where $D_{xx}$ denotes the central difference operator for the second-order derivative, defined as $(D_{xx} \varphi)_i(t) = \frac{\varphi_{i+1}(t) + \varphi_{i-1}(t) - 2 \varphi_{i}(t)}{\Delta x^2}$.

\textbf{Temporal discretization.}
The next step involves discretizing the time domain $[0,T]$. We define the uniform temporal grid points as $t_k$ with a grid size $\Delta t$, setting $\Delta t = \frac{T}{n_t-1}$ and assigning $t_i = (i-1)\Delta t$ for $i=1,\dots, n_t$.
For temporal discretization, particularly for $\partial_t \varphi$, we adopt an implicit Euler scheme, which circumvents the restrictive CFL condition encountered in explicit Euler discretization.
Mirroring our approach in spatial discretization, we simplify notation by using $f_{i,k}$ and $\hat L_{i,k}$ for the functions $f(x_i, T-t_k,\cdot)$ and $\hat L(x_i, T-t_k,\cdot)$, respectively. Our goal is to determine the values of functions at the grid points, namely $\varphi_{i,k} = \varphi(x_i, t_k)$, $\rho_{i,k} = \rho(x_i, t_k)$, $\alpha_{1,i,k} = \alpha_1(x_i, t_k)$, and $\alpha_{2,i,k} = \alpha_2(x_i, t_k)$.
Upon dividing by $\Delta t$, the saddle point problem is reformulated as:
\begin{equation} \label{eqt:saddle_det_1d_fully}
{\scriptsize
\begin{split}
\min_{\substack{\varphi\\ \varphi_{i,1}=g(x_i)}}\max_{\rho\geq 0, \alpha}  \sum_{k=2}^{n_t} \sum_{i=1}^{n_x} \rho_{i,k}\Bigg((D_t^-\varphi)_{i,k} - f_{i,k}(\alpha_{1,i,k})_+ (D_x^+\varphi)_{i,k} - f_{i,k}(\alpha_{2,i,k})_- (D_x^-\varphi)_{i,k} \\
- \hat L_{i,k}\left(\alpha_{1,i,k}, \alpha_{2,i,k}\right) \Bigg) - \frac{c}{\Delta t}\sum_{i=1}^{n_x} \varphi_{i,n_t},
\end{split}
}
\end{equation}
targeting optimization over $\varphi_{i,k}$, $\rho_{i,k}$, $\alpha_{1,i,k}$, and $\alpha_{2,i,k}$ for grid indices $i=1,\dots, n_x$ and time steps $k = 2,\dots, n_t$.
Henceforth, $D_t^-$ and $D_t^+$ will represent the implicit and explicit Euler discretizations, respectively, defined as $(D_t^-\varphi)_{i,k} = \frac{\varphi_{i,k} - \varphi_{i,k-1}}{\Delta t}$ and $(D_t^+\varphi)_{i,k} = \frac{\varphi_{i,k+1} - \varphi_{i,k}}{\Delta t}$.

The proposed update for the $\ell$-th iteration is then formulated as follows
\begin{equation*}
\begin{adjustbox}{width=0.99\textwidth}$
\begin{dcases}
\rho_{i,k}^{\ell+1} = \Big(\rho_{i,k}^{\ell} + \tau_\rho\big((D_t^-{\tilde\varphi^\ell})_{i,k} - f_{i,k}(\alpha^{\ell}_{1,i,k})_+ (D_x^+ \tilde\varphi^\ell)_{i,k} - f_{i,k}(\alpha^\ell_{2,i,k})_- (D_x^- \tilde\varphi^\ell)_{i,k} - \hat L_{i,k}(\alpha^\ell_{1,i,k}, \alpha^\ell_{2,i,k})\big)\Big)_+.\\
(\alpha^{\ell+1}_{1,i,k}, \alpha^{\ell+1}_{2,i,k}) = \argmin_{\alpha_1, \alpha_2\in\R^m}\Big\{f_{i,k}(\alpha_1)_+ (D_x^+\tilde\varphi^\ell)_{i,k} + f_{i,k}(\alpha_2)_- (D_x^-\tilde\varphi^\ell)_{i,k} + \hat L_{i,k}(\alpha_1, \alpha_2)\\
\quad\quad\quad\quad\quad\quad\quad\quad
\quad\quad\quad\quad\quad\quad\quad\quad+ \frac{\rho_{i,k}^{\ell+1}}{2\tau_\alpha} \left(\|\alpha_1 - \alpha^\ell_{1,i,k}\|^2 + \|\alpha_2 - \alpha^\ell_{2,i,k}\|^2\right)\Big\}.\\
\varphi^{\ell+1}_{i,k} = \varphi^{\ell}_{i,k} + \tau_{\varphi}(I-D_{tt}-D_{xx})^{-1}\left((D_t^+{\rho}^{\ell+1})_{i,k} - D_x^-(f_{i,k}(\alpha^{\ell+1}_{1,i,k})_+\rho^{\ell+1}_{i,k}) - D_x^+(f_{i,k}(\alpha^{\ell+1}_{2,i,k})_-\rho^{\ell+1}_{i,k})\right).\\
\tilde\varphi^{\ell+1} = 2\varphi^{\ell+1} - \varphi^\ell,
\end{dcases}
$\end{adjustbox}
\end{equation*}
where $D_{tt}$ represent the central difference operators for second-order derivative, defined as $(D_{tt}\varphi)_{i,k} = \frac{\varphi_{i,k+1} + \varphi_{i,k-1} - 2\varphi_{i,k}}{\Delta t^2}$.
The Fast Fourier Transform (FFT) is utilized to numerically calculate $(I - D_{tt} - D_{xx})^{-1}$. 
For additional information regarding the numerical updates, we refer readers to~\cite{jacobs2019solving,liu2023primal,meng2023primal}.

Note that if the numerical Lagrangian is chosen to be $L_{x,t}(\alpha_{1}) + L_{x,t}(\alpha_{2})$, then the updating step for $\alpha$ can be computed in parallel for $\alpha_{1}$ and $\alpha_{2}$.

\subsubsection{Two-dimensional problems}\label{sec:discretization_det_2d}

The approach to discretizing two-dimensional problems (where $n=2$) mirrors that of one-dimensional cases. Initially, we focus on spatial discretization, establishing an upwind scheme. Subsequently, for time discretization, we employ an implicit Euler method, facilitating the use of a larger temporal grid size $\Delta t$.

\textbf{Spatial discretization.}
In our analysis, we utilize regular grid structures within the domain $\Omega$. This study specifically focuses on rectangular domains $\Omega$ with periodic boundary conditions, although our methodology is adaptable to various boundary conditions. Suppose $\Omega = [a_1, b_1]\times [a_2, b_2]$, with grid sizes defined as $\Delta x = \frac{b_1-a_1}{n_x}$ and $\Delta y = \frac{b_2-a_2}{n_y}$. The grid points are denoted by $x_{i,j} = (a_1 + (i-1)\Delta x, a_2 + (j-1)\Delta y)$.
As with the one-dimensional scenario, handling positive and negative velocities in two dimensions necessitates multiple $\alpha$ variables. Given that $f$ outputs values in $\R^2$, we require four $\alpha$ variables, unlike the two needed for one-dimensional cases. The components of $f$ are represented by $f_1$ and $f_2$ for the first and second components, respectively.
For each component $i=1,2$, we introduce $\alpha_{i1}$ and $\alpha_{i2}$ to capture the positive and negative aspects of $f_i$. Our goal is to calculate the function values at the grid points; thus, we aim at computing $\varphi_{i,j}(t) = \varphi(x_{i,j},t)$, $\rho_{i,j}(t) = \rho(x_{i,j},t)$, and $\alpha_{IJ, i,j}(t) = \alpha_{IJ}(x_{i,j},t)$ for $I,J=1,2$.
For ease of reference, the functions $f_1(x_{i,j}, T-t,\cdot)$, $f_2(x_{i,j}, T-t,\cdot)$, and $\hat L(x_{i,j}, T-t,\cdot)$ are abbreviated as $f_{1,i,j,t}$, $f_{2,i,j,t}$, and $\hat L_{i,j,t}$, respectively. Here, $\hat L\colon \Omega\times [0,T]\times (\R^m)^4\to \R$ denotes the numerical Lagrangian, which will be defined later, similar to the approach taken in one-dimensional scenarios.
The formulation of the saddle point problem, divided by $\Delta x \Delta y$, is as follows:
\begin{equation}\label{eqt:saddle_point_det_2d_semi}
\begin{adjustbox}{width=0.99\textwidth}$
\begin{split}
\min_{\substack{\varphi\\ \varphi_{i,j}(0)=g(x_{i,j})}}\max_{\rho\geq 0, \alpha}  \int_0^T \sum_{i=1}^{n_x} \sum_{j=1}^{n_y} \rho_{i,j}(t)\Bigg(\dot\varphi_{i,j}(t) - f_{1,i,j,t}(\alpha_{11,i,j}(t))_+ (D_x^+\varphi)_i(t) 
- f_{1,i,j,t}(\alpha_{12,i,j}(t))_- (D_x^-\varphi)_i(t) \\
- f_{2,i,j,t}(\alpha_{21,i,j}(t))_+ (D_y^+\varphi)_i(t) - f_{2,i,j,t}(\alpha_{22,i,j}(t))_- (D_y^-\varphi)_i(t) \\
- \hat L_{i,j,t}\left(\alpha_{11,i,j}(t), \alpha_{12,i,j}(t), \alpha_{21,i,j}(t), \alpha_{22,i,j}(t)\right) \Bigg)dt - c\sum_{i=1}^{n_x}\sum_{j=1}^{n_y} \varphi_{i,j}(T),
\end{split} 
$\end{adjustbox}
\end{equation}
where $D_y^-$ and $D_y^+$ represent the left and right finite difference schemes in the $y$-dimension, analogous to $D_x^-$ and $D_x^+$ in the $x$-dimension.

Consider a stationary point $(\varphi, \rho, \alpha)$ in this saddle point problem. If we further assume $\rho_{i,j}(t) > 0$ for any $t\in [0,T]$, then the first order optimality condition is
\begin{equation}\label{eqt:first_order_optimality_semi2d}
\begin{adjustbox}{width=0.99\textwidth}$
\begin{dcases}
\dot{\varphi}_{i,j}(t) - f_{1,i,j,t}(\alpha_{11,i,j}(t))_+ (D_x^+ \varphi)_{i,j}(t) - f_{1,i,j,t}(\alpha_{12,i,j}(t))_- (D_x^- \varphi)_{i,j}(t) \\
\quad\quad - f_{2,i,j,t}(\alpha_{21,i,j}(t))_+ (D_y^+ \varphi)_{i,j}(t) - f_{2,i,j,t}(\alpha_{22,i,j}(t))_- (D_y^- \varphi)_{i,j}(t) \\
\quad\quad- \hat L_{i,j,t}(\alpha_{11,i,j}(t), \alpha_{12,i,j}(t), \alpha_{21,i,j}(t), \alpha_{22,i,j}(t)) = 0,\\
(\alpha_{11,i,j}(t), \alpha_{12,i,j}(t), \alpha_{21,i,j}(t), \alpha_{22,i,j}(t)) = \argmin_{\alpha_{11}, \alpha_{12},\alpha_{21}, \alpha_{22}\in\R^m}\{f_{1,i,j,t}(\alpha_{11})_+ (D_x^+\varphi)_{i,j}(t) \\
\quad\quad+ f_{1,i,j,t}(\alpha_{12})_- (D_x^- \varphi)_{i,j}(t) + f_{2,i,j,t}(\alpha_{21})_+ (D_y^+\varphi)_{i,j}(t) \\
\quad\quad + f_{2,i,j,t}(\alpha_{22})_- (D_y^- \varphi)_{i,j}(t) + \hat L_{i,j,t}(\alpha_{11}, \alpha_{12},\alpha_{21}, \alpha_{22})\}, \\
\dot{\rho}_{i,j}(t) - D_x^-(f_{1,i,j,t}(\alpha_{11,i,j}(t))_+\rho_{i,j}(t)) - D_x^+(f_{1,i,j,t}(\alpha_{12,i,j}(t))_-\rho_{i,j}(t))\\
\quad\quad - D_y^-(f_{2,i,j,t}(\alpha_{21,i,j}(t))_+\rho_{i,j}(t)) - D_y^+(f_{2,i,j,t}(\alpha_{22,i,j}(t))_-\rho_{i,j}(t))= 0.
\end{dcases}
$\end{adjustbox}
\end{equation}
Integrating the first two equations of this optimality condition yields:
\begin{equation*}
\begin{adjustbox}{width=0.99\textwidth}$
\begin{split}
\dot{\varphi}_{i,j}(t) +\sup_{\alpha_{11},\alpha_{12},\alpha_{21},\alpha_{22}\in\R^m} \{-f_{1,i,j,t}(\alpha_{11})_+ (D_x^+ \varphi)_{i,j}(t) - f_{1,i,j,t}(\alpha_{12})_- (D_x^- \varphi)_{i,j}(t)
-f_{2,i,j,t}(\alpha_{21})_+ (D_y^+ \varphi)_{i,j}(t) \\
- f_{2,i,j,t}(\alpha_{22})_- (D_y^- \varphi)_{i,j}(t) 
- \hat L_{i,j,t}(\alpha_{11},\alpha_{12},\alpha_{21},\alpha_{22})\} = 0.
\end{split}
$\end{adjustbox}
\end{equation*}
This integration results in a semi-discrete formulation for the HJ PDE~\eqref{eqt:cont_HJ_initial}, with the numerical Hamiltonian being
\begin{equation}\label{eqt:numericalH2d}
\begin{adjustbox}{width=0.99\textwidth}$
\begin{split}
\hat H(x,t, D_x^+ \varphi, D_x^- \varphi, D_y^+ \varphi, D_y^- \varphi) = \sup_{\alpha_{11}, \alpha_{12},\alpha_{21}, \alpha_{22}\in \R^m} \{-f_{1,x,t}(\alpha_{11})_+ D_x^+ \varphi - f_{1,x,t}(\alpha_{12})_- D_x^- \varphi\\
-f_{2,x,t}(\alpha_{21})_+ D_y^+ \varphi - f_{2,x,t}(\alpha_{22})_- D_y^- \varphi - \hat L_{x,t}(\alpha_{11}, \alpha_{12},\alpha_{21}, \alpha_{22})\}.
\end{split}
$\end{adjustbox}
\end{equation}
Here, $f_{1,x,t}$ and $f_{2,x,t}$ refer to the first and second components of the function $f_{x,t}$, respectively.
For the numerical Hamiltonian to be considered valid, it must meet the criteria outlined below:
\begin{enumerate}
    \item Monotonicity: Specifically, $\hat H$ should be non-increasing with respect to $D_x^+ \varphi$ and $D_y^+ \varphi$, and should be non-decreasing with respect to $D_x^- \varphi$ and $D_y^- \varphi$. This property is ensured by the formulation specified in~\eqref{eqt:numericalH2d}.
    \item Consistency: 
    It is required that $\hat H(x,t, p_1, p_1,p_2,p_2) = H(x,t,p)$ for every $x\in \Omega$, $t\in [0,T]$, and $p=(p_1,p_2)\in\R^2$. This consistency must be considered when selecting $\hat L$. 
    When $f$ exhibits a linear relationship with $\alpha$, $L$ is non-negative and convex with $1$-coerciveness over $\alpha$ and satisfies $L_{x,t}(0) = 0$, and the Hamiltonian $H$ is separable into $H_{x,t}(p_1,p_2) = H_{1}(x,t,p_1) + H_{2}(x,t,p_2)$ for some functions $H_1$ and $H_2$ applicable across all $x\in\Omega$, $t\in [0,T]$, and for any $p_1,p_2\in\R$, then the selection for $\hat L(x,t,\alpha_{11}, \alpha_{12},\alpha_{21}, \alpha_{22}) = L(x,t,\alpha_{11}) + L(x,t,\alpha_{12}) + L(x,t,\alpha_{21}) + L(x,t,\alpha_{22})$ is appropriate.
    For more technical details, refer to Appendix~\ref{sec:appendix_consistency_2d}.
\end{enumerate}

Following this discretization approach, the update for the $\ell$-th iteration is as follows: 
\begin{equation*}
\begin{adjustbox}{width=0.99\textwidth}$
\begin{dcases}
\rho_{i,j}^{\ell+1}(t) = \Big(\rho_{i,j}^{\ell}(t) + \tau_\rho\big(\dot{\tilde\varphi}^\ell_{i,j}(t) - f_{1,i,j,t}(\alpha^{\ell}_{11,i,j}(t))_+ (D_x^+ \tilde\varphi^\ell)_{i,j}(t) \\
\quad\quad\quad\quad
- f_{1,i,j,t}(\alpha^\ell_{12,i,j}(t))_- (D_x^- \tilde\varphi^\ell)_{i,j}(t) -f_{2,i,j,t}(\alpha^{\ell}_{21,i,j}(t))_+ (D_y^+ \tilde\varphi^\ell)_{i,j}(t) \\
\quad\quad\quad\quad
- f_{2,i,j,t}(\alpha^\ell_{22,i,j}(t))_- (D_y^- \tilde\varphi^\ell)_{i,j}(t) - \hat L_{i,j,t}(\alpha^\ell_{1,i}(t), \alpha^\ell_{2,i}(t))\big)\Big)_+.\\
(\alpha^{\ell+1}_{11,i,j}(t), \alpha^{\ell+1}_{12,i,j}(t), \alpha^{\ell+1}_{21,i,j}(t), \alpha^{\ell+1}_{22,i,j}(t)) = \argmin_{\alpha_{11}, \alpha_{12}, \alpha_{21}, \alpha_{22}\in\R^m}\Big\{f_{1,i,j,t}(\alpha_{11})_+ (D_x^+\tilde\varphi^\ell)_{i,j}(t) \\
\quad\quad\quad\quad
+ f_{1,i,j,t}(\alpha_{12})_- (D_x^-\tilde\varphi^\ell)_{i,j}(t) + f_{2,i,j,t}(\alpha_{21})_+ (D_y^+\tilde\varphi^\ell)_{i,j}(t) + f_{2,i,j,t}(\alpha_{22})_- (D_y^-\tilde\varphi^\ell)_{i,j}(t)\\
\quad\quad\quad\quad
 + \hat L_{i,j,t}(\alpha_{11}, \alpha_{12},\alpha_{21}, \alpha_{22}) + \frac{\rho_{i,j}^{\ell+1}(t)}{2\tau_\alpha} \big( \|\alpha_{11} - \alpha^\ell_{11,i,j}(t)\|^2 + \|\alpha_{12} - \alpha^\ell_{12,i,j}(t)\|^2\\
\quad\quad\quad\quad
 +  \|\alpha_{21} - \alpha^\ell_{21,i,j}(t)\|^2 +  \|\alpha_{22} - \alpha^\ell_{22,i,j}(t)\|^2\big)\Big\}.\\
\varphi^{\ell+1}_{i,j}(t) = \varphi^{\ell}_{i,j}(t) + \tau_{\varphi}(I-\partial_{tt}-D_{xx} - D_{yy})^{-1}\Big(\dot{\rho}^{\ell+1}_{i,j}(t) - D_x^-\left(f_{1,i,j,t}(\alpha^{k+1}_{11,i,j}(t))_+\rho^{\ell+1}_{i,j}(t)\right) \\
\quad\quad\quad\quad - D_x^+\left(f_{1,i,j,t}(\alpha^{k+1}_{12,i,j}(t))_-\rho^{\ell+1}_{i,j}(t)\right) - D_y^-\left(f_{2,i,j,t}(\alpha^{k+1}_{21,i,j}(t))_+\rho^{\ell+1}_{i,j}(t)\right) \\
\quad\quad\quad\quad - D_y^+\left(f_{2,i,j,t}(\alpha^{k+1}_{22,i,j}(t))_-\rho^{\ell+1}_{i,j}(t)\right) \Big).\\
\tilde\varphi^{\ell+1} = 2\varphi^{\ell+1} - \varphi^\ell,
\end{dcases}
$\end{adjustbox}
\end{equation*}
where $D_{yy}$ denotes the central difference operators for the second-order derivative in the $y$-dimension, analogous to the $D_{xx}$ and $D_{tt}$ operators for the $x$-dimension and $t$-dimension, respectively.

\textbf{Temporal discretization.}
We proceed to discretize the time domain, employing an implicit Euler method for $\dot \varphi$ to facilitate the use of larger time steps. The time step size is defined as $\Delta t = \frac{T}{n_t-1}$, with grid points marked by $t_k = (k-1)\Delta t$ for $k=1,\dots, n_t$.
The functions $f_I(x_{i,j}, T-t_k,\cdot)$ (for $I=1,2$) and $\hat L(x_{i,j}, T-t_k,\cdot)$ are represented as $f_{I,i,j,k}$ and $\hat L_{i,j,k}$, respectively. Our objective is to compute the values of these functions at the grid points, specifically $\varphi_{i,j,k}$, $\rho_{i,j,k}$, and $\alpha_{IJ,i,j,k}$ for $I,J=1,2$.
Following this, the formulation of the saddle point challenge, normalized by $\Delta t$, is presented as:
\begin{equation*}
\begin{adjustbox}{width=0.99\textwidth}$
\begin{split}
\min_{\substack{\varphi\\ \varphi_{i,j,1}=g(x_{i,j})}}\max_{\rho\geq 0, \alpha}  \sum_{k=2}^{n_t} \sum_{i=1}^{n_x} \sum_{j=1}^{n_y} \rho_{i,j,k}\Bigg((D_t^-\varphi)_{i,j,k} - f_{1,i,j,k}(\alpha_{11,i,j,k})_+ (D_x^+\varphi)_{i,j,k} 
- f_{1,i,j,k}(\alpha_{12,i,j,k})_- (D_x^-\varphi)_{i,j,k} \\
- f_{2,i,j,k}(\alpha_{21,i,j,k})_+ (D_y^+\varphi)_{i,j,k} - f_{2,i,j,k}(\alpha_{22,i,j,k})_- (D_y^-\varphi)_{i,j,k} \\
- \hat L_{i,j,k}\left(\alpha_{11,i,j,k}, \alpha_{12,i,j,k}, \alpha_{21,i,j,k}, \alpha_{22,i,j,k}\right) \Bigg) - \frac{c}{\Delta t}\sum_{i=1}^{n_x} \sum_{j=1}^{n_y} \varphi_{i,j,n_t}.
\end{split}
$\end{adjustbox}
\end{equation*}
The update process for the $\ell$-th iteration in the proposed algorithm is then established as
\begin{equation*}
\begin{adjustbox}{width=0.99\textwidth}$
\begin{dcases}
\rho_{i,j,k}^{\ell+1} = \Big(\rho_{i,j,k}^{\ell} + \tau_\rho\big((D_t^-{\tilde\varphi^\ell})_{i,j,k} - f_{1,i,j,k}(\alpha^{\ell}_{11,i,j,k})_+ (D_x^+ \tilde\varphi^\ell)_{i,j,k} - f_{1,i,j,k}(\alpha^\ell_{12,i,j,k})_- (D_x^- \tilde\varphi^\ell)_{i,j,k} \\
\quad\quad\quad\quad- f_{2,i,j,k}(\alpha^{\ell}_{21,i,j,k})_+ (D_y^+ \tilde\varphi^\ell)_{i,j,k} - f_{2,i,j,k}(\alpha^\ell_{22,i,j,k})_- (D_y^- \tilde\varphi^\ell)_{i,j,k}\\
\quad\quad\quad\quad
- \hat L_{i,j,k}(\alpha^\ell_{11,i,j,k}, \alpha^\ell_{12,i,j,k},\alpha^\ell_{21,i,j,k}, \alpha^\ell_{22,i,j,k})\big)\Big)_+.\\
(\alpha^{\ell+1}_{11,i,j,k}, \alpha^{\ell+1}_{12,i,j,k}, \alpha^{\ell+1}_{21,i,j,k}, \alpha^{\ell+1}_{22,i,j,k}) = \argmin_{\alpha_{11}, \alpha_{12}, \alpha_{21}, \alpha_{22}\in\R^m}\{f_{1,i,j,k}(\alpha_{11})_+ (D_x^+\tilde\varphi^\ell)_{i,j,k} \\
\quad\quad\quad\quad
+ f_{1,i,j,k}(\alpha_{12})_- (D_x^-\tilde\varphi^\ell)_{i,j,k} + f_{2,i,j,k}(\alpha_{21})_+ (D_y^+\tilde\varphi^\ell)_{i,j,k} \\
\quad\quad\quad\quad
+ f_{2,i,j,k}(\alpha_{22})_- (D_y^-\tilde\varphi^\ell)_{i,j,k} + \hat L_{i,j,k}(\alpha_{11}, \alpha_{12}, \alpha_{21}, \alpha_{22})\\
\quad\quad\quad\quad
+ \frac{\rho_{i,j,k}^{\ell+1}}{2\tau_\alpha} \left( \|\alpha_{11} - \alpha^\ell_{11,i,j,k}\|^2 + \|\alpha_{12} - \alpha^\ell_{12,i,j,k}\|^2 +  \|\alpha_{21} - \alpha^\ell_{21,i,j,k}\|^2 +  \|\alpha_{22} - \alpha^\ell_{22,i,j,k}\|^2\right)\}.\\
\varphi^{\ell+1}_{i,j,k} = \varphi^{\ell}_{i,j,k} + \tau_{\varphi}(I-D_{tt}-D_{xx} - D_{yy})^{-1}\Big((D_t^+{\rho}^{\ell+1})_{i,j,k} - D_x^-\left(f_{1,i,j,k}(\alpha^{k+1}_{11,i,j,k})_+\rho^{\ell+1}_{i,j,k}\right) \\
\quad\quad\quad\quad - D_x^+\left(f_{1,i,j,k}(\alpha^{k+1}_{12,i,j,k})_-\rho^{\ell+1}_{i,j,k}\right) - D_y^-\left(f_{2,i,j,k}(\alpha^{k+1}_{21,i,j,k})_+\rho^{\ell+1}_{i,j,k}\right) \\
\quad\quad\quad\quad - D_y^+\left(f_{2,i,j,k}(\alpha^{k+1}_{22,i,j,k})_-\rho^{\ell+1}_{i,j,k}\right) \Big).\\
\tilde\varphi^{\ell+1} = 2\varphi^{\ell+1} - \varphi^\ell.
\end{dcases}
$\end{adjustbox}
\end{equation*}
It is important to note that selecting the numerical Lagrangian $\hat L$ as $L_{x,t}(\alpha_{11}) + L_{x,t}(\alpha_{12}) + L_{x,t}(\alpha_{21}) + L_{x,t}(\alpha_{22})$ enables parallel computation of the update steps for each of the variables $\alpha_{11}$, $\alpha_{12}$, $\alpha_{21}$, and $\alpha_{22}$.

\subsection{Numerical examples}\label{sec:numerics_det}
This section focuses on the optimal control problem as described in~\eqref{eqt:oc_problem}, characterized by the dynamics 
\begin{equation*}
    f(x,t,\alpha) = A(x,t)\alpha + b(x,t),
\end{equation*}
where $A(x,t)$ represents an $n\times m$ matrix, and $b(x,t)$ is a vector in $\Rn$. It is further assumed that the Lagrangian $L$ exhibits convexity with respect to $\alpha$. 
The corresponding Hamiltonian $H$ in~\eqref{eqt:cont_HJ_initial} is detailed as follows
\begin{equation*}
\begin{split}
H(x,t,p) &= \sup_{\alpha\in\Rn} \{-\langle A(x,T-t)\alpha + b(x,T-t), p\rangle - L(x,T-t,\alpha)\}\\
&= -\langle b(x,T-t),p\rangle +\sup_{\alpha\in\Rn} \{\langle \alpha, -A(x,T-t)^Tp\rangle - L(x,T-t,\alpha)\}\\
&= -\langle b(x,T-t),p\rangle + L^*(x, T-t, -A(x,T-t)^Tp),
\end{split}
\end{equation*}
with $L^*(x,T-t,\cdot)$ denoting the Fenchel-Legendre transform of $L(x,T-t,\cdot)$ for any $x\in \Omega$ and $t\in [0,T]$.

In the cases where $n=1$, the updates of $\alpha$ in~\eqref{eqt:pdhg_det_semi} are computed by
\begin{equation}\label{eqt:numerics_det_alpha_update}
\begin{adjustbox}{width=0.99\textwidth}$
\begin{split}
\alpha^{\ell+1}_{1,i}(t)
& = \argmin_{f_{i,t}(\alpha)\geq 0}\left\{(A(x,T-t)\alpha + b(x,T-t)) (D_x^+\tilde\varphi^\ell)_i(t) + L_{i,t}(\alpha) + \frac{\rho_i^{\ell+1}(t)}{2\tau_\alpha} \|\alpha - \alpha^\ell_{1,i}(t)\|^2\right\}\\
& = \argmin_{f_{i,t}(\alpha)\geq 0}\left\{L_{i,t}(\alpha) + \frac{\rho_i^{\ell+1}(t)}{2\tau_\alpha} \left\|\alpha - \alpha^\ell_{1,i}(t) + \frac{\tau_\alpha (D_x^+\tilde\varphi^\ell)_i(t)}{\rho_i^{\ell+1}(t)} A(x,T-t)^T\right\|^2\right\}, \\
\alpha^{\ell+1}_{2,i}(t) 
&= \argmin_{f_{i,t}(\alpha)\leq 0}\left\{(A(x,T-t)\alpha + b(x,T-t)) (D_x^-\tilde\varphi^\ell)_i(t) + L_{i,t}(\alpha) + \frac{\rho_i^{\ell+1}(t)}{2\tau_\alpha} \|\alpha - \alpha^\ell_{2,i}(t)\|^2\right\}\\
&= \argmin_{f_{i,t}(\alpha)\leq 0}\left\{L_{i,t}(\alpha) + \frac{\rho_i^{\ell+1}(t)}{2\tau_\alpha} \left\|\alpha - \alpha^\ell_{2,i}(t) + \frac{\tau_\alpha (D_x^-\tilde\varphi^\ell)_i(t)}{\rho_i^{\ell+1}(t)} A(x,T-t)^T  \right\|^2\right\},
\end{split}
$\end{adjustbox}
\end{equation}
which are $m$-dimensional proximal point problems of the functions $L_{x,t}$ and can be computed in parallel for different $x,t$. The two-dimensional $\alpha$-updates are similar as the one-dimensional case.

In each example, we begin by numerically solving for \(\varphi\), \(\alpha_1\), and \(\alpha_2\) in the one-dimensional case, or $\varphi$, \(\alpha_{11}\), \(\alpha_{12}\), \(\alpha_{21}\), and \(\alpha_{22}\) in the two-dimensional case. Following this, applying the forward Euler method to~\eqref{eqt:ode_openloop}, we obtain several optimal trajectories \(\gamma^*\) and the optimal open-loop controls \(\alpha^*\), which provide the solutions to~\eqref{eqt:oc_problem} with \(t=0\) and different initial conditions $x$. For the purpose of demonstration, in one-dimensional examples, we illustrate the level sets of \(\varphi\) and \(\alpha_1 + \alpha_2\) within the spatial-temporal domain. We also plot the graphs of the functions \(\gamma^*\) and \(\alpha^*\). In the case of two-dimensional examples, the level sets of \(\varphi(\cdot, t)\) and each component of \(\alpha_{11}(\cdot, t) + \alpha_{12}(\cdot, t) + \alpha_{21}(\cdot, t) + \alpha_{22}(\cdot, t)\) are depicted in the spatial domain for several times \(t \in [0,T]\). We display the paths of \(\gamma^*\) and \(\alpha^*\) by plotting their values in the spatial domain.

\subsubsection{Quadratic Hamiltonian with spatial dependent coefficients}\label{sec:eg1}
\begin{figure}[htbp]
    \centering
    \begin{subfigure}{0.45\textwidth}
        \centering \includegraphics[width=\textwidth]{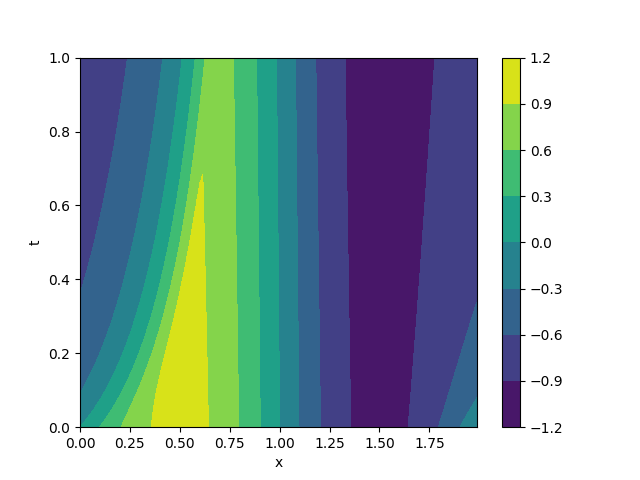}
        \caption{$\varphi$ in the $xt$-space}
    \end{subfigure}
    \hfill
    \begin{subfigure}{0.45\textwidth}
        \centering \includegraphics[width=\textwidth]{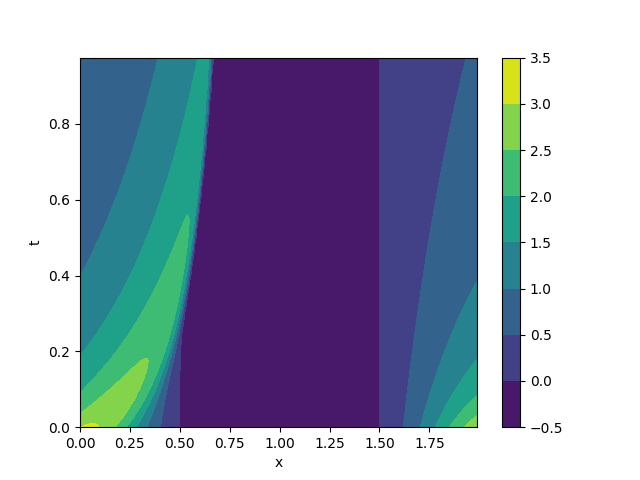}
        \caption{$\alpha$ in the $xt$-space}
    \end{subfigure}\\
    \begin{subfigure}{0.45\textwidth}
        \centering \includegraphics[width=\textwidth]{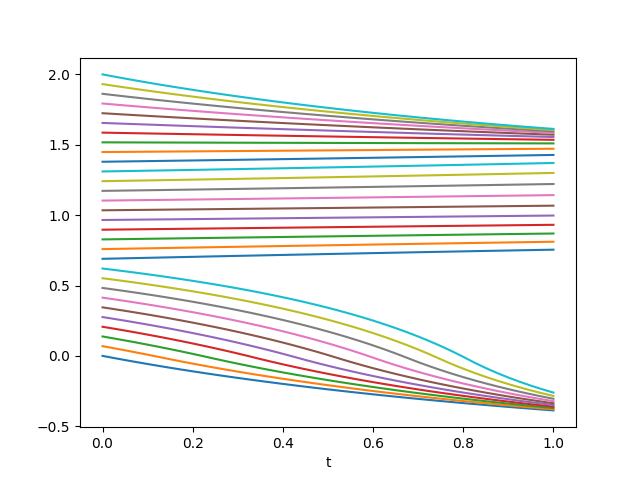}
        \caption{Optimal trajectories $s\mapsto \gamma^*(s)$}
    \end{subfigure}
    \hfill
    \begin{subfigure}{0.45\textwidth}
        \centering \includegraphics[width=\textwidth]{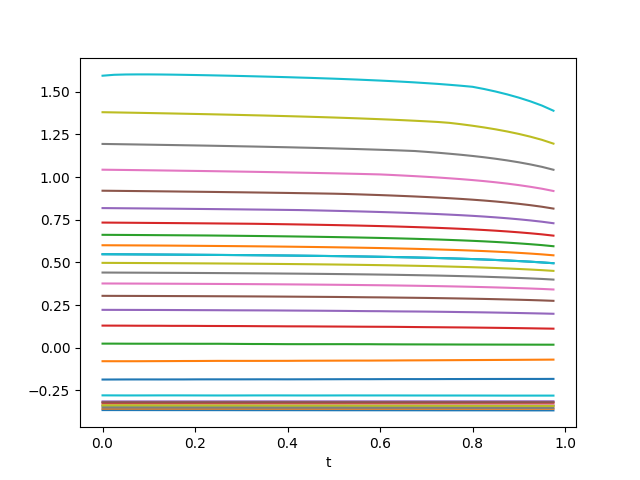}
        \caption{Optimal controls $s\mapsto \alpha^*(s)$}
    \end{subfigure}
    
    \caption{
    Visualization of the solution for the one-dimensional scenario discussed in Section~\ref{sec:eg1}, using $n_t = 41$ and $n_x = 160$ grid points. Figures (a) and (b) showcase the level sets of the solution $\varphi$ to the HJ PDE~\eqref{eqt:cont_HJ_initial}, along with the corresponding function $\alpha$ from~\eqref{eqt:feedback_initial}, which represents the time reversal of the feedback control function. Figures (c) and (d) depict several optimal paths $s\mapsto \gamma^*(s)$ and their associated open-loop optimal controls $s\mapsto \alpha^*(s)$. These paths and control trajectories are the solutions to the optimal control problem~\eqref{eqt:oc_problem}, each beginning from a unique initial condition $x$.} 
    \label{fig:eg1_1d}
\end{figure}

\begin{figure}[htbp]
    \centering
    \begin{subfigure}{0.45\textwidth}
        \centering \includegraphics[width=\textwidth]{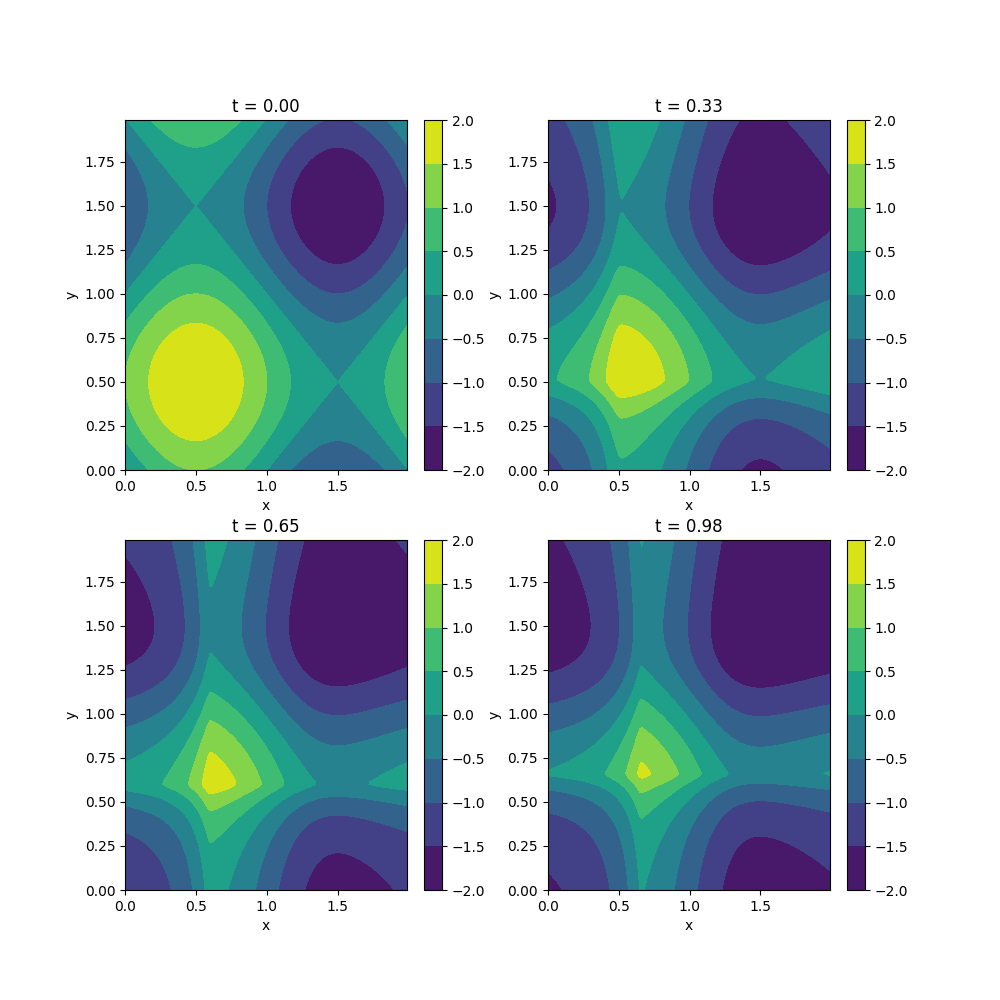}
        \caption{$(x,y)\mapsto \varphi(x,y,t)$ at different $t$}
    \end{subfigure}
    \\
    \begin{subfigure}{0.45\textwidth}
        \centering \includegraphics[width=\textwidth]{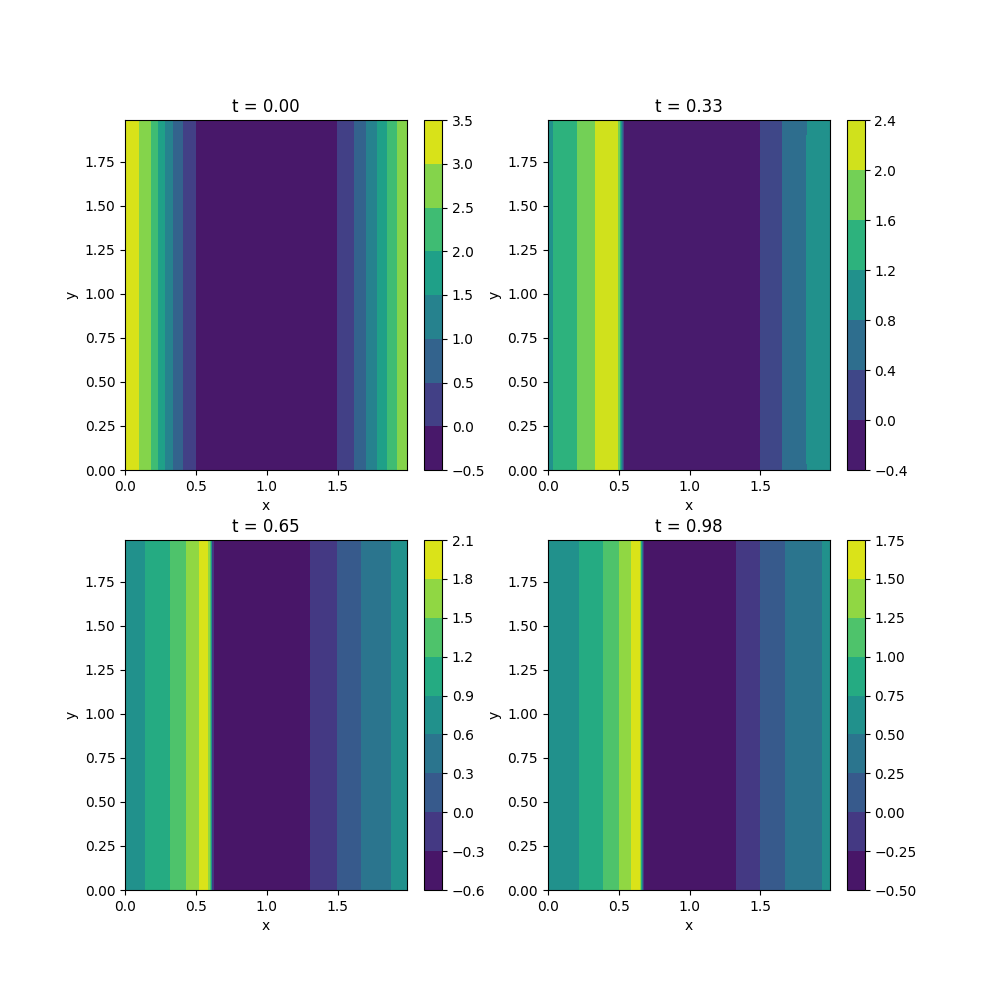}
        \caption{The first component of $(x,y)\mapsto \alpha(x,y,t)$ at different $t$}
    \end{subfigure}\hfill
    \begin{subfigure}{0.45\textwidth}
        \centering \includegraphics[width=\textwidth]{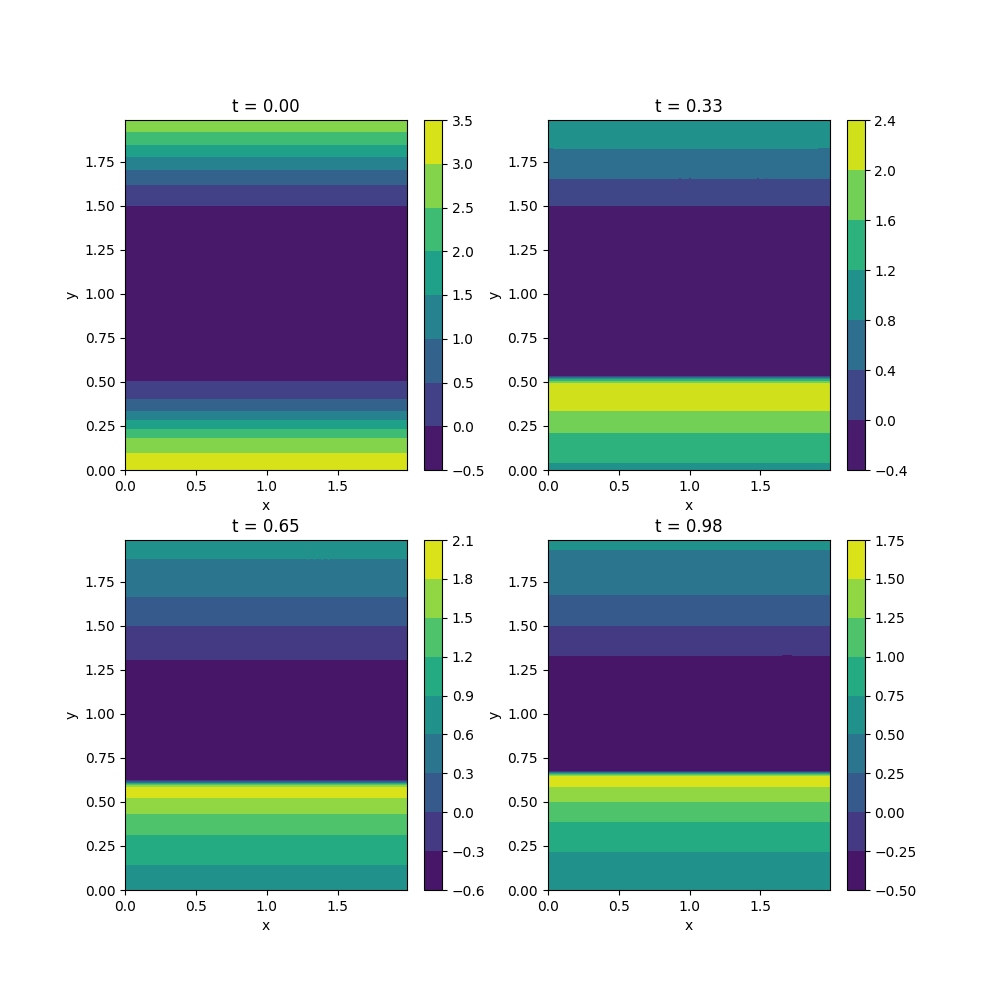}
        \caption{The second component of $(x,y)\mapsto \alpha(x,y,t)$ at different $t$}
    \end{subfigure}\\
    \begin{subfigure}{0.45\textwidth}
        \centering \includegraphics[width=\textwidth]{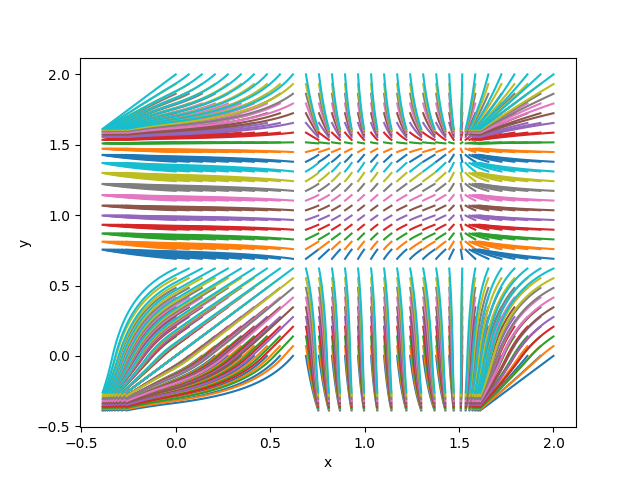}
        \caption{Optimal trajectories $\gamma^*$}
    \end{subfigure}
    \hfill
    \begin{subfigure}{0.45\textwidth}
        \centering \includegraphics[width=\textwidth]{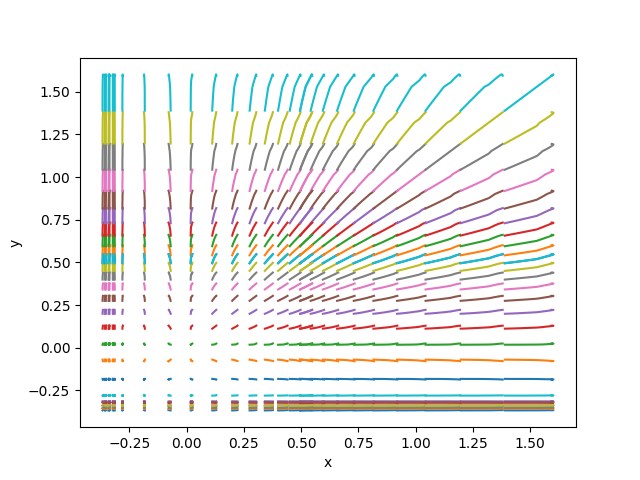}
        \caption{Optimal controls $\alpha^*$}
    \end{subfigure}
    
    \caption{Depiction of the two-dimensional solution as discussed in Section~\ref{sec:eg1}, utilizing $n_t = 41$ and $n_x = n_y = 160$ grid points. Figure (a) illustrates the level sets of the solution $\varphi(\cdot, t)$ to the HJ PDE~\eqref{eqt:cont_HJ_initial} at different times $t$. Figures (b) and (c) show the first and second components, respectively, of the associated function $\alpha(\cdot, t)$ from~\eqref{eqt:feedback_initial} at various times $t$, which depict the time reversal of the feedback control function. Figures (d) and (e) present several optimal trajectories $s\mapsto \gamma^*(s)$ along with their corresponding open-loop optimal controls $s\mapsto \alpha^*(s)$. These trajectories and control strategies solve the optimal control problem specified in~\eqref{eqt:oc_problem}, starting from distinct initial conditions $x$. Notably, both $\gamma^*$ and $\alpha^*$ take values in $\R^2$. For visualization, they are plotted within the spatial domain, excluding the time dimension for clarity.} 
    \label{fig:eg1_2d}
\end{figure}

In this section, we explore an optimal control problem characterized by dynamics 
$f(x,t,\alpha)_I = -(|x_I-1|^2 + 0.1)\alpha_I$
for $I=1,\dots, n$, with a Lagrangian $L(x,t,\alpha) = \frac{1}{2}|\alpha|^2$. Accordingly, the Hamiltonian defined in~\eqref{eqt:cont_HJ_initial} is expressed as $H(x,t,p) = \frac{1}{2}\sum_{I=1}^n(|x_I-1|^2 + 0.1)^2p_I^2$. The terminal cost function is $g(x) = \sum_{I=1}^n\sin \pi x_I$, over the spatial domain $[0,2]^n$ under periodic boundary conditions. The control dimension $m$ is equal to the state dimension $n$.

We adopt the numerical Lagrangian $\hat L$ as $\hat L_{x,t}(\alpha_1,\alpha_2) = L_{x,t}(\alpha_1) + L_{x,t}(\alpha_2)$ in one-dimensional cases and $\hat L_{x,t}(\alpha_{11},\alpha_{12}, \alpha_{21},\alpha_{22}) = L_{x,t}(\alpha_{11}) + L_{x,t}(\alpha_{12}) + L_{x,t}(\alpha_{21}) + L_{x,t}(\alpha_{22})$ for two-dimensional scenarios. As discussed in Section~\ref{sec:appendix_consistency}, these choices for the numerical Lagrangians yield consistent Hamiltonians in this example.
For the one-dimensional scenario, the $\alpha$ updates as outlined in~\eqref{eqt:numerics_det_alpha_update} are calculated by:
\begin{equation}\label{eqt:numerics_det_eg1_alpha}
\begin{adjustbox}{width=0.99\textwidth}$
\begin{split}
\alpha^{\ell+1}_{1,i}(t)
& = \argmin_{\alpha\leq 0}\left\{\frac{1}{2}\alpha^2 + \frac{\rho_i^{\ell+1}(t)}{2\tau_\alpha} \left(\alpha - \alpha^\ell_{1,i}(t) - \frac{\tau_\alpha}{\rho_i^{\ell+1}(t)} (|x_i-1|^2 + 0.1)(D_x^+\tilde\varphi^\ell)_i(t)\right)^2\right\}, \\
&= \min\left\{0, \frac{\rho_i^{\ell+1}(t)\alpha^\ell_{1,i}(t) + \tau_\alpha (|x_i-1|^2 + 0.1)(D_x^+\tilde\varphi^\ell)_i(t)}{\rho_i^{\ell+1}(t) + \tau_\alpha}\right\},
\\
\alpha^{\ell+1}_{2,i}(t) 
&= \argmin_{\alpha\geq 0}\left\{\frac{1}{2}\alpha^2 + \frac{\rho_i^{\ell+1}(t)}{2\tau_\alpha} \left(\alpha - \alpha^\ell_{2,i}(t) - \frac{\tau_\alpha}{\rho_i^{\ell+1}(t)} (|x_i-1|^2 + 0.1) (D_x^-\tilde\varphi^\ell)_i(t) \right)^2\right\}
\\
&= \max\left\{0, \frac{\rho_i^{\ell+1}(t)\alpha^\ell_{2,i}(t) + \tau_\alpha (|x_i-1|^2 + 0.1)(D_x^-\tilde\varphi^\ell)_i(t)}{\rho_i^{\ell+1}(t) + \tau_\alpha}\right\}.
\end{split}
$\end{adjustbox}
\end{equation}
Two-dimensional scenarios follow a similar computation process.
In our numerical experiment, we selected $n_t = 41$ and $n_x = n_y = 160$ as the number of grid points. 
Utilizing the implicit Euler method for time discretization allows us to circumvent the CFL condition, enabling the adoption of a larger time step $\Delta t$.

The solution of the one-dimensional problem ($n=m=1$) is illustrated in Figure~\ref{fig:eg1_1d}, showcasing the level sets of our numerical solutions $\varphi$ and $\alpha_1+\alpha_2$ in (a) and (b), respectively. Optimal trajectories and controls for evenly spaced initial conditions within $\Omega$ are depicted in (c) and (d).
The solution for the two-dimensional problem ($n=m=2$) follows a similar methodology and is depicted in Figure~\ref{fig:eg1_2d}.

\subsubsection{One-homogeneous Hamiltonian with spatial dependent coefficients} \label{sec:eg2}

\begin{figure}[htbp]
    \centering
    \begin{subfigure}{0.45\textwidth}
        \centering \includegraphics[width=\textwidth]{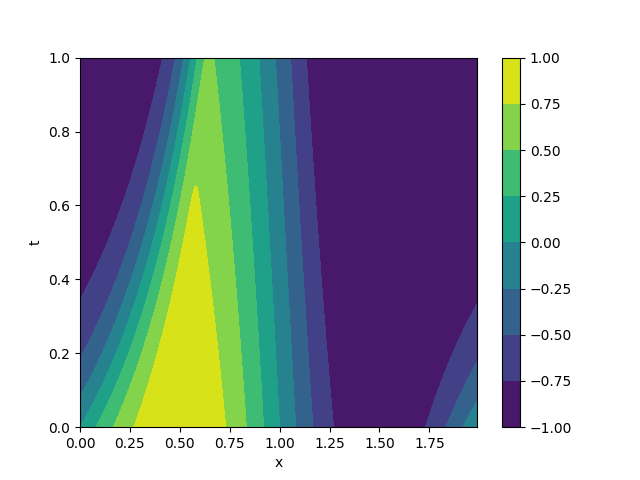}
        \caption{$\varphi$ in the $xt$-space}
    \end{subfigure}
    \hfill
    \begin{subfigure}{0.45\textwidth}
        \centering \includegraphics[width=\textwidth]{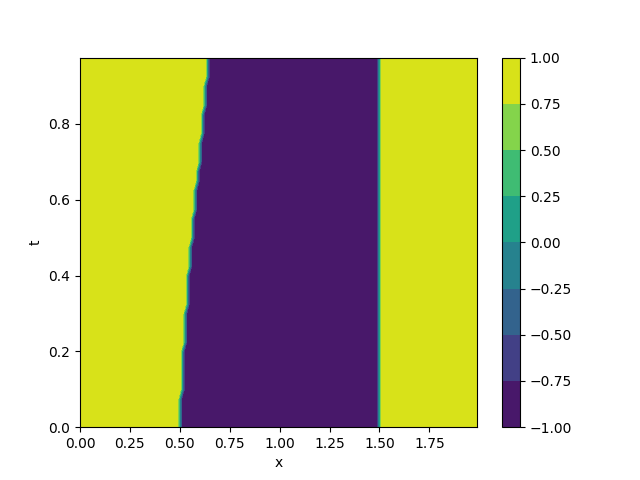}
        \caption{$\alpha$ in the $xt$-space}
    \end{subfigure}\\
    \begin{subfigure}{0.45\textwidth}
        \centering \includegraphics[width=\textwidth]{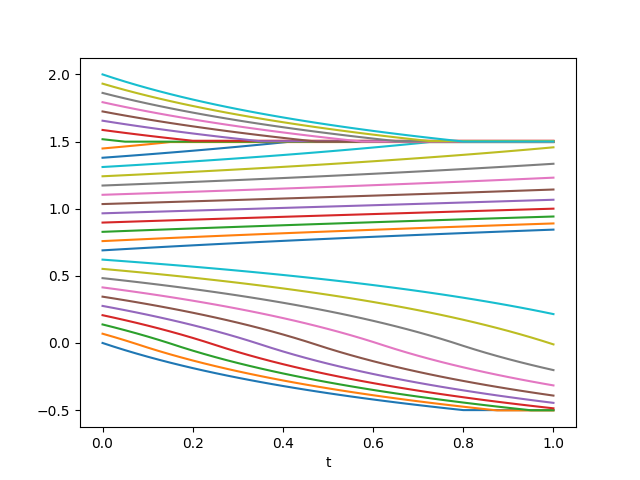}
        \caption{Optimal trajectories $s\mapsto \gamma^*(s)$}
    \end{subfigure}
    \hfill
    \begin{subfigure}{0.45\textwidth}
        \centering \includegraphics[width=\textwidth]{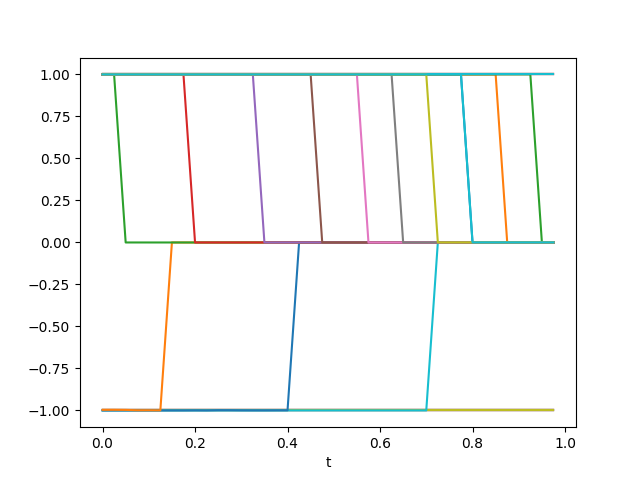}
        \caption{Optimal controls $s\mapsto \alpha^*(s)$}
    \end{subfigure}
    
    \caption{Visualization of the solution for the one-dimensional scenario discussed in Section~\ref{sec:eg2}, using $n_t = 41$ and $n_x = 160$ grid points. Figures (a) and (b) showcase the level sets of the solution $\varphi$ to the HJ PDE~\eqref{eqt:cont_HJ_initial}, along with the corresponding function $\alpha$ from~\eqref{eqt:feedback_initial}, which represents the time reversal of the feedback control function. Figures (c) and (d) depict several optimal paths $s\mapsto \gamma^*(s)$ and their associated open-loop optimal controls $s\mapsto \alpha^*(s)$. These paths and control trajectories are the solutions to the optimal control problem~\eqref{eqt:oc_problem}, each beginning from a unique initial condition $x$.} 
    \label{fig:eg2_1d}
\end{figure}

\begin{figure}[htbp]
    \centering
    \begin{subfigure}{0.45\textwidth}
        \centering \includegraphics[width=\textwidth]{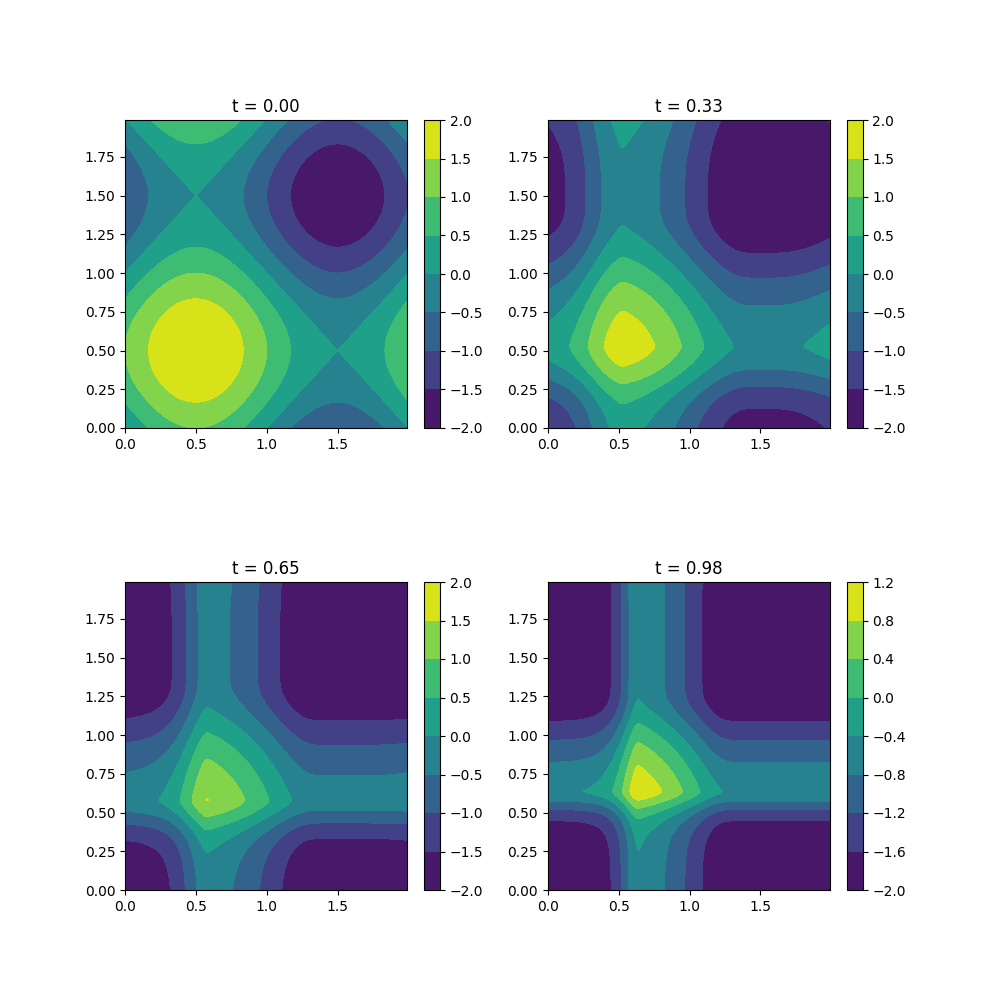}
        \caption{$(x,y)\mapsto \varphi(x,y,t)$ at different $t$}
    \end{subfigure}
    \\
    \begin{subfigure}{0.45\textwidth}
        \centering \includegraphics[width=\textwidth]{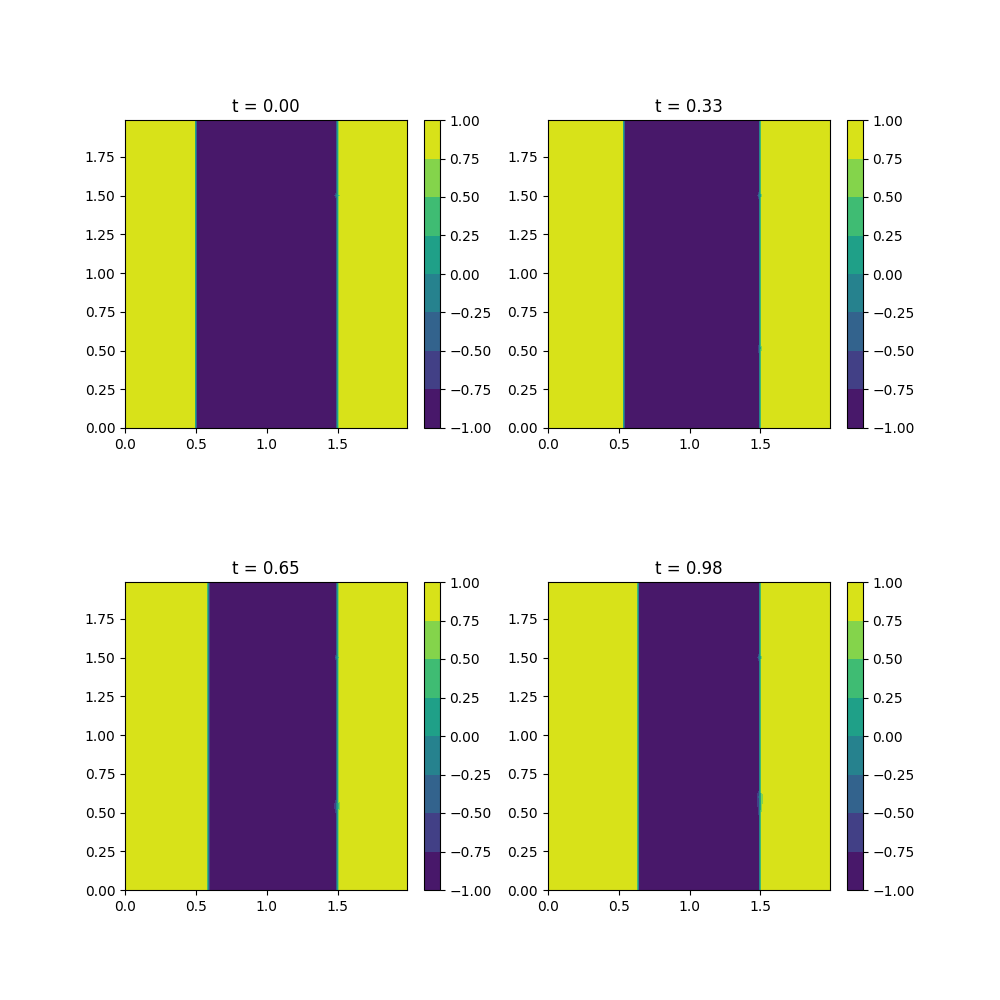}
        \caption{The first component of $(x,y)\mapsto \alpha(x,y,t)$ at different $t$}
    \end{subfigure}\hfill
    \begin{subfigure}{0.45\textwidth}
        \centering \includegraphics[width=\textwidth]{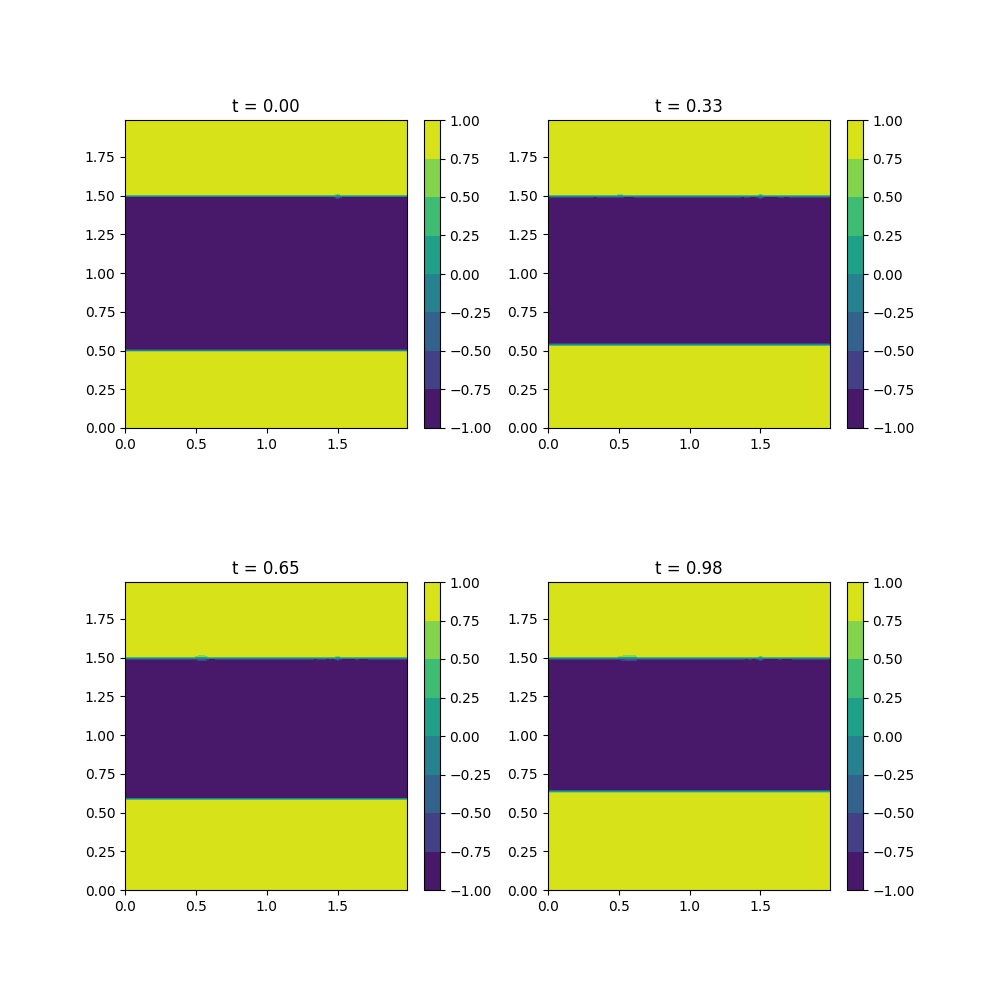}
        \caption{The second component of $(x,y)\mapsto \alpha(x,y,t)$ at different $t$}
    \end{subfigure}\\
    \begin{subfigure}{0.45\textwidth}
        \centering \includegraphics[width=\textwidth]{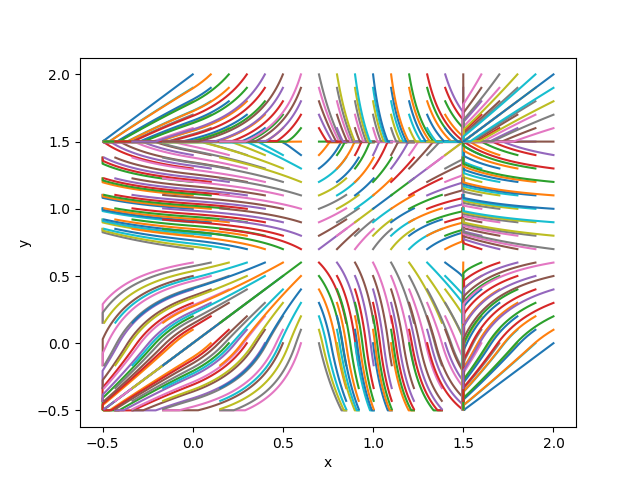}
        \caption{Optimal trajectories $\gamma^*$}
    \end{subfigure}
    \hfill
    \begin{subfigure}{0.45\textwidth}
        \centering \includegraphics[width=\textwidth]{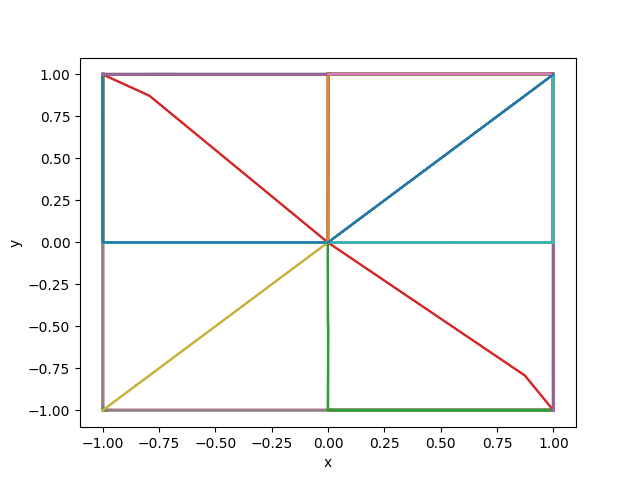}
        \caption{Optimal controls $\alpha^*$}
    \end{subfigure}
    
    \caption{
    Depiction of the two-dimensional solution as discussed in Section~\ref{sec:eg2}, utilizing $n_t = 41$ and $n_x = n_y = 160$ grid points. Figure (a) illustrates the level sets of the solution $\varphi(\cdot, t)$ to the HJ PDE~\eqref{eqt:cont_HJ_initial} at different times $t$. Figures (b) and (c) show the first and second components, respectively, of the associated function $\alpha(\cdot, t)$ from~\eqref{eqt:feedback_initial} at various times $t$, which depict the time reversal of the feedback control function. Figures (d) and (e) present several optimal trajectories $s\mapsto \gamma^*(s)$ along with their corresponding open-loop optimal controls $s\mapsto \alpha^*(s)$. These trajectories and control strategies solve the optimal control problem specified in~\eqref{eqt:oc_problem}, starting from distinct initial conditions $x$. Notably, both $\gamma^*$ and $\alpha^*$ take values in $\R^2$. For visualization, they are plotted within the spatial domain, excluding the time dimension for clarity.} 
    \label{fig:eg2_2d}
\end{figure}

We maintain the same spatial domain $\Omega$ and utilize the functions $f$ and $g$ as delineated in Section~\ref{sec:eg1}. The Lagrangian is set to $L(x,t,\alpha) = \ind_{B_\infty}(\alpha)$, where $\ind$ signifies the indicator function, and $B_\infty$ represents the unit ball in the $\ell^\infty$ space, defined as $B_\infty = \{\alpha\in \R^m : \|\alpha_I\|_\infty \leq 1\}$. Recall that the indicator function $\ind_C(x)$ associated with a set $C$ is set to $0$ when $x$ falls within $C$ and is assigned $+\infty$ for values of $x$ outside $C$. This choice of Lagrangian imposes a constraint on the control values. The Hamiltonian described in~\eqref{eqt:cont_HJ_initial} becomes $H(x,t,p) = \sum_{I=1}^n(|x_I-1|^2 + 0.1) |p_I|$, which is $1$-homogeneous with respect to $p$.
We select the numerical Lagrangian $\hat L$ as $\hat L_{x,t}(\alpha_1,\alpha_2) = L_{x,t}(\alpha_1) + L_{x,t}(\alpha_2)$ for the one-dimensional scenarios and $\hat L_{x,t}(\alpha_{11},\alpha_{12}, \alpha_{21},\alpha_{22}) = L_{x,t}(\alpha_{11}) + L_{x,t}(\alpha_{12}) + L_{x,t}(\alpha_{21}) + L_{x,t}(\alpha_{22})$ for two-dimensional cases. As discussed in Section~\ref{sec:appendix_consistency}, these selections ensure consistency in the resulting Hamiltonians.

For one-dimensional scenarios, the $\alpha$ updates detailed in~\eqref{eqt:numerics_det_alpha_update} are computed as follows:
\begin{equation*}
\begin{adjustbox}{width=0.99\textwidth}$
\begin{split}
\alpha^{\ell+1}_{1,i}(t)
& = \argmin_{\alpha\in [-1,0]}\left\{\frac{\rho_i^{\ell+1}(t)}{2\tau_\alpha} \left(\alpha - \alpha^\ell_{1,i}(t) - \frac{\tau_\alpha}{\rho_i^{\ell+1}(t)} (|x_i-1|^2 + 0.1)(D_x^+\tilde\varphi^\ell)_i(t)\right)^2\right\}, \\
&= \max\left\{ -1, \min\left\{0, \alpha^\ell_{1,i}(t) + \frac{\tau_\alpha}{\rho_i^{\ell+1}(t)} (|x_i-1|^2 + 0.1)(D_x^+\tilde\varphi^\ell)_i(t)\right\}\right\},
\\
\alpha^{\ell+1}_{2,i}(t) 
&= \argmin_{\alpha\in [0,1]}\left\{\frac{\rho_i^{\ell+1}(t)}{2\tau_\alpha} \left(\alpha - \alpha^\ell_{2,i}(t) - \frac{\tau_\alpha}{\rho_i^{\ell+1}(t)} (|x_i-1|^2 + 0.1) (D_x^-\tilde\varphi^\ell)_i(t) \right)^2\right\}
\\
&= \min\left\{1, \max\left\{0, \alpha^\ell_{2,i}(t) + \frac{\tau_\alpha}{\rho_i^{\ell+1}(t)} (|x_i-1|^2 + 0.1) (D_x^-\tilde\varphi^\ell)_i(t) \right\}\right\}.
\end{split}
$\end{adjustbox}
\end{equation*}
The methodology for two-dimensional scenarios follows a similar approach. As with the previous example, we opted for $n_t = 41$ and $n_x = n_y = 160$ for the grid points. The implicit Euler method for time discretization demonstrates its benefit by allowing for a larger time step, $\Delta t$.

The solution for the one-dimensional case ($n=m=1$) is depicted in Figure~\ref{fig:eg2_1d}. Initially, we employ the proposed method to calculate $\varphi$, $\alpha_1 + \alpha_2$ at the grid points, as illustrated in Figures~\ref{fig:eg2_1d} (a) and (b). Subsequently, we determine the optimal trajectories $\gamma^*$ and optimal open-loop controls $\alpha^*$ in~\eqref{eqt:oc_problem} for $t=0$ with varying initial conditions $x$ uniformly distributed in $\Omega$. These are showcased in Figures~\ref{fig:eg2_1d} (c) and (d).
The process for solving the two-dimensional problem ($n=m=2$) mirrors the one-dimensional case and is exhibited in Figure~\ref{fig:eg2_2d}.

In contrast to the numerical results presented in Section~\ref{sec:eg1}, the control function $\alpha$ in this scenario lacks smoothness, predominantly adopting values of $-1$, $0$, or $1$. This phenomenon is evidenced by the discontinuities observed in Figures~\ref{fig:eg2_1d} (b) and~\ref{fig:eg2_2d} (b)-(c), as well as in the depicted paths of the optimal controls in Figures~\ref{fig:eg2_1d}~(d) and~\ref{fig:eg2_2d} (e).
The emergence of this distinct bang-bang control pattern is attributable to the control function's formulation as provided in~\eqref{eqt:feedback_initial}. Specifically, in the context of this example, the control function is defined as
\begin{equation*}
\alpha(x,t) \in \argmin_{\alpha\in B_\infty} \sum_{I=1}^n -(|x_I-1|^2 + 0.1) \partial_{x_I}\varphi(x,t) \alpha_I.
\end{equation*}
Consequently, when the condition $(|x_I-1|^2 + 0.1) \partial_{x_I}\varphi(x,t) > 1$ is satisfied, it results in $\alpha_I(x,t) = 1$. Conversely, if $(|x_I-1|^2 + 0.1) \partial_{x_I}\varphi(x,t) < -1$, then $\alpha_I(x,t)$ equals $-1$. This mechanism underpins the prevalence of these specific control values in the visual representations.

In this example, the optimal control problem is described as:
\begin{equation*}
    \begin{split}
\min_{\alpha(\cdot)}\left\{g(\gamma(T)) \colon \gamma(t)=x, \dot{\gamma}(s) = f(\gamma(s),s,\alpha(s)), \alpha(s)\in B_\infty, \forall s\in(t,T) \right\}.
\end{split}
\end{equation*}
This setup lacks a running cost and instead imposes a constraint on control. This means that if a feasible control leads to a trajectory that starts at $\gamma(t) = x$ and ends at a point $\gamma (T)$ that minimizes $g$ over $\Omega$, such a control is optimal. Consequently, it's possible to have multiple optimal controls. In our numerical results, we predominantly encounter controls exhibiting a bang-bang behavior, initially adopting values of $1$ or $-1$ and subsequently transitioning to zero. The corresponding trajectories initially move at maximum speed in the desired direction and then halt upon reaching a minimizer of $g$. It is an interesting question to characterize the controls selected by the saddle point problem in~\eqref{eqt:saddle_oc} and through our numerical algorithms.

\subsubsection{Newton mechanics}\label{sec:numerical_newton_det}

\begin{figure}[htbp]
    \centering
    \begin{subfigure}{0.45\textwidth}
        \centering \includegraphics[width=\textwidth]{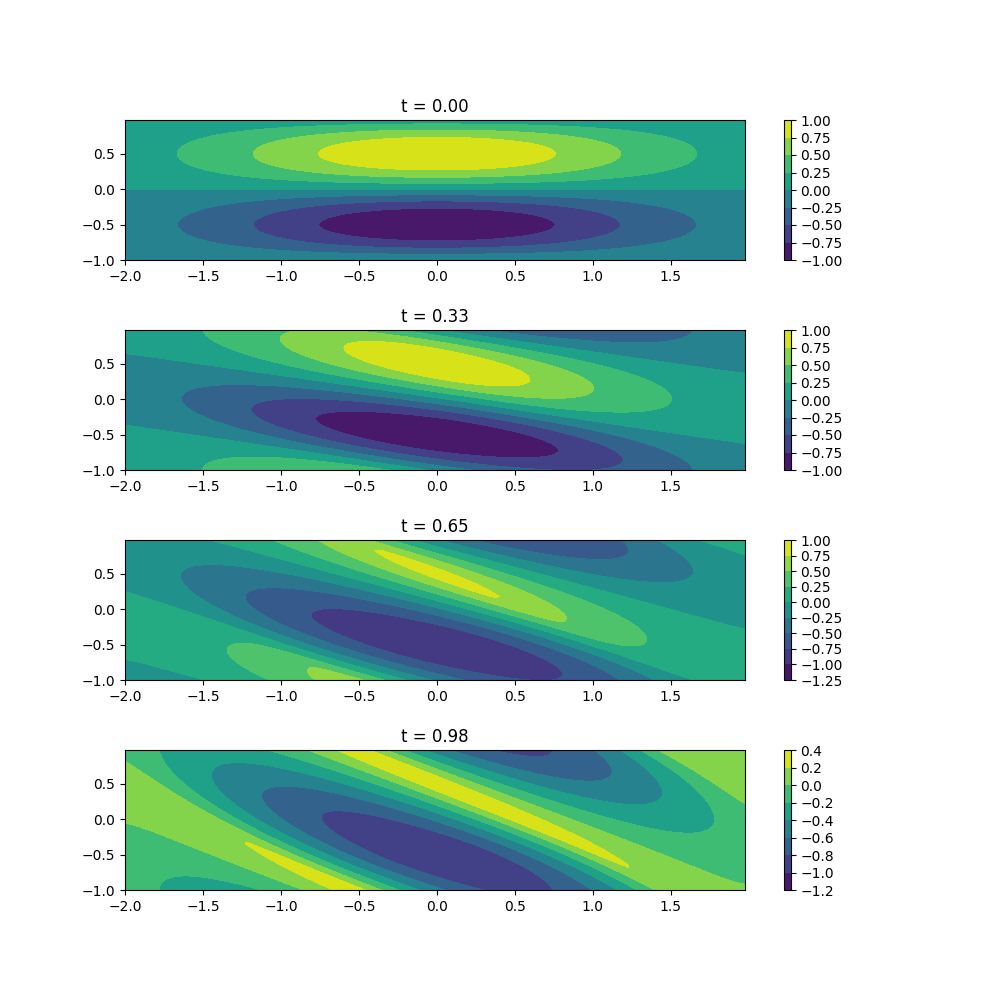}
        \caption{$(x,y)\mapsto \varphi(x,y,t)$ at different $t$}
    \end{subfigure}
    \hfill
    \begin{subfigure}{0.45\textwidth}
        \centering \includegraphics[width=\textwidth]{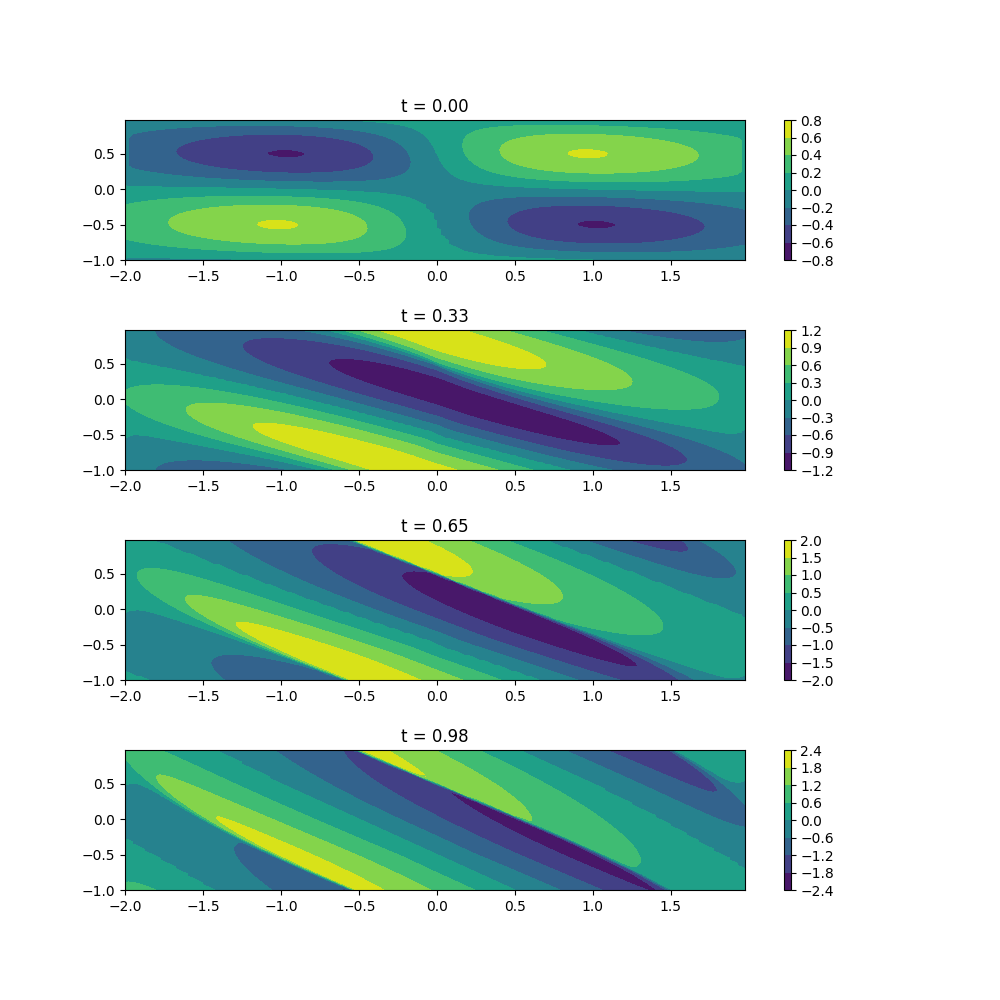}
        \caption{$(x,y)\mapsto \alpha(x,y,t)$ at different $t$}
    \end{subfigure}\\
    \begin{subfigure}{0.45\textwidth}
        \centering \includegraphics[width=\textwidth]{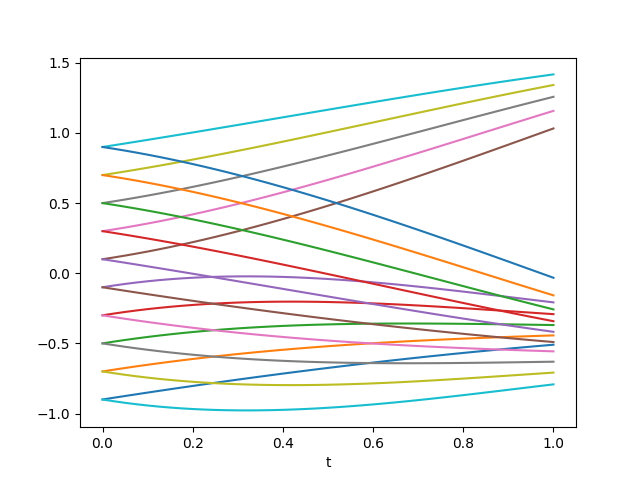}
        \caption{Optimal positions $\gamma^*_2$}
    \end{subfigure}\hfill
    \begin{subfigure}{0.45\textwidth}
        \centering \includegraphics[width=\textwidth]{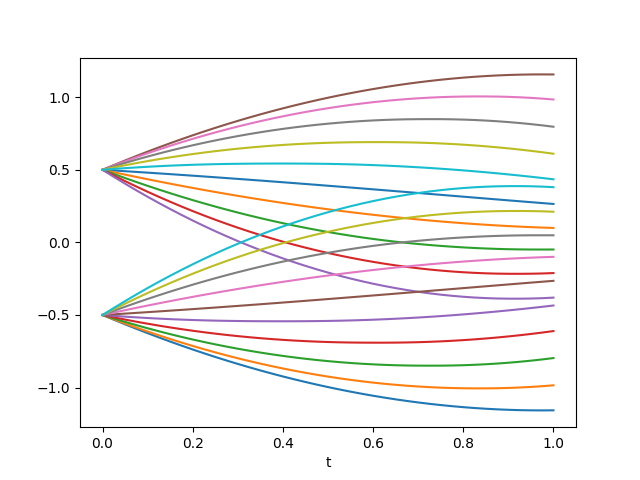}
        \caption{Optimal velocities $\gamma^*_1$}
    \end{subfigure}
    \\
    \begin{subfigure}{0.45\textwidth}
        \centering \includegraphics[width=\textwidth]{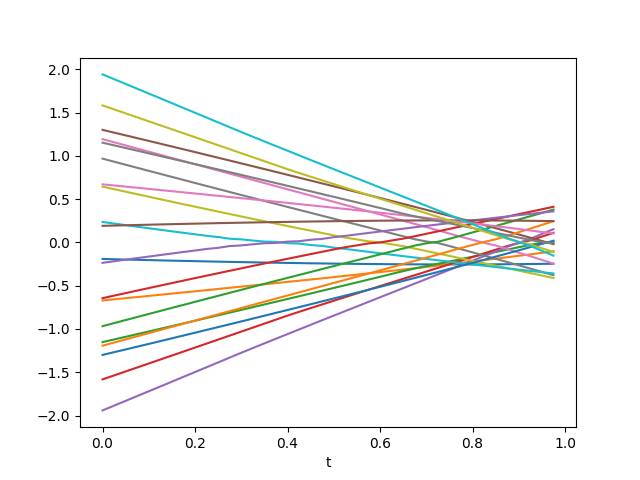}
        \caption{Optimal controls $\alpha^*$}
    \end{subfigure}
    
    \caption{
    Visualization of the solution for example discussed in Section~\ref{sec:numerical_newton_det}, using $n_t = 41$, $n_x = 160$, and $n_y=80$ grid points. Figures (a) and (b) showcase the level sets of the solution $\varphi$ to the HJ PDE~\eqref{eqt:cont_HJ_initial}, along with the corresponding function $\alpha$ from~\eqref{eqt:feedback_initial}, which represents the time reversal of the feedback control function. Figures (c), (d), (e) depict the first and second components of several optimal paths $s\mapsto \gamma^*(s)$ and their associated open-loop optimal controls $s\mapsto \alpha^*(s)$. These paths and control trajectories are the solutions to the optimal control problem~\eqref{eqt:oc_problem}, each beginning from a unique initial condition $x$.
    } 
    \label{fig:eg_newton}
\end{figure}

This example delves into an optimal control problem underpinned by Newtonian mechanics. With the settings $m=1$, $n=2$, and the dynamic equation $f(x_1, x_2, t,\alpha) = [\alpha, x_1]^T$, we derive the ODE constraints as $\dot \gamma_1(s) = \alpha(s)$ and $\dot \gamma_2(s) = \gamma_1(s)$. From a physical standpoint, $\gamma_2$ is interpreted as the position, $\gamma_1$ as the velocity, and the control variable $\alpha$ as the acceleration. The Lagrangian $L$ is defined as the quadratic function $L(x, t, \alpha) = \frac{1}{2}|\alpha|^2$.
The Hamiltonian in~\eqref{eqt:cont_HJ_initial}, is formulated as
\begin{equation*}
H(x,t,p) = \sup_{\alpha\in\R} \{-\langle f(x,\alpha), p\rangle - L(x,t,\alpha) \} =  \frac{1}{2}|p_1|^2 - x_1p_2.
\end{equation*}
The terminal cost function $g$ is expressed as
\begin{equation*}
g(x_1, x_2) = \exp\left(-\frac{x_1^2}{2}\right) \sin(\pi x_2).
\end{equation*}
The domain for the position variable $x_2$ is set to $[-1,1]$ with periodic boundary conditions, whereas the domain for the velocity variable $x_1$ is $[-2,2]$ with Neumann boundary conditions.

Given that the function $f_2$ does not depend on the control variable, we can disregard the variables $\alpha_{21}$ and $\alpha_{22}$ and select $\hat L_{x,t}(\alpha_{11}, \alpha_{12}) = L_{x,t}(\alpha_{11}) + L_{x,t}(\alpha_{12})$.
Thus, the saddle point problem is formulated as
\begin{equation*}
\begin{adjustbox}{width=0.99\textwidth}$
\begin{split}
\min_{\substack{\varphi\\ \varphi_{i,j}(0)=g(x_{i,j})}}\max_{\rho\geq 0, \alpha_{11}, \alpha_{12}}  \int_0^T \sum_{i=1}^{n_x} \sum_{j=1}^{n_y} \rho_{i,j}(t)\Bigg(\dot\varphi_{i,j}(t) - (\alpha_{11,i,j}(t))_+ (D_x^+\varphi)_{i,j}(t) 
- (\alpha_{12,i,j}(t))_- (D_x^-\varphi)_{i,j}(t)\\
- ((x_1)_{i,j})_+ (D_y^+\varphi)_{i,j}(t)
- ((x_1)_{i,j})_- (D_y^-\varphi)_{i,j}(t) 
- \frac{1}{2}\alpha_{11,i,j}(t)^2- \frac{1}{2}\alpha_{12,i,j}(t)^2 \Bigg)dt - c\sum_{i=1}^{n_x}\sum_{j=1}^{n_y} \varphi_{i,j}(T).
\end{split} 
$\end{adjustbox}
\end{equation*}
Updates for $\varphi$ and $\rho$ proceed as outlined in~\eqref{eqt:pdhg_det_semi}, with $\alpha$ updates following the approach described in~\eqref{eqt:numerics_det_alpha_update}. Specifically, we have
\begin{equation*}
\begin{split}
\alpha^{\ell+1}_{11,i,j}(t)
&= \max\left\{0, \frac{\rho_{i,j}^{\ell+1}(t)\alpha^\ell_{11,i,j}(t) - \tau_\alpha (D_x^+\tilde\varphi^\ell)_{i,j}(t)}{\rho_{i,j}^{\ell+1}(t) + \tau_\alpha}\right\},
\\
\alpha^{\ell+1}_{12,i,j}(t) 
&= \min\left\{0, \frac{\rho_{i,j}^{\ell+1}(t)\alpha^\ell_{12,i,j}(t) - \tau_\alpha (D_x^-\tilde\varphi^\ell)_{i,j}(t)}{\rho_{i,j}^{\ell+1}(t) + \tau_\alpha}\right\}.
\end{split}
\end{equation*}
In the numerical experiment, we chose $n_t = 41$, $n_x = 160$, and $n_y = 80$ as the grid specifications, circumventing the CFL condition through the application of the implicit Euler method for temporal discretization.

Figure~\ref{fig:eg_newton} (a) and (b) display the level sets for the solution $\varphi$ of the HJ PDE~\eqref{eqt:cont_HJ_initial}, along with the associated function $\alpha$ described in~\eqref{eqt:feedback_initial}, captured at various times.
For the initial conditions \(x = (x_1, x_2)\) specified in the optimal control problem~\eqref{eqt:oc_problem}, we define the position variable \(x_2\) as points spanning the range \([-1,1]\), and set the velocity \(x_1\) to either \(-0.5\) or \(0.5\).
The calculated optimal controls and trajectories are subsequently illustrated in Figure~\ref{fig:eg_newton} (c), (d), and (e). Since the functions $\gamma_1^*$, $\gamma_2^*$, and $\alpha^*$ take values in $\R$, these three functions are plotted across the spatial-temporal dimensions.

\section{Stochastic optimal control problems}\label{sec:stochastic_oc}
The methods discussed in Section~\ref{sec:oc_problem} are also applicable to stochastic optimal control problems.
In this section, we consider the following stochastic optimal control problems:
\begin{equation}\label{eqt:soc_problem}
\begin{adjustbox}{width=0.99\textwidth}$
\begin{split}
    \min_{\alpha}\left\{\E\left[\int_t^T
    L(X_s,s, \alpha_s) ds + g(X_T)\right] \colon X_t=x, dX_s = f(X_s,s,\alpha_s)ds + \sqrt{2\epsilon} dW_s \right\},
\end{split}
$\end{adjustbox}
\end{equation}
where $W_s$ represents Brownian motion in $\Rn$, the control $\alpha_s$ is an adapted process, and $f$, $L$, and $g$ function as the drift, Lagrangian, and terminal cost, respectively, akin to their definitions in Section~\ref{sec:oc_problem}.
The minimal value of problem~\eqref{eqt:soc_problem} is represented by $\phi(x,t)$, which satisfies a specific viscous HJ PDE as follows:
\begin{equation*}
\begin{adjustbox}{width=0.99\textwidth}$
\begin{dcases}
\frac{\partial \phi(x,t)}{\partial t} + \inf_{\alpha\in \R^m} \{\langle f(x, t, \alpha), \nabla_x \phi(x,t)\rangle + L(x,t,\alpha)\} + \epsilon\Delta_x \phi(x,t) = 0, & x\in \Omega, t\in [0,T],\\
\phi(x,T) = g(x), & x\in \Omega.
\end{dcases}
$\end{adjustbox}
\end{equation*}
Applying a time reversal technique, we arrive at a viscous HJ PDE with an initial condition as:
\begin{equation}\label{eqt:visc_HJ_initial}
\begin{adjustbox}{width=0.99\textwidth}$
\begin{dcases}
\frac{\partial \varphi(x,t)}{\partial t} + \sup_{\alpha\in \R^m} \{-\langle f(x, t, \alpha), \nabla_x \varphi(x,t)\rangle - L(x,t,\alpha)\} = \epsilon\Delta_x \varphi(x,t), & x\in \Omega, t\in [0,T],\\
\varphi(x,0) = g(x), & x\in \Omega.
\end{dcases}
$\end{adjustbox}
\end{equation}
We introduce the function $\alpha\colon \Omega \times [0,T] \to \R^m$ defined as:
\begin{equation} \label{eqt:visc_alp_initial}
\alpha(x,t) = \argmax_{\alpha\in \R^m} \{-\langle f(x, t, \alpha), \nabla_x \varphi(x,t)\rangle - L(x,t,\alpha)\}.
\end{equation}
Upon reversing time, the function $\alpha(x, T-t)$ serves as the feedback control function.
Consequently, the optimal trajectories and controls are determined through the following computations:
\begin{equation} \label{eqt:sde_openloop}
\begin{dcases}
d \gamma^*_s = f(\gamma^*_s, s, \alpha(\gamma^*_s, T-s))ds + \sqrt{2\epsilon}dW_s, & s\in (t,T),\\
\gamma^*_t = x, \\
\alpha^*_s = \alpha(\gamma^*_s, T-s), & s\in [t,T].
\end{dcases}
\end{equation}
For an in-depth exploration of the linkage between stochastic optimal control problems and viscous HJ PDEs, refer to~\cite{Yong1999Stochastic}.

\subsection{Saddle point formulation}
To solve the viscous HJ PDE, we employ methodologies akin to those used for the first-order HJ PDE as outlined in Section~\ref{sec:oc_problem}. This approach involves framing the PDE as a constraint within an optimization problem, which has the objective function $-c\int_{\Omega} \varphi(x,T) dx$, and subsequently introducing a Lagrange multiplier $\rho$. Following computations akin to those for the first-order case yield the subsequent saddle point formula:
\begin{equation}\label{eqt:saddle_soc}
\begin{adjustbox}{width=0.99\textwidth}$
\begin{split}
\min_{\substack{\varphi\\ \varphi(x,0)=g(x)}}\max_{\rho\geq 0, \alpha}  \int_0^T \int_\Omega \rho(x,t)\Bigg(\frac{\partial \varphi(x,t)}{\partial t} -\langle f_{x,t}(\alpha(x,t)), \nabla_x \varphi(x,t)\rangle - L_{x,t}(\alpha(x,t)) \\
- \epsilon\Delta_x \varphi(x,t) \Bigg)dxdt - c\int_\Omega \varphi(x,T)dx.
\end{split}
$\end{adjustbox}
\end{equation}
For a stationary point $(\varphi, \rho, \alpha)$, where $\rho(x,t) > 0$ across all $x\in \Omega$ and $t\in[0,T]$, the first-order optimality conditions are given as follows:
\begin{equation*}
\begin{dcases}
\partial_t \varphi(x,t) + H(x,t, \nabla_x \varphi(x,t)) = \epsilon\Delta_x \varphi(x,t),\\
\partial_t \rho(x,t) + \nabla_x \cdot (\nabla_p H(x,t,\nabla_x \varphi(x,t))\rho(x,t)) + \epsilon\Delta_x \rho(x,t) = 0,\\
\varphi(x,0) = g(x),\quad  \rho(x,T) = c.
\end{dcases}
\end{equation*}

Echoing the insights from Remark~\ref{rem:connection_mfg}, this saddle point problem similarly establishes a connection to an MFG problem. Specifically, with a time reversal, the formulations for $\varphi(x,T-t)$ and $\tilde \rho(x,t) = \rho(x,T-t)$ align with the first-order optimality conditions of the following MFG problem:
\begin{equation*}
\begin{split}
\min_{\tilde\alpha}\Big\{\int_0^T\int_{\Omega}
L(x,s, \tilde\alpha(x,s))\tilde \rho(x,s) dxds + \int_\Omega g(x)\tilde \rho(x,T) dx \colon \quad\quad\quad\quad\\
\partial_t\tilde \rho(x,s) + \nabla_x \cdot (f(x,s,\tilde\alpha(x,s))\tilde\rho(x,s)) = \epsilon\Delta_x \tilde\rho(x,s) \ 
 \forall x\in\Omega, s\in [0,T],\\
 \tilde \rho(x,0) = c \ \forall x\in \Omega\Big\}.
\end{split}
\end{equation*}
Given that the constraint represents a Fokker-Planck equation associated with drifted Brownian motion, the resulting density is inherently positive throughout, thereby naturally fulfilling the $\rho > 0$ assumption.

The process of discretizing the saddle point problem as described in~\eqref{eqt:saddle_soc} follows a similar approach to that outlined in Section~\ref{sec:discretization_det}, with an additional step to incorporate the diffusion term using a second-order centered difference scheme. The specifics of this approach are elaborated in Appendix~\ref{sec:discretization_soc}.

The primary distinction between the formulations in~\eqref{eqt:saddle_soc} and~\eqref{eqt:saddle_oc} lies in the inclusion of the diffusion term. This addition introduces challenges for the convergence analysis detailed in Appendix~\ref{appendix:conv}, notably due to the coupling term $-\int_0^T\int_\Omega \epsilon \rho \Delta_x \varphi dx dt$, which lacks continuity for $\rho \in L^2$ and $\varphi \in H^1$.
This challenge is mitigated with the implementation of discretization. Through spatial discretization, a bilinear continuous operator emerges. However, its norm depends on the spatial grid dimensions $\Delta x$ in one dimension or $\Delta x\Delta y$ in two dimensions, leading to reduced step sizes as the granularity of spatial discretization increases.
For illustration, let's consider the one-dimensional scenario, noting that the two-dimensional case follows a similar pattern. After discretization, the space for $\varphi$ is represented as $X = \R^{n_t\times n_x}$, equipped with an inner product $\langle\cdot, \cdot\rangle_X$ defined by
\begin{equation*}
\langle \phi, \varphi\rangle_X = \sum_{i=1}^{n_x} \sum_{k=1}^{n_t} \left(\phi_{i,k}\varphi_{i,k} + (D_x^+\phi)_{i,k}(D_x^+ \varphi)_{i,k}\right) + \sum_{i=1}^{n_x} \sum_{k=2}^{n_t} (D_t^-\phi)_{i,k}(D_t^- \varphi)_{i,k} .
\end{equation*}
The spaces for $\rho$ and $\alpha$ align with Euclidean spaces. Consequently, the discretization of the term $-\int_0^T\int_\Omega \epsilon \rho \Delta_x \varphi dx dt$, divided by $\Delta x\Delta t$, is presented as $-\sum_{k=2}^{n_t}\sum_{i=1}^{n_x} \epsilon \rho_{i,k} (D_{xx} \varphi)_{i,k}$, showcasing bilinearity and continuity for $\varphi\in X$ and $\rho\in \R^{(n_t-1)\times n_x}$. However, the norm of this bilinear operator is inversely proportional to $\Delta x$, 
suggesting that the step sizes $\tau\varphi$, $\tau_\rho$, and $\tau_\alpha$ in the algorithm depend on $\Delta x$ and may lead to slower convergence compared to scenarios where $\epsilon = 0$.
Exploring solutions to this slower convergence requires further research. Our experiments with $H^2$-preconditioning for $\varphi$ or $H^1$-preconditioning for $\rho$ did not result in faster convergence.

\subsection{Numerical examples}  \label{sec:numerics_sto}
This section presents numerical examples under settings akin to those in Section~\ref{sec:numerics_det}, with a diffusion coefficient $\epsilon = 0.1$. Given that the functions $f$ and $L$ remain the same with those previously described, the updates for $\alpha$ remain unchanged. For each example, we initially deploy our proposed methodology to calculate the functions $\varphi$ and $\alpha$. Subsequently, we employ the Euler–Maruyama method on~\eqref{eqt:sde_openloop} to generate samples of the optimal processes $\gamma^*_s$ and controls $\alpha^*_s$.

\subsubsection{Quadratic Hamiltonian with spatial dependent coefficients}
\label{sec:sto_eg1}

\begin{figure}[htbp]
    \centering
    \begin{subfigure}{0.45\textwidth}
        \centering \includegraphics[width=\textwidth]{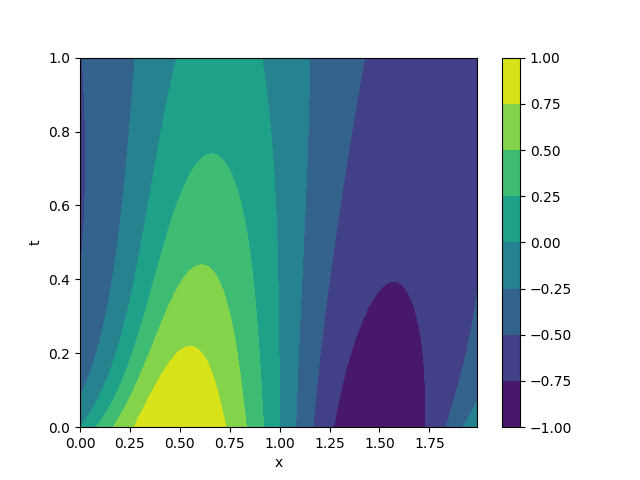}
        \caption{$\varphi$ in the $xt$-space}
    \end{subfigure}
    \hfill
    \begin{subfigure}{0.45\textwidth}
        \centering \includegraphics[width=\textwidth]{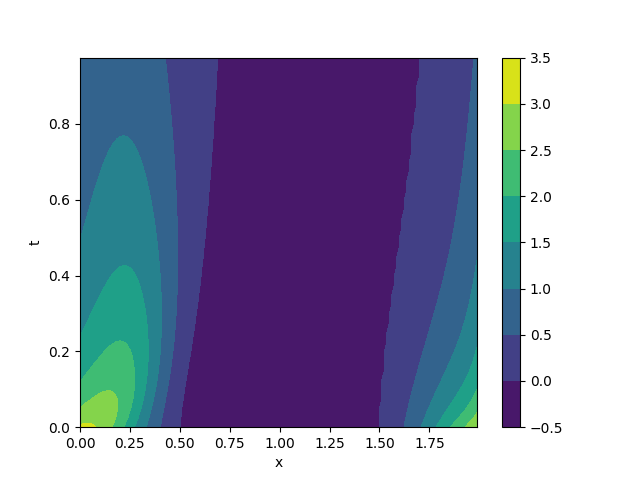}
        \caption{$\alpha$ in the $xt$-space}
    \end{subfigure}\\
    \begin{subfigure}{0.45\textwidth}
        \centering \includegraphics[width=\textwidth]{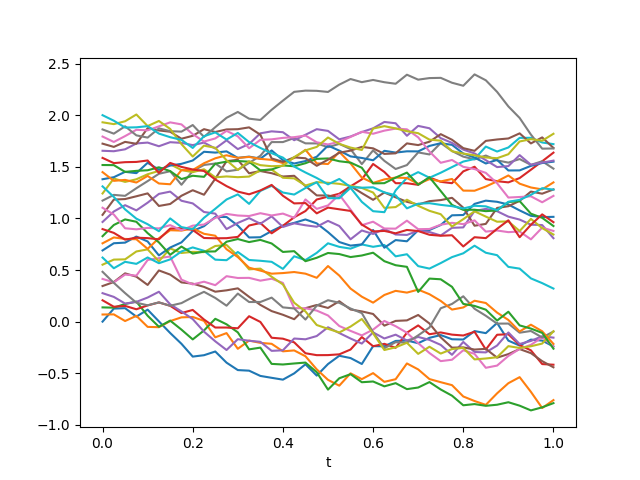}
        \caption{Samples of optimal trajectories $\gamma^*_s$}
    \end{subfigure}
    \hfill
    \begin{subfigure}{0.45\textwidth}
        \centering \includegraphics[width=\textwidth]{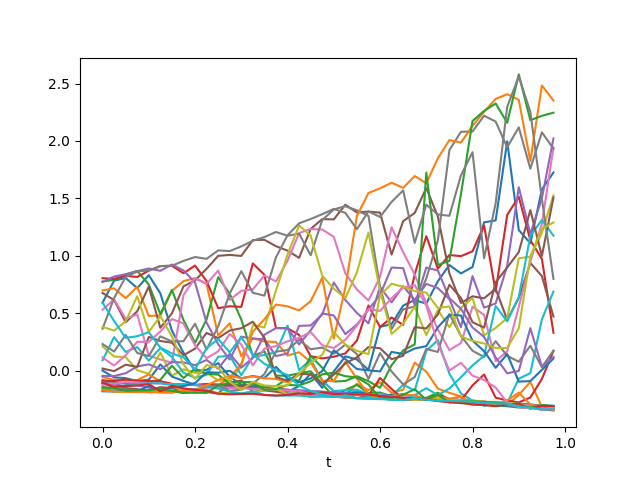}
        \caption{Samples of optimal controls $\alpha^*_s$}
    \end{subfigure}
    
    \caption{
    Visualization of the solution for the one-dimensional scenario discussed in Section~\ref{sec:sto_eg1}, using $n_t = 41$ and $n_x = 160$ grid points. Figures (a) and (b) showcase the level sets of the solution $\varphi$ to the viscous HJ PDE~\eqref{eqt:visc_HJ_initial}, along with the corresponding function $\alpha$ from~\eqref{eqt:visc_alp_initial}, which represents the time reversal of the feedback control function. Figures (c) and (d) depict several samples of optimal paths $\gamma^*_s$ and their associated open-loop optimal controls $\alpha^*_s$. These paths and control trajectories are the solutions to the stochastic optimal control problem~\eqref{eqt:soc_problem}, each beginning from a unique initial condition $x$.} 
    \label{fig:eg1_sto_1d}
\end{figure}

\begin{figure}[htbp]
    \centering
    \begin{subfigure}{0.45\textwidth}
        \centering \includegraphics[width=\textwidth]{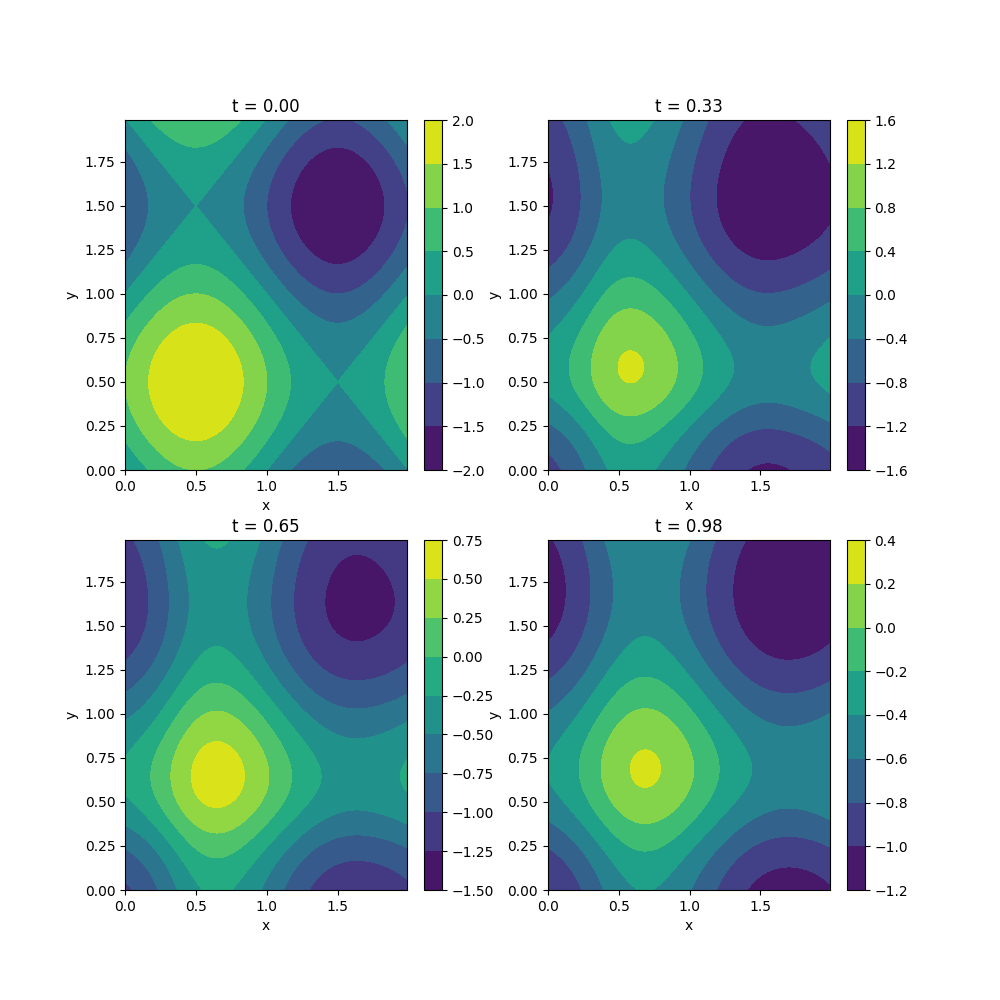}
        \caption{$(x,y)\mapsto \varphi(x,y,t)$ at different $t$}
    \end{subfigure}
    \\
    \begin{subfigure}{0.45\textwidth}
        \centering \includegraphics[width=\textwidth]{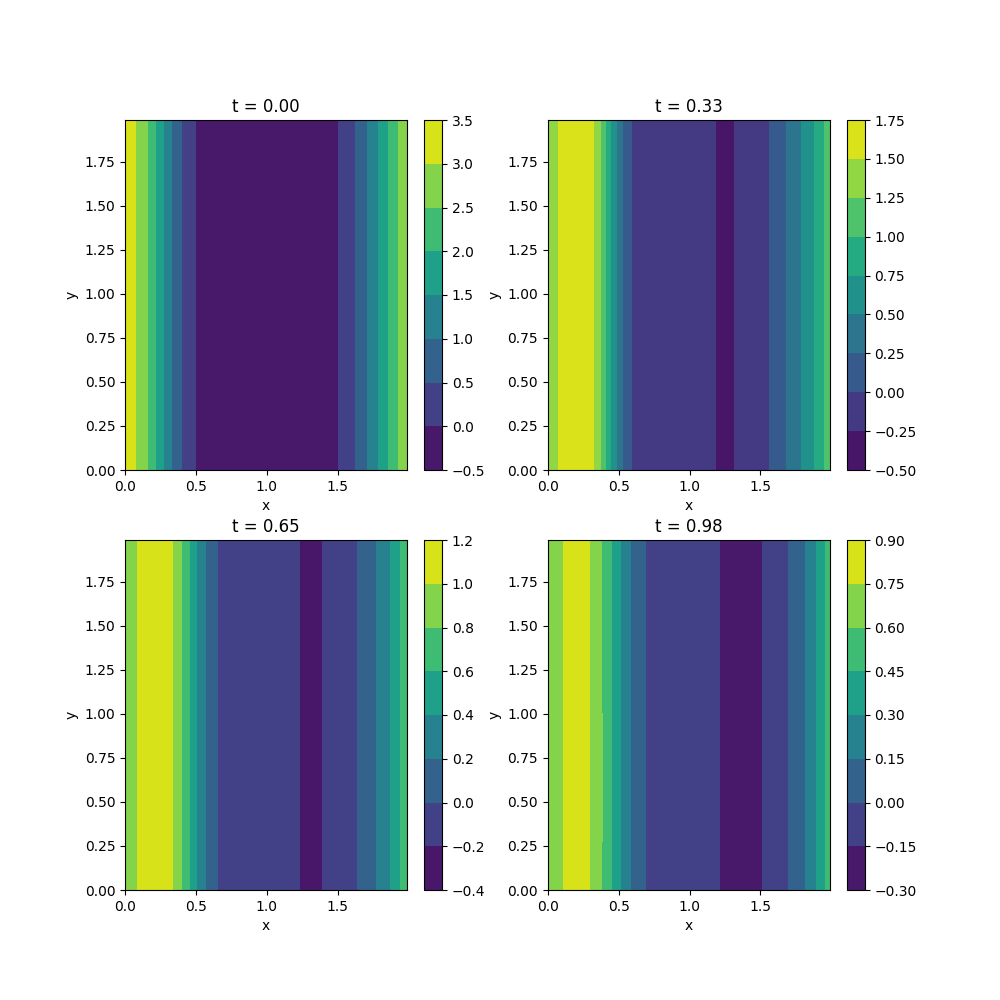}
        \caption{The first component of $(x,y)\mapsto \alpha(x,y,t)$ at different $t$}
    \end{subfigure}\hfill
    \begin{subfigure}{0.45\textwidth}
        \centering \includegraphics[width=\textwidth]{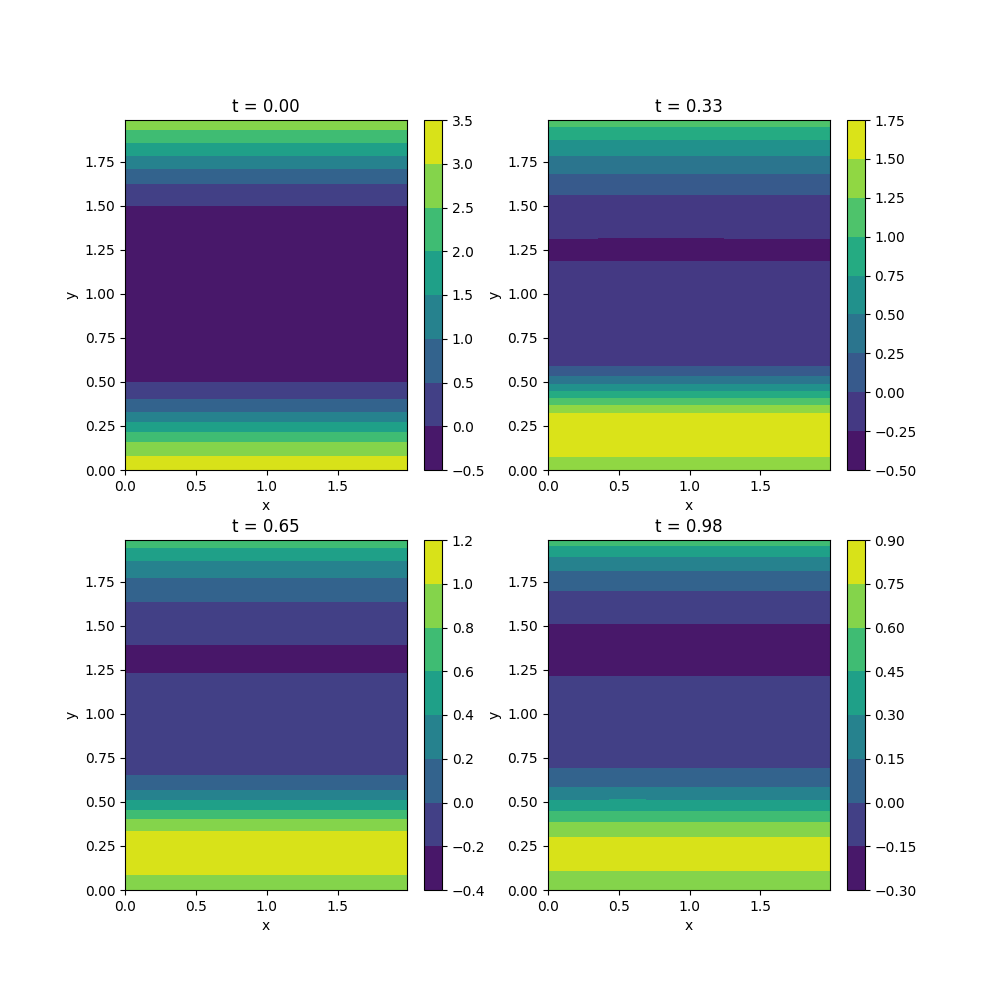}
        \caption{The second component of $(x,y)\mapsto \alpha(x,y,t)$ at different $t$}
    \end{subfigure}\\
    \begin{subfigure}{0.45\textwidth}
        \centering \includegraphics[width=\textwidth]{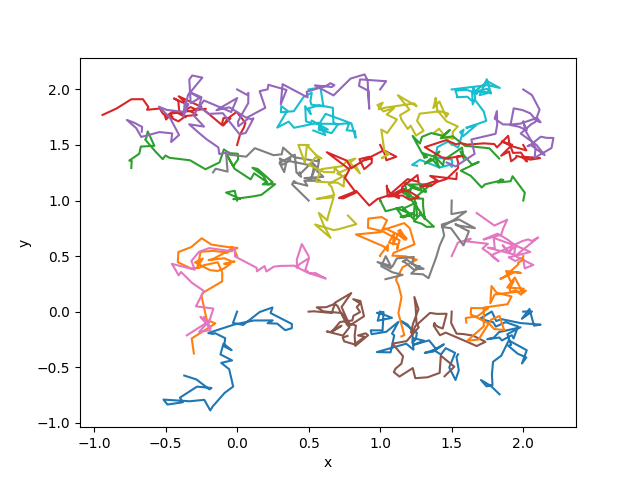}
        \caption{Samples of optimal trajectories $\gamma^*_s$}
    \end{subfigure}
    \hfill
    \begin{subfigure}{0.45\textwidth}
        \centering \includegraphics[width=\textwidth]{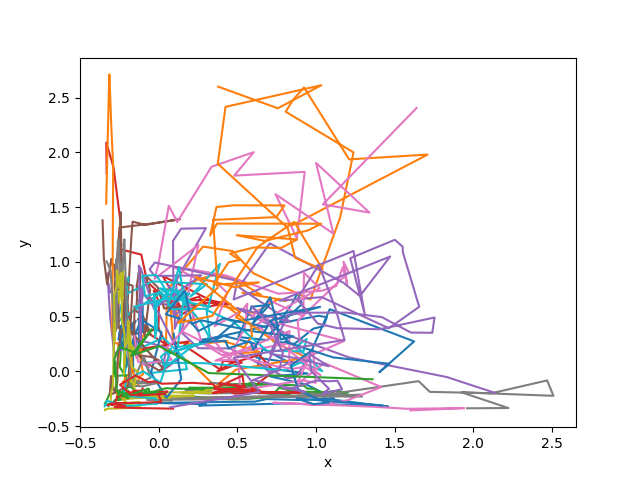}
        \caption{Samples of optimal controls $\alpha^*_s$}
    \end{subfigure}
    
    \caption{
    Depiction of the two-dimensional solution as discussed in Section~\ref{sec:sto_eg1}, utilizing $n_t = 41$ and $n_x = n_y = 160$ grid points. Figure (a) illustrates the level sets of the solution $\varphi(\cdot, t)$ to the viscous HJ PDE~\eqref{eqt:visc_HJ_initial} at different times $t$. Figures (b) and (c) show the first and second components, respectively, of the associated function $\alpha(\cdot, t)$ from~\eqref{eqt:visc_alp_initial} at various times $t$, which depict the time reversal of the feedback control function. Figures (d) and (e) present several samples of optimal trajectories $\gamma^*_s$ along with their corresponding open-loop optimal controls $\alpha^*_s$. These trajectories and control strategies solve the stochastic optimal control problem specified in~\eqref{eqt:soc_problem}, starting from distinct initial conditions $x$. Notably, both $\gamma^*_s$ and $\alpha^*_s$ take values in $\R^2$. For visualization, they are plotted within the spatial domain, excluding the time dimension for clarity.} 
    \label{fig:eg1_sto_2d}
\end{figure}

We employ the same spatial domain $\Omega$ and functions \(f, g, L\) as in Section~\ref{sec:eg1}. The solution of the one-dimensional case is depicted in Figure~\ref{fig:eg1_sto_1d}, while the solution for the two-dimensional problem is illustrated in Figure~\ref{fig:eg1_sto_2d}. Upon examining the level set diagrams of \(\varphi\) and \(\alpha\), it becomes evident that the numerical solutions exhibit greater smoothness compared to those presented in Section~\ref{sec:eg1}. Despite the stochastic nature of the optimal trajectories, which allows us to only display certain samples, there is a discernible pattern where the trajectories seem to cluster around the minimum value of \(g\) at \(t=T\).

\subsubsection{One-homogeneous Hamiltonian with spatial dependent coefficients}\label{sec:numerics_sto_eg2}
\begin{figure}[htbp]
    \centering
    \begin{subfigure}{0.45\textwidth}
        \centering \includegraphics[width=\textwidth]{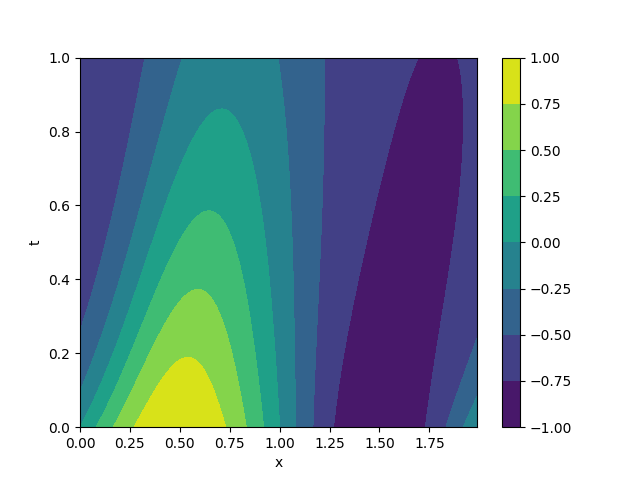}
        \caption{$\varphi$ in the $xt$-space}
    \end{subfigure}
    \hfill
    \begin{subfigure}{0.45\textwidth}
        \centering \includegraphics[width=\textwidth]{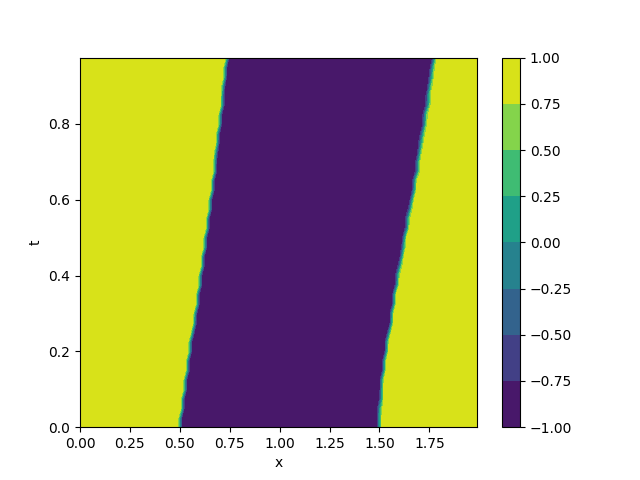}
        \caption{$\alpha$ in the $xt$-space}
    \end{subfigure}\\
    \begin{subfigure}{0.45\textwidth}
        \centering \includegraphics[width=\textwidth]{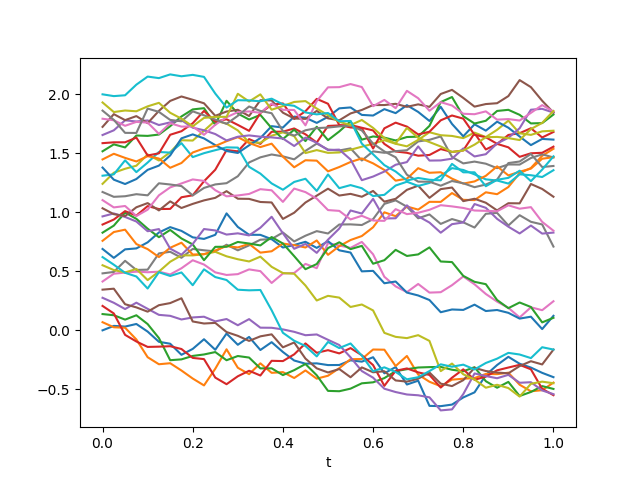}
        \caption{Samples of optimal trajectories $\gamma^*_s$}
    \end{subfigure}
    \hfill
    \begin{subfigure}{0.45\textwidth}
        \centering \includegraphics[width=\textwidth]{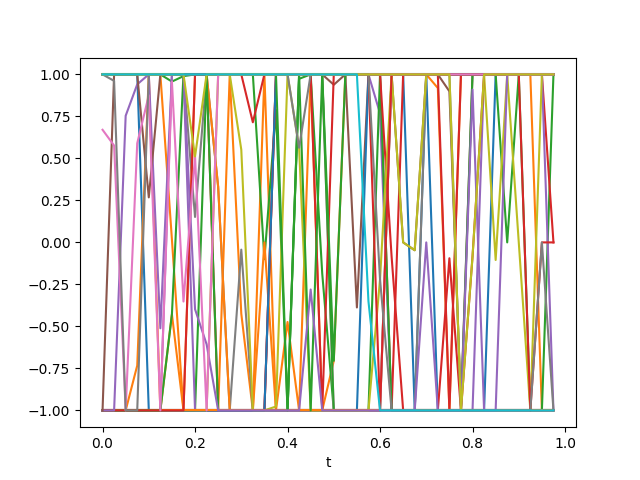}
        \caption{Samples of optimal controls $\alpha^*_s$}
    \end{subfigure}
    
    \caption{
    Visualization of the solution for the one-dimensional scenario discussed in Section~\ref{sec:numerics_sto_eg2}, using $n_t = 41$ and $n_x = 160$ grid points. Figures (a) and (b) showcase the level sets of the solution $\varphi$ to the viscous HJ PDE~\eqref{eqt:visc_HJ_initial}, along with the corresponding function $\alpha$ from~\eqref{eqt:visc_alp_initial}, which represents the time reversal of the feedback control function. Figures (c) and (d) depict several samples of optimal paths $\gamma^*_s$ and their associated open-loop optimal controls $\alpha^*_s$. These paths and control trajectories are the solutions to the stochastic optimal control problem~\eqref{eqt:soc_problem}, each beginning from a unique initial condition $x$.} 
    \label{fig:eg2_sto_1d}
\end{figure}

\begin{figure}[htbp]
    \centering
    \begin{subfigure}{0.45\textwidth}
        \centering \includegraphics[width=\textwidth]{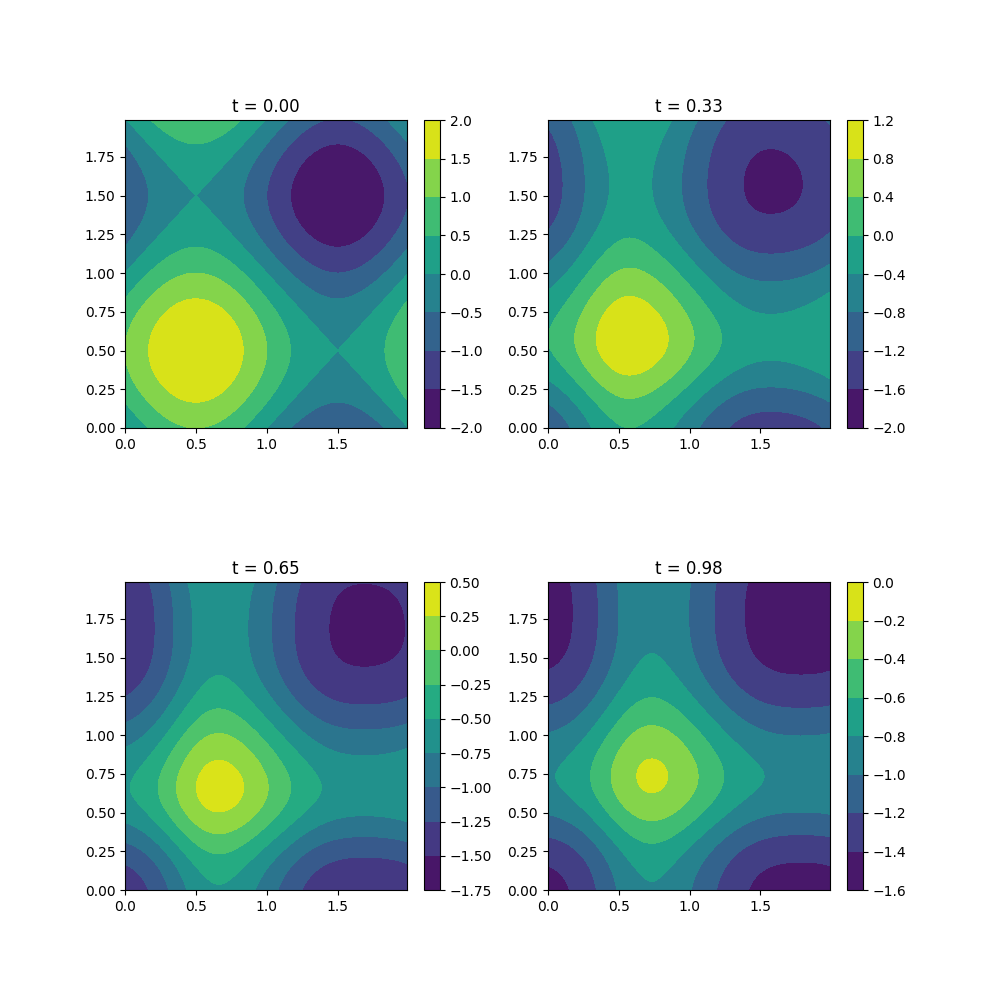}
        \caption{$(x,y)\mapsto \varphi(x,y,t)$ at different $t$}
    \end{subfigure}
    \\
    \begin{subfigure}{0.45\textwidth}
        \centering \includegraphics[width=\textwidth]{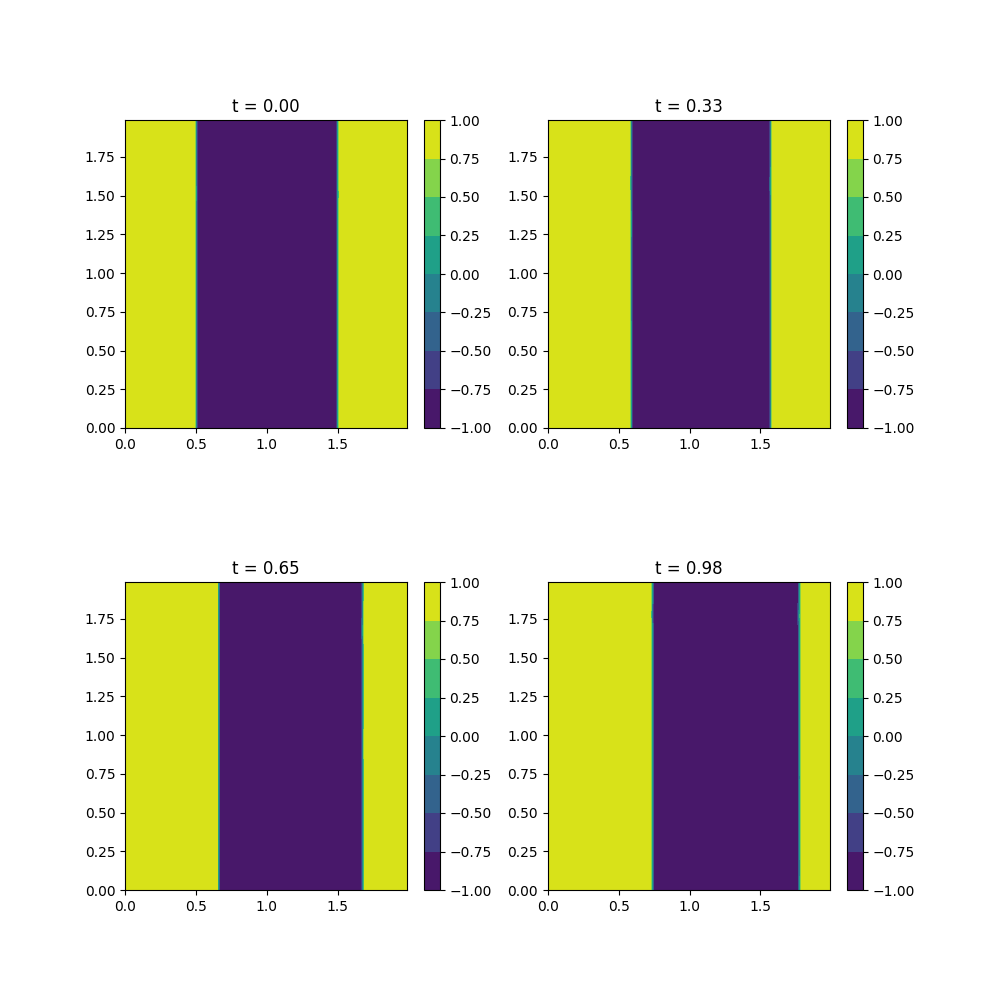}
        \caption{The first component of $(x,y)\mapsto \alpha(x,y,t)$ at different $t$}
    \end{subfigure}\hfill
    \begin{subfigure}{0.45\textwidth}
        \centering \includegraphics[width=\textwidth]{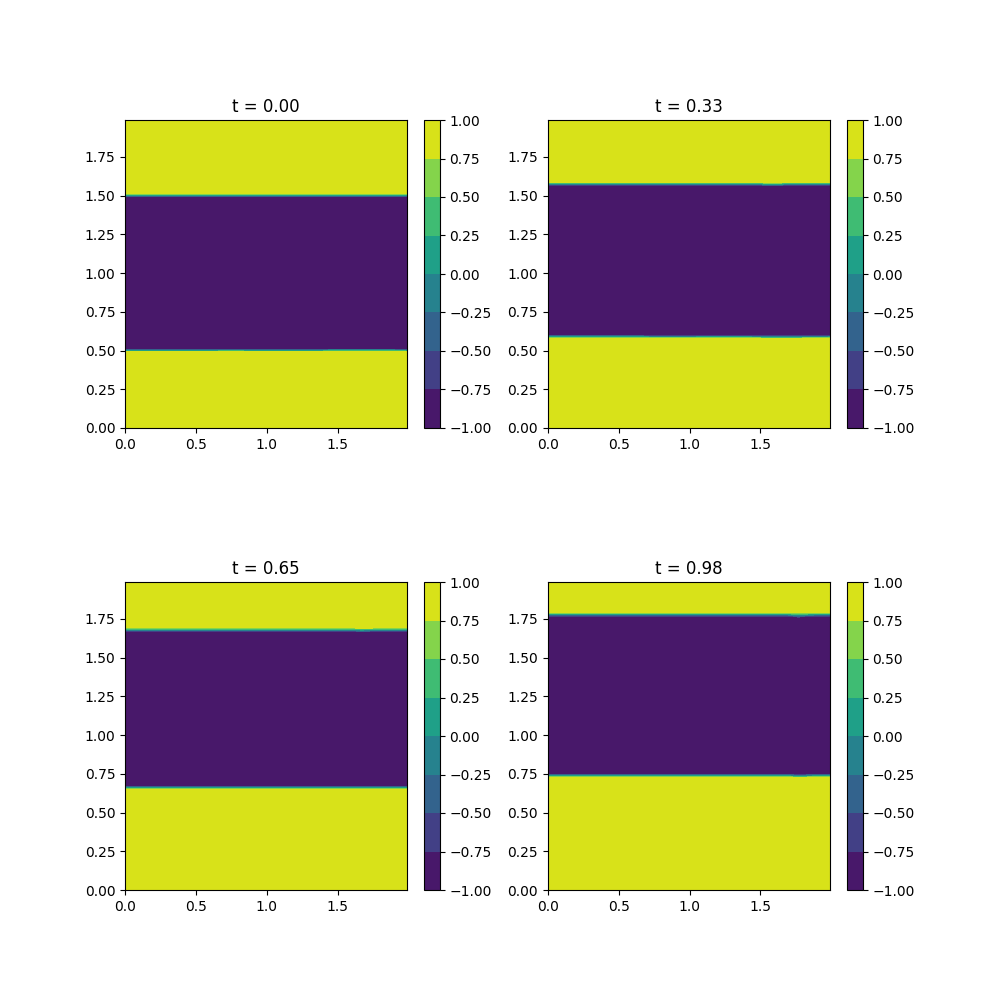}
        \caption{The second component of $(x,y)\mapsto \alpha(x,y,t)$ at different $t$}
    \end{subfigure}\\
    \begin{subfigure}{0.45\textwidth}
        \centering \includegraphics[width=\textwidth]{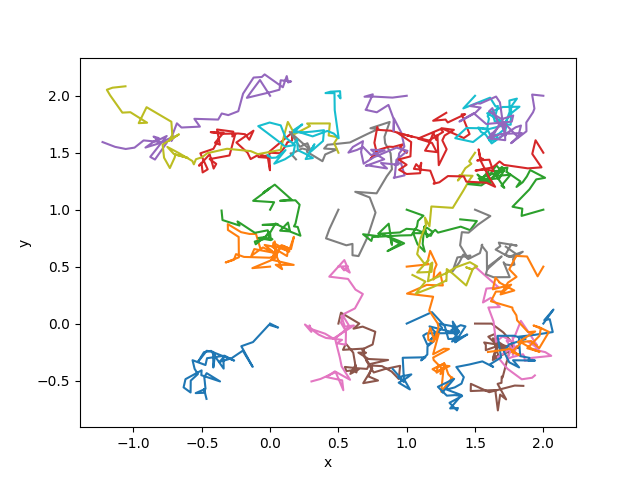}
        \caption{Samples of optimal trajectories $\gamma^*_s$}
    \end{subfigure}
    \hfill
    \begin{subfigure}{0.45\textwidth}
        \centering \includegraphics[width=\textwidth]{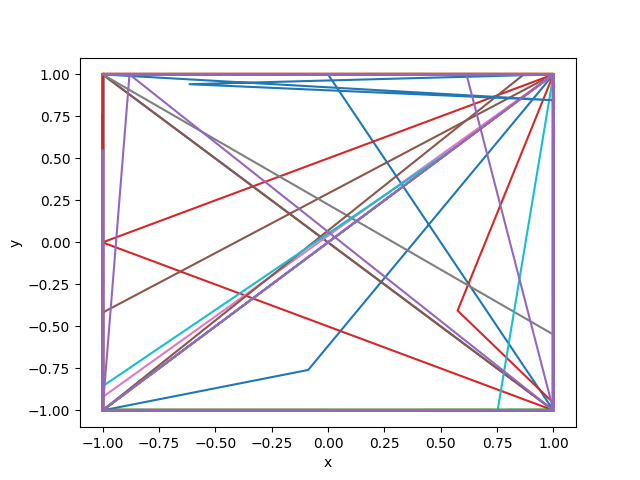}
        \caption{Samples of optimal controls $\alpha^*_s$}
    \end{subfigure}
    
    \caption{
    Depiction of the two-dimensional solution as discussed in Section~\ref{sec:numerics_sto_eg2}, utilizing $n_t = 41$ and $n_x = n_y = 160$ grid points. Figure (a) illustrates the level sets of the solution $\varphi(\cdot, t)$ to the viscous HJ PDE~\eqref{eqt:visc_HJ_initial} at different times $t$. Figures (b) and (c) show the first and second components, respectively, of the associated function $\alpha(\cdot, t)$ from~\eqref{eqt:visc_alp_initial} at various times $t$, which depict the time reversal of the feedback control function. Figures (d) and (e) present several samples of optimal trajectories $\gamma^*_s$ along with their corresponding open-loop optimal controls $\alpha^*_s$. These trajectories and control strategies solve the stochastic optimal control problem specified in~\eqref{eqt:soc_problem}, starting from distinct initial conditions $x$. Notably, both $\gamma^*_s$ and $\alpha^*_s$ take values in $\R^2$. For visualization, they are plotted within the spatial domain, excluding the time dimension for clarity.} 
    \label{fig:eg2_sto_2d}
\end{figure}

This example follows a similar setup to that described in Section~\ref{sec:eg2}, with a diffusion parameter of $\epsilon = 0.1$. The outcomes for the one-dimensional scenario are depicted in Figure~\ref{fig:eg2_sto_1d}, and those for the two-dimensional scenario are presented in Figure~\ref{fig:eg2_sto_2d}. The solution $\varphi$ exhibits greater smoothness compared to scenarios lacking a diffusion term. Nonetheless, the function $\alpha$ still demonstrates discontinuities. This occurrence is attributable to $\alpha$ employing the same expression~\eqref{eqt:visc_alp_initial} as found in the deterministic optimal control problem in Section~\ref{sec:oc_problem}, thus the rationale from Section~\ref{sec:eg2} remains applicable. It should be noted that, due to the element of randomness, the trajectories for optimal control display a higher frequency of jumps in comparison to the deterministic configuration.

\subsubsection{Newton mechanics}\label{sec:numerics_newton_sto}
\begin{figure}[htbp]
    \centering
    \begin{subfigure}{0.45\textwidth}
        \centering \includegraphics[width=\textwidth]{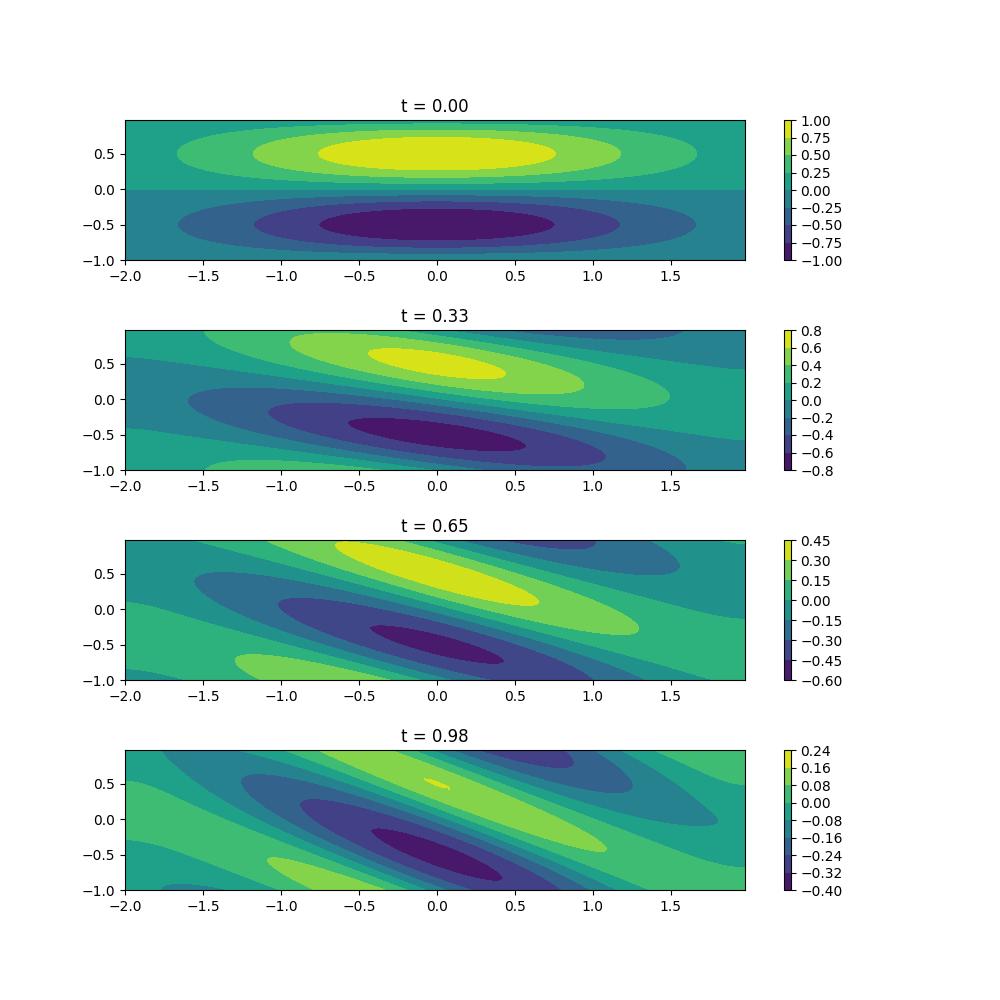}
        \caption{$(x,y)\mapsto \varphi(x,y,t)$ at different $t$}
    \end{subfigure}
    \hfill
    \begin{subfigure}{0.45\textwidth}
        \centering \includegraphics[width=\textwidth]{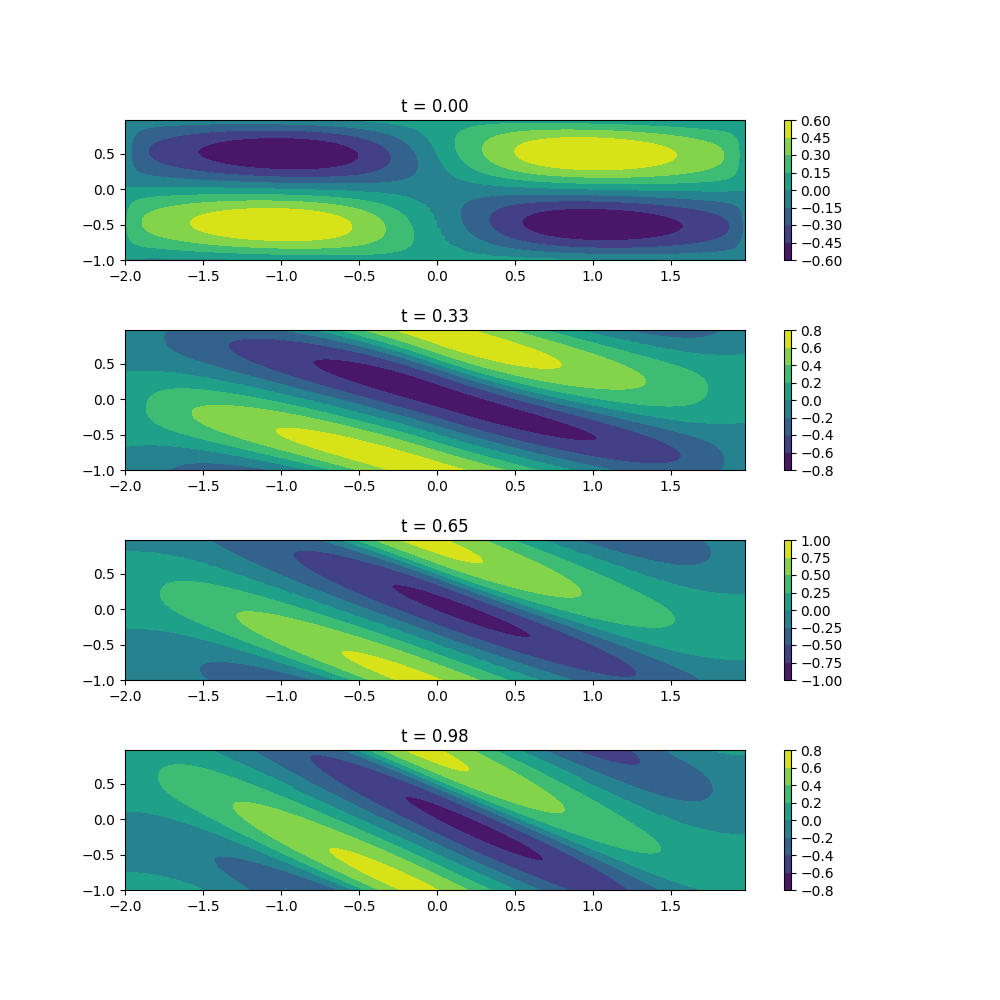}
        \caption{$(x,y)\mapsto \alpha(x,y,t)$ at different $t$}
    \end{subfigure}\\
    \begin{subfigure}{0.45\textwidth}
        \centering \includegraphics[width=\textwidth]{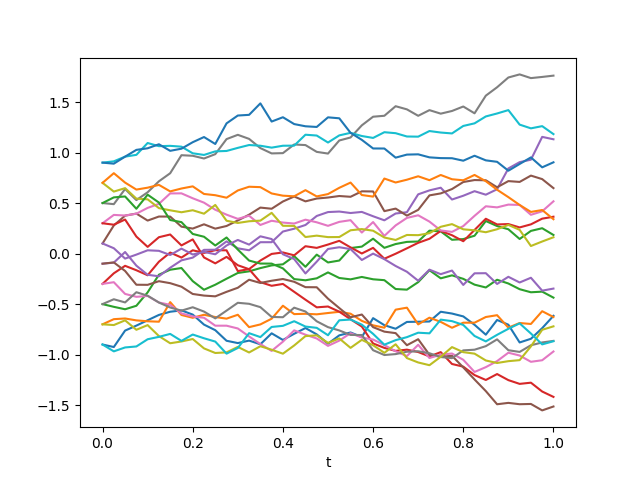}
        \caption{Optimal position samples $(\gamma^*_2)_s$}
    \end{subfigure}\hfill
    \begin{subfigure}{0.45\textwidth}
        \centering \includegraphics[width=\textwidth]{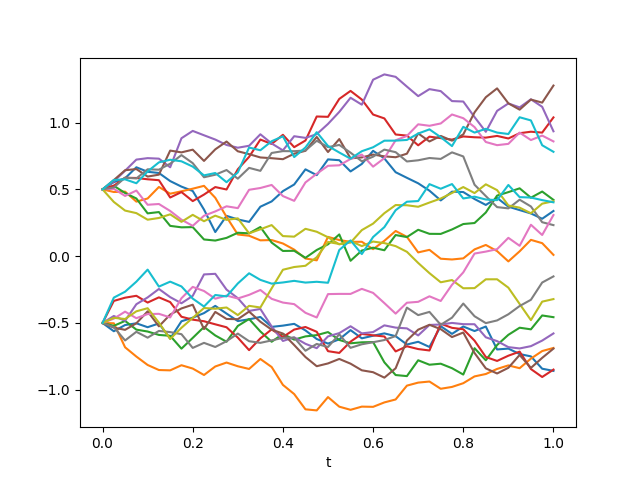}
        \caption{Optimal velocity samples $(\gamma^*_1)_s$}
    \end{subfigure}
    \\
    \begin{subfigure}{0.45\textwidth}
        \centering \includegraphics[width=\textwidth]{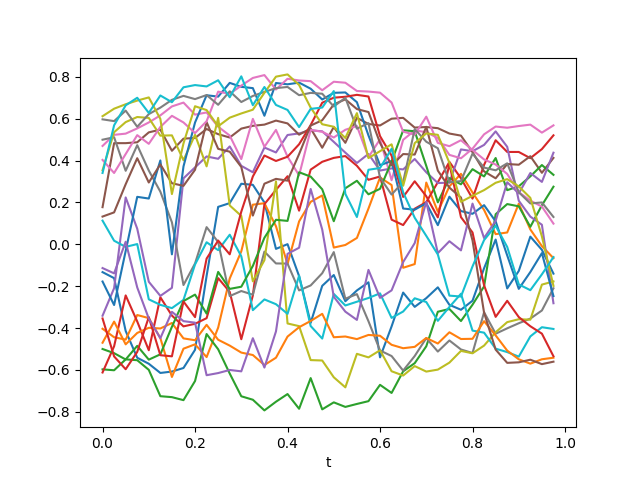}
        \caption{Optimal control samples $\alpha^*_s$}
    \end{subfigure}
    
    \caption{
    Visualization of the solution for example discussed in Section~\ref{sec:numerics_newton_sto}, using $n_t = 41$, $n_x = 160$, and $n_y=80$ grid points. Figures (a) and (b) showcase the level sets of the solution $\varphi$ to the viscous HJ PDE~\eqref{eqt:visc_HJ_initial}, along with the corresponding function $\alpha$ from~\eqref{eqt:visc_alp_initial}, which represents the time reversal of the feedback control function. Figures (c), (d), (e) depict the first and second components of several samples of optimal paths $\gamma^*_s$ and their associated open-loop optimal controls $\alpha^*_s$. These paths and control trajectories are the solutions to the stochastic optimal control problem~\eqref{eqt:soc_problem}, each beginning from a unique initial condition $x$.} 
    \label{fig:eg_newton_sto}
\end{figure}

This section employs the same functions $f$, $L$, and $g$ detailed in Section~\ref{sec:numerical_newton_det}. We present the numerical solution for a diffusion coefficient of $\epsilon = 0.1$ in Figure~\ref{fig:eg_newton_sto}.

\section{Summary}\label{sec:summary}

In this study, we present a framework for solving optimal control problems and HJ PDEs through a novel saddle point formulation, which is solved using an algorithm derived from the preconditioned PDHG method. This approach extends to solve stochastic optimal control problems and associated viscous HJ PDEs. Through a series of numerical examples, we illustrate the efficacy of our algorithm, especially in handling Hamiltonians which depend on spatial and temporal variables. A key advantage of our method is the use of implicit time discretization, which allows for larger time steps than explicit methods, thereby enhancing computational speed.
The core strength of our method lies in its straightforward saddle point formulation, linearly related to the solution $\varphi$ of the HJ PDEs, a feature facilitated by leveraging the inherent relationship between optimal control problems and HJ PDEs. This simplicity enables the method to support explicit updates or to benefit from parallel computation, making it versatile and efficient.

Although our method is of first order, it shows potential as an initial step towards methods of higher accuracy, especially in applications that require small errors. Future work could explore integrating this method with higher-order schemes or employing it within machine learning frameworks, where its formulation could inspire novel loss function designs for solving optimal control problems and related HJ PDEs.

\begin{acknowledgements}
We would like to thank Dr. Donggun Lee for his valuable discussions. Additionally, we are thankful to Dr. Levon Nurbekyan for his assistance with the analysis presented in Appendix~\ref{appendix:conv}.
\end{acknowledgements}

\section*{Conflict of interest}

The authors declare that they have no conflict of interest.

\bibliographystyle{spmpsci}      % mathematics and physical sciences
\bibliography{references}

\appendix

\section{The convergence of the PDHG algorithms in Section~\ref{sec:saddle_point_cont}}\label{appendix:conv}

In this section, we show the convergence of the PDHG algorithm to a saddle point in certain cases.

Let $\Omega\subset \Rn$ be a bounded rectangular domain with periodic boundary condition.
We consider the case when $f_{x,t}$ is an affine function for any $x\in \Omega$ and $t\in [0,T]$, i.e., $f(x,t,\alpha) = A(x,t)\alpha + b(x,t)$ for continuous functions $A\colon \Omega\times [0,T]\to \R^{n\times m}$ (where $\R^{n\times m}$ denotes the set of matrices with $n$ rows and $m$ columns) and $b\colon \Omega\times [0,T]\to\Rn$.
Assume $L_{x,t}\colon \R^m\to [0,+\infty]$ is proper, convex, lower semi-continuous, $1$-coercive, and satisfies $L_{x,t}(0) = 0$ for any $x\in \Omega$ and $t\in [0,T]$. 
To simplify the notation, we use $A_{x,t}$ and $b_{x,t}$ to denote $A(x,T-t)$ and $b(x,T-t)$. Then, we have $f_{x,t}(\alpha) = A_{x,t}\alpha + b_{x,t}$.
In this case, the saddle point problem~\eqref{eqt:saddle_oc} becomes 
\begin{equation}\label{eqt:saddle_point_rem21}
\begin{split}
\min_{\substack{\varphi\\ \varphi(x,0)=g(x)}}\max_{\rho\geq 0, \alpha}  \int_0^T \int_\Omega \rho(x,t)\Big(\frac{\partial \varphi(x,t)}{\partial t} -\langle A_{x,t}\alpha(x,t) + b_{x,t}, \nabla_x \varphi(x,t)\rangle\\ 
- L_{x,t}(\alpha(x,t)) \Big)dxdt - c\int_\Omega \varphi(x,T)dx.
\end{split}
\end{equation}
With a change of variable $m(x,t) = \rho(x,t)\alpha(x,t)$, the problem becomes 
\begin{equation}\label{eqt:saddle_m_cont}
\begin{adjustbox}{width=0.99\textwidth}$
\begin{split}
\min_{\substack{\varphi\\ \varphi(x,0)=g(x)}}\max_{\rho\geq 0, m}  \int_0^T \int_\Omega \rho(x,t)\frac{\partial \varphi(x,t)}{\partial t} - m(x,t)^TA_{x,t}^T \nabla_x \varphi(x,t) - \rho(x,t)b_{x,t}^T \nabla_x \varphi(x,t)\\ 
- \rho(x,t)L_{x,t}\left(\frac{m(x,t)}{\rho(x,t)}\right) dxdt - c\int_\Omega \varphi(x,T)dx,
\end{split}
$\end{adjustbox}
\end{equation}
where $\rho L_{x,t}(\frac{m}{\rho})$ is defined to be $\ind_{\{0\}}(m)$ if $\rho$ is zero in orde to get a lower semi-continuous function (by~\cite[Prop.X.1.2.1]{Urruty1993Convex} and the assumption that $L_{x,t}$ is $1$-coercive which implies $\dom L_{x,t}^* = \R^m$). Recall that $\ind_C$ denotes the indicator function of a set $C$ which takes value $0$ at points in $C$ and $+\infty$ otherwise.
Let $X = H^1(\Omega\times [0,T])$ and $Y = L^2(\Omega\times [0,T])\times L^2(\Omega\times [0,T]; \R^m)$. 
Now, we prove the convergence of the PDHG algorithm applied to~\eqref{eqt:saddle_m_cont} for $\varphi \in X$ and $(\rho,m)\in Y$.

Define the operators $K\colon X\times Y\to \R$, $G\colon X\to \R\cup\{+\infty\}$, and $\hat F\colon Y\to \R\cup\{+\infty\}$ as follows
\begin{equation*}
\begin{adjustbox}{width=0.99\textwidth}$
\begin{split}
K(\varphi, \rho,m) &= \int_0^T \int_\Omega \rho(x,t)\frac{\partial \varphi(x,t)}{\partial t} - m(x,t)^TA_{x,t}^T \nabla_x \varphi(x,t) - \rho(x,t)b_{x,t}^T \nabla_x \varphi(x,t) dxdt,\\
G(\varphi) &= -c\int_\Omega \varphi(x,T)dx + \ind_{\{\varphi \colon \varphi(\cdot,0) = g\}}(\varphi), \\
\hat F(\rho, m) &= \int_0^T\int_{\Omega} \rho(x,t)L_{x,t}\left(\frac{m(x,t)}{\rho(x,t)}\right) + \ind_{[0,+\infty)} (\rho(x,t)) dxdt.
\end{split}
$\end{adjustbox}
\end{equation*}
It is straightforward to check that $K$ is bilinear with respect to $\varphi$ and $(\rho, m)$, and $G$, $\hat F$ are convex and proper functions (recall that a function is proper means it is not constantly $+\infty$).

Now, we prove that $G$ is lower semi-continuous on $X$. Let $\varphi_j$ be a sequence converging to $\varphi$ in $X$. By trace theorem, $\varphi_j(\cdot, T)$ converges to $\varphi(\cdot, T)$ in $L^2(\Omega)$. Since $\Omega$ is bounded, we have $\lim_{j\to\infty}\int_{\Omega} \varphi_j(x,T)dx = \int_{\Omega} \varphi(x,T)dx$. Similarly, by trace theorem, $\varphi_j(\cdot, 0)$ converges to $\varphi(\cdot, 0)$ in $L^2(\Omega)$, and hence there exists a subsequence of $\varphi_j$, denoted by $\varphi_{j_i}$, such that $\varphi_{j_i}(\cdot, 0)$ converges to $\varphi(\cdot, 0)$ almost everywhere. As a result, $\liminf_{j\to\infty}\ind_{\{\varphi \colon \varphi(\cdot,0) = g\}}(\varphi_j) \geq \ind_{\{\varphi \colon \varphi(\cdot,0) = g\}}(\varphi)$ holds. Then, we conclude that $G$ is lower semi-continuous.

We also prove that $\hat F$ is lower semi-continuous on $Y$. Let $(\rho_j,m_j)$ be a convergent sequence in $Y$ with limit $(\rho,m)$. We need to prove $\liminf_{j\to\infty}\hat F(\rho_j, m_j)\geq \hat F(\rho, m)$. Without loss of generality, we can assume $\rho_j$ is non-negative almost everywhere for all $j$. Since $L$ is non-negative, $\rho_{j}(x,t)L_{x,t}(\frac{m_{j}(x,t)}{\rho_{j}(x,t)})$ is non-negative almost everywhere for all $j$. Since $(\rho_j, m_j)$ converges to $(\rho, m)$ in $L^2$, there exists a subsequence $(\rho_{j_i},m_{j_i})$ converging to $(\rho,m)$ pointwisely almost everywhere. According to the lower semi-continuity of $L_{x,t}$ (and using the $1$-coercivity of $L_{x,t}$ in the case of $\rho(x,t)=0$), we obtain 
\begin{equation*}
\liminf_{i\to\infty} \rho_{j_i}(x,t)L_{x,t}\left(\frac{m_{j_i}(x,t)}{\rho_{j_i}(x,t)}\right)\geq \rho(x,t)L_{x,t}\left(\frac{m(x,t)}{\rho(x,t)}\right) \text{ a.e. } x\in\Omega, t\in [0,T].
\end{equation*}
Then, we get $\liminf_{j\to\infty}\hat F(\rho_j, m_j)\geq \hat F(\rho, m)$ by Fatou's lemma. 
Define $F$ as the Fenchel-Legendre transform of $\hat F$, i.e., $F = \hat F^*$.
Since $\hat F$ is proper, convex, and lower semi-continuous on $Y$, by~\cite[Thm 13.37]{Bauschke2011Convex}, we have $\hat F = F^*$.

Now, we prove the continuity of $K$ as follows
\begin{equation*}
\begin{split}
|K(\varphi, \rho, m)| &\leq \max\left\{1, \sup_{x,t} \|b_{x,t}\|_\infty\right\}\|\rho\|_{L^2}\|\varphi\|_{H^1} + \sup_{x,t, i,j} |A_{x,t,i,j}|\|m\|_{L^2}\|\varphi\|_{H^1}\\
&\leq \left(\max\left\{1, \sup_{x,t} \|b_{x,t}\|_\infty\right\}^2 + \sup_{x,t, i,j} |A_{x,t,i,j}|^2\right) \|(\rho, m)\|_{L^2}\|\varphi\|_{H^1}.
\end{split}
\end{equation*}
Therefore, $K$ is continuous with operator norm $\|K\|\leq \max\{1, \sup_{x,t} \|b_{x,t}\|_\infty\}^2 + \sup_{x,t, i,j} |A_{x,t,i,j}|^2$.

Then, the saddle point problem~\eqref{eqt:saddle_m_cont} is written as the standard form for PDHG algorithm:
\begin{equation*}
\min_{\varphi \in X} \max_{(\rho,m)\in Y} K(\varphi, \rho, m) + G(\varphi) - F^*(\rho, m),
\end{equation*}
where $K$ is bilinear and continuous with respect to $\varphi$ and $(\rho, m)$, $G$ and $F$ are both proper, convex and lower semi-continuous. 
The $\ell$-th iteration of the PDHG updates is given by
\begin{equation} \label{eqt:pdhg_det_cont_appendix}
\begin{adjustbox}{width=0.99\textwidth}$
\begin{dcases}
(\rho^{\ell+1}, m^{\ell+1}) = \argmax_{(\rho, m)\in Y} \left\{K(\tilde\varphi^\ell, \rho, m) - F^*(\rho, m) - \frac{1}{2\tau_{\rho,m}}\left(\|\rho - \rho^\ell\|^2_{L^2} +\|m - m^\ell\|^2_{L^2}\right)\right\},\\
\varphi^{\ell+1} = \argmin_{\varphi \in X} \left\{K(\varphi, \rho^{\ell+1}, m^{\ell+1}) + G(\varphi) + \frac{1}{2\tau_\varphi}\|\varphi - \varphi^\ell\|_{H^1}^2\right\},\\
\tilde\varphi^{\ell+1} = 2\varphi^{\ell+1} - \varphi^\ell.
\end{dcases}
$\end{adjustbox}
\end{equation}
Then, the convergence of these updates follows from~\cite{darbon2021accelerated} if the stepsizes $\tau_{\varphi}$ and $\tau_{\rho, m}$ are small enough. To be specific, assume the problem~\eqref{eqt:saddle_m_cont} has a saddle point, the stepsizes satisfy $\tau_{\varphi}\tau_{\rho, m} < \frac{1}{\|K\|^2}$, where $\|K\|$ is estimated above, and let $\varphi^j, \rho^\ell, m^\ell$ be the functions given by the algorithm~\eqref{eqt:pdhg_det_cont_appendix}. 
Then, the sequence $\{\bar \varphi^{\ell}, \bar \rho^{\ell}, \bar m^{\ell}\}_\ell$ has a subsequence which converges weakly to a saddle point of~\eqref{eqt:saddle_m_cont}, where $\bar \varphi^{\ell}, \bar \rho^{\ell}, \bar m^{\ell}$ are the average functions defined by
$(\bar \varphi^{\ell}, \bar \rho^{\ell}, \bar m^{\ell}) = \frac{1}{\ell}\sum_{j=1}^{\ell} (\varphi^j, \rho^j, m^j)$.

\begin{rem}[Connection of the PDHG updates above with our proposed updates]\label{rem:A1_conv_pdhg_rho_alp}
First, we work on the details about the update for $\varphi$ as follows. After simplifying the updates for $\varphi$, we get
\begin{equation*}
\begin{adjustbox}{width=0.99\textwidth}$
\begin{split}
\varphi^{\ell+1} &= \argmin_{\substack{\varphi \in X\\\varphi(\cdot,0)=g}} \left\{\int_0^T \int_\Omega \rho\frac{\partial \varphi}{\partial t} - m(x,t)^TA_{x,t}^T \nabla_x \varphi - \rho(x,t)b_{x,t}^T \nabla_x \varphi dxdt -c\int_\Omega \varphi(x,T)dx + \frac{1}{2\tau_\varphi}\|\varphi - \varphi^\ell\|_{H^1}^2\right\}\\
&= \argmin_{\substack{\varphi \in X\\\varphi(\cdot,0)=g}} \left\{\int_0^T \int_\Omega \left(-\partial_t \rho^{\ell+1} + \nabla_x\cdot (A_{x,t}m^{\ell+1} + \rho^{\ell+1}b_{x,t})\right)\varphi dxdt +\int_\Omega (\rho^{\ell+1}(x,T)-c) \varphi(x,T)dx + \frac{1}{2\tau_\varphi}\|\varphi - \varphi^\ell\|_{H^1}^2\right\}.
\end{split}
$\end{adjustbox}
\end{equation*}
If $\rho^{\ell+1}$ satisfies the terminal condition $\rho^{\ell+1}(x,T) = c$ for all $x\in \Omega$, the first order optimality condition gives
\begin{equation*}
0 = \tau_\varphi \left(-\partial_t \rho^{\ell+1} + \nabla_x\cdot (A_{x,t}m^{\ell+1} + \rho^{\ell+1}b_{x,t})\right) + (I - \Delta )(\varphi^{\ell+1} - \varphi^{\ell}),
\end{equation*}
which implies 
\begin{equation*}
\varphi^{\ell+1} = \varphi^{\ell} + \tau_{\varphi}(I-\Delta)^{-1}\left(\partial_t \rho^{\ell+1} - \nabla_x \cdot (A_{x,t}m^{\ell+1}(x,t) + \rho^{\ell+1}(x,t)b_{x,t})\right).
\end{equation*}
With a change of variable $\alpha = \frac{m}{\rho}$, this update on $\varphi$ is the same as the third line in~\eqref{eqt:pdhg_det_cont}.

Then, we focus on the update for $(\rho, m)$ (the first line in~\eqref{eqt:pdhg_det_cont_appendix}), which, after the change of variable, becomes
\begin{equation*}
\max_{\rho, \alpha} \left\{\mathcal{L}(\tilde\varphi^{\ell}, \rho, \alpha) - \frac{1}{2\tau_{\rho,m}}\|\rho - \rho^{\ell}\|_{L^2}^2 - \frac{1}{2\tau_{\rho,m}}\|\rho\alpha - \rho^{\ell}\alpha^{\ell}\|_{L^2}^2\right\},
\end{equation*}
where $\mathcal{L}$ is the objective function in~\eqref{eqt:saddle_point_rem21}. To numeically solve this problem, we apply the alternative updates on $\rho$ and $\alpha$. Fixing $\alpha = \alpha^\ell$, the update for $\rho$ gives the first line in~\eqref{eqt:pdhg_det_cont} with $\tau_{\rho} = \tau_{\rho,m}$. Then, fixing $\rho=\rho^{\ell+1}$ and approximating $\rho^{\ell}\alpha^{\ell}$ by $\rho^{\ell+1}\alpha^{\ell}$, the update for $\alpha$ is
\begin{equation*}
\begin{adjustbox}{width=0.99\textwidth}$
\begin{split}
\alpha^{\ell+1} &= \argmax_{\alpha} \left\{\mathcal{L}(\tilde\varphi^{\ell}, \rho^{\ell+1}, \alpha) -  \frac{1}{2\tau_{\rho,m}}\|\rho^{\ell+1}\alpha - \rho^{\ell+1}\alpha^{\ell}\|_{L^2}^2\right\}
\\
&=\argmax_{\alpha}\int_0^T \int_\Omega -\rho^{\ell+1}(x,t) \langle A_{x,t}\alpha(x,t), \nabla_x \varphi^\ell(x,t)\rangle
- \rho^{\ell+1}(x,t)L_{x,t}(\alpha(x,t)) \\
&\quad\quad\quad\quad\quad\quad\quad\quad - \frac{\rho^{\ell+1}(x,t)^2}{2\tau_{\rho,m}}(\alpha(x,t) - \alpha^{\ell}(x,t))^2dxdt.
\end{split}
$\end{adjustbox}
\end{equation*}
This optimization problem can be solved in parallel for any $(x,t)$ by updating $\alpha(x,t)$ as follows
\begin{equation*}
\begin{adjustbox}{width=0.99\textwidth}$
\begin{split}
\alpha^{\ell+1}(x,t) &= \argmax_{\alpha\in\R^m}\left\{-\rho^{\ell+1}(x,t) \langle A_{x,t}\alpha, \nabla_x \varphi^\ell(x,t)\rangle
- \rho^{\ell+1}(x,t)L_{x,t}(\alpha) - \frac{\rho^{\ell+1}(x,t)^2}{2\tau_{\rho,m}}(\alpha - \alpha^{\ell}(x,t))^2\right\}\\
&= \argmin_{\alpha\in\R^m}\left\{\langle A_{x,t}\alpha, \nabla_x \varphi^\ell(x,t)\rangle
+L_{x,t}(\alpha) + \frac{\rho^{\ell+1}(x,t)}{2\tau_{\rho,m}}(\alpha - \alpha^{\ell}(x,t))^2\right\},
\end{split}
$\end{adjustbox}
\end{equation*}
which gives the second line in~\eqref{eqt:pdhg_det_cont} with $\tau_{\alpha} = \tau_{\rho,m}$. This also provides an intuition for choosing the specific penalty for $\alpha$.
\end{rem}

\begin{rem}
[PDHG convergence for discretized problems]
Note that it is not easy to apply the analysis above to the discretized saddle point problem in Section~\ref{sec:discretization_det}. The main difficulty is, after spatial discretization, the function $K$ is no longer bilinear and smooth, which requires PDHG analysis for non-linear and non-smooth coupling term. Moreover, how to analyze the convergence of the proposed algorithm for more general $f$ and $L$ (for instance, when $f$ is not affine with respect to $\alpha$) requires more study.
\end{rem}

\section{Consistency of the numerical Hamiltonian}\label{sec:appendix_consistency}

In Appendix~\ref{appendix:conv}, we explore how the numerical algorithm converges to a saddle point, necessitating a subsequent examination of how this saddle point relates to the solution of the HJ PDE. The linkage between the saddle point and the HJ PDE solution becomes apparent through the first-order optimality conditions outlined in~\eqref{eqt:first_order_optimality_cont}, when the stationary point \(\rho\) is positive.
Furthermore, by referencing the optimality conditions specified for discretized scenarios in~\eqref{eqt:first_order_optimality_semi1d} and~\eqref{eqt:first_order_optimality_semi2d}, we derive a discretized formulation for solving HJ PDEs, alongside the associated numerical Hamiltonians in~\eqref{eqt:numericalH1d} and~\eqref{eqt:numericalH2d}.
This section aims to validate the consistency of the numerical Hamiltonians, particularly in scenarios where the function \(f\) exhibits linearity in relation to \(\alpha\), and under a set of additional assumptions.

\subsection{One-dimensional cases} \label{sec:appendix_consistency_1d}
In this section, we consider the one-dimensional cases ($n=1$). Assume $f$ is linear with respect to $\alpha$, i.e., $f(x,t,\alpha) = a(x,t)^T \alpha$ where $a$ is a continuous function from $\Omega\times [0,T]$ to $\R^m$. Assume $L_{x,t}\colon \R^m\to\R$ is non-negative, convex, and $1$-coercive with $L_{x,t}(0) = 0$ for any $x\in \Omega$ and $t\in [0,T]$. In this case, we define the numerical Lagrangian by $\hat L_{x,t}(\alpha_1,\alpha_2) = L_{x,t}(\alpha_1) + L_{x,t}(\alpha_2)$. We are going to show the consistency of the corresponding numerical Hamiltonian. 
The corresponding numerical Hamiltonian $\hat H$ is 
\begin{equation*}
\begin{adjustbox}{width=0.99\textwidth}$
\begin{split}
\hat H(x,t,p^+, p^-) &= \sup_{\alpha_1, \alpha_2\in \R^m} \{-(a(x,t)^T\alpha_{1})_+ p^+ - (a(x,t)^T\alpha_{2})_- p^- - L_{x,t}(\alpha_{1})- L_{x,t}(\alpha_{2})\}\\
&= \sup_{\alpha_1\in \R^m} \{-(a(x,t)^T\alpha_{1})_+ p^+ - L_{x,t}(\alpha_{1})\} + \sup_{\alpha_2\in \R^m} \{-(a(x,t)^T\alpha_{2})_- p^- - L_{x,t}(\alpha_{2})\}.
\end{split}
$\end{adjustbox}
\end{equation*}
To show the consistency of the numerical Hamiltonian, we need to show that $\hat H(x,t,p, p) = H(x,t,p)$ for any $x\in \Omega$, $t\in [0,T]$ and $p\in \R$ where $H$ is the Hamiltonian defined in~\eqref{eqt:def_H}. 

Let $x\in \Omega$, $t\in [0,T]$ and $p\in\R$ be arbitrary numbers. To simplify notations, we ignore the $x,t$ dependency in $a$ and $L$, and denote $a(x,t)$ by $a$ and $L_{x,t}$ by $L$. 
Define $h_+$ by $h_+ = \sup_{\alpha\in \R^m} \{-(a^T\alpha)_+ p - L(\alpha)\}$ and $h_-$ by $h_- = \sup_{\alpha\in \R^m} \{-(a^T\alpha)_+ p - L(\alpha)\}$. Then, we need to show $h_+ + h_- = H(x,t,p)$.
Since $0 = L(0) = \min_{\alpha\in \R^m}L(\alpha)$, we have
\begin{equation}\label{eqt:appendix_def_hp_hm_1d}
h_+ = \sup_{\alpha: a^T\alpha \geq 0} \{-p a^T\alpha - L(\alpha)\}, \quad\quad h_- = \sup_{\alpha: a^T\alpha \leq 0} \{-p a^T\alpha - L(\alpha)\}.
\end{equation}
Since $L$ is convex and $1$-coercive, then the optimizer in~\eqref{eqt:def_H} exists. Denote the set of optimizers by $\mathcal{A}$, which is a closed convex set. Define $\mathcal{A}_+$ and $\mathcal{A}_-$ by $\mathcal{A}_+ = \{\alpha\in \mathcal{A}\colon a^T\alpha \geq 0\}$ and $\mathcal{A}_- = \{\alpha\in \mathcal{A}\colon a^T\alpha \leq 0\}$. 

If $\mathcal{A}_+$ is non-empty, then we have $h_+ = H(x,t,p)$. If $\mathcal{A}_+$ is empty, we show that $h_+$ equals zero. Denote the objective function in the minimization problem in~\eqref{eqt:def_H} by $F$, i.e., $F(\alpha) = -p a^T\alpha - L(\alpha)$. Note that $F$ is also the objective functions in the optimization problems in~\eqref{eqt:appendix_def_hp_hm_1d}. Since $L$ is convex and $1$-coercive, the minimizer in the first optimization problem in~\eqref{eqt:appendix_def_hp_hm_1d} exists. If there is a minimizer $\bar\alpha$ satisfying $a^T\bar\alpha > 0$, then there is a neighborhood of $\bar\alpha$, denoted by $\mathcal{N}$, such that $\bar\alpha$ is the minimizer of the function $F$ in $\mathcal{N}$. In other words, $\bar\alpha$ is a local minimizer of $F$, and hence it is also a global minimizer by convexity. 
Therefore, we have $\bar\alpha \in \mathcal{A}_+$, which contradicts to the assumption that $\mathcal{A}_+$ is empty. As a result, we have $h_+ = \sup_{\alpha: a^T\alpha = 0} \{-p a^T\alpha - L(\alpha)\} = 0$ in the case when $\mathcal{A}_+$ is empty. Then, we obtain
\begin{equation*}
h_+ = \begin{dcases}
H(x,t,p) & \mathcal{A}_+ \neq \emptyset, \\
0 & \mathcal{A}_+ = \emptyset.
\end{dcases}
\end{equation*}
With a similar argument, we can get
\begin{equation*}
h_- = \begin{dcases}
H(x,t,p) & \mathcal{A}_- \neq \emptyset, \\
0 & \mathcal{A}_- = \emptyset.
\end{dcases}
\end{equation*}

If either $\mathcal{A}_+$ or $\mathcal{A}_-$ is empty, then $h_+ + h_- = H(x,t,p)$ follows directly from these two formulas.
If both $\mathcal{A}_+$ and $\mathcal{A}_-$ are non-empty, by convexity of $\mathcal{A}$, there exists $\bar\alpha\in \mathcal{A}$ satisfying $a^T\bar\alpha = 0$, which implies $H(x,t,p) = \sup_{\alpha: a^T\alpha = 0} \{-p a^T\alpha - L(\alpha)\} = 0$, and hence $h_+ + h_- = 2H(x,t,p) = H(x,t,p)$. Therefore, we conclude that $\hat H$ is consistent. 

\begin{rem}\label{rem:A1}
The above argument can be applied to more general cases. Let $f$ be in the form of $f(x,t,\alpha) = a(x,t)^T \alpha + b(x,t)$, where $a$ is a continuous function from $\Omega\times [0,T]$ to $\R^m$ and $b$ is a continuous function from $\Omega\times [0,T]$ to $\R$.
Let $L$ be a continuous function such that $L_{x,t}$ is convex and $1$-coercive for any $x\in \Omega$ and $t\in [0,T]$. Further assume that $\argmin_{\alpha\in \R^m} L_{x,t}(\alpha) \cap \{\alpha \in \R^m \colon f(x,t,\alpha) = 0\}$ is non-empty for any $x\in \Omega$ and $t\in [0,T]$. Then, define the numerical Lagrangian $\hat L$ by $\hat L_{x,t}(\alpha_1, \alpha_2) = L_{x,t}(\alpha_1) + L_{x,t}(\alpha_2) - \min_{\alpha\in \R^m} L_{x,t}(\alpha)$. With a similar argument as above, we can prove that the corresponding numerical Hamiltonian defined by~\eqref{eqt:numericalH1d} is consistent.
\end{rem}

\begin{rem}
We can further remove the requirement in Remark~\ref{rem:A1} that the minimal value of $L_{x,t}$ can be achieved on a root of $f_{x,t}$. In other words, we only need to assume the continuity of $f$ and $L$, and assume that $f$ is affine on $\alpha$, and $L_{x,t}$ is convex and $1$-coercive for any $x\in \Omega$ and $t\in [0,T]$. Note that the $1$-coercivity assumption of $L_{x,t}$ is for the existence of the optimizer in~\eqref{eqt:def_H}. It can be replaced by any other assumption which guarantees the existence of the minimizers. The assumption that $f$ is affine on $\alpha$ guarantees the convexity of the optimization problem in~\eqref{eqt:def_H} so that we can use convex analysis techniques in the proof. Under this general setup, we need to modify the saddle point problem~\eqref{eqt:saddle_point_det_1d_semi} to 
\begin{equation*}
\begin{adjustbox}{width=0.99\textwidth}$
\begin{split}
\min_{\substack{\varphi\\ \varphi_i(0)=g(x_i)}}\max_{\substack{\rho\geq 0, \alpha\\
f_{i,t}(\alpha_{1,i}(t))\geq 0 \forall i,t\\
f_{i,t}(\alpha_{2,i}(t))\leq 0 \forall i,t}}  \int_0^T \sum_{i=1}^{n_x} \rho_i(t)\Bigg(\dot\varphi_i(t) - f_{i,t}(\alpha_{1,i}(t)) (D_x^+\varphi)_i(t) - f_{i,t}(\alpha_{2,i}(t)) (D_x^-\varphi)_i(t) \\
- \hat L_{i,t}\left(\alpha_{1,i}(t), \alpha_{2,i}(t)\right) \Bigg)dt - c\sum_{i=1}^{n_x} \varphi_i(T).
\end{split}
$\end{adjustbox}
\end{equation*}
Then, we define the numerical Lagrangian by
\begin{equation*}
\hat L_{x,t}(\alpha_1, \alpha_2) = L_{x,t}(\alpha_1) + L_{x,t}(\alpha_2) - \min_{\substack{\alpha\in \R^m\\f(x,t,\alpha) = 0}} L_{x,t}(\alpha).
\end{equation*}
The consistency of the corresponding numerical Hamiltonian defined by
\begin{equation*}
\hat H(x,t,p^+, p^-) = \sup_{\substack{\alpha_1, \alpha_2\in \R^m\\ f_{x,t}(\alpha_{1})\geq 0\\
f_{x,t}(\alpha_{2})\leq 0}} \{-f_{x,t}(\alpha_{1}) p^+ - f_{x,t}(\alpha_{2}) p^- - \hat L_{x,t}(\alpha_1, \alpha_2)\}
\end{equation*}
can be proved similarly as above using the formulas
\begin{equation*}
\begin{split}
\sup_{\substack{\alpha\in \R^m\\ f_{x,t}(\alpha)\geq 0}} \{-f_{x,t}(\alpha) p - L_{x,t}(\alpha)\} &= \begin{dcases}
H(x,t,p) & \mathcal{A}_{+,x,t}\neq \emptyset, \\
-\min_{\substack{\alpha\in \R^m\\f(x,t,\alpha) = 0}} L_{x,t}(\alpha) & \mathcal{A}_{+,x,t}= \emptyset,
\end{dcases}\\
\sup_{\substack{\alpha\in \R^m\\ f_{x,t}(\alpha)\leq 0}} \{-f_{x,t}(\alpha) p - L_{x,t}(\alpha)\} &= \begin{dcases}
H(x,t,p) & \mathcal{A}_{-,x,t}\neq \emptyset, \\
-\min_{\substack{\alpha\in \R^m\\f(x,t,\alpha) = 0}} L_{x,t}(\alpha) & \mathcal{A}_{-,x,t}= \emptyset,
\end{dcases}
\end{split}
\end{equation*}
where $\mathcal{A}_{+,x,t}$ ($\mathcal{A}_{-,x,t}$ resp.) are the set containing the minimizers $\bar\alpha$ in~\eqref{eqt:def_H} which satisfy $f_{x,t}(\bar\alpha)\geq 0$ ($f_{x,t}(\bar\alpha)\leq 0$ resp.). 
Then, the implicit time discretization and PDHG updates can also be applied to this modified saddle point problem. 
\end{rem}

\subsection{Two-dimensional cases}\label{sec:appendix_consistency_2d}
In this section, we consider the case when the spatial domain is two-dimensional ($n=2$). 
We assume that $f$ is a continuous function in the form $f(x,t,\alpha) = \begin{pmatrix} a_1^T(x,t)\alpha \\ a_2^T(x,t) \alpha\end{pmatrix}$, where $a_1$ and $a_2$ are continuous functions from $\Omega\times [0,T]$ to $\R^{m}$. 
Let $L$ be a continuous non-negative function. Assume $L_{x,t}\colon \R^m\to\R$ is convex and $1$-coercive with $L_{x,t}(0) = 0$ for any $x\in \Omega$ and $t\in [0,T]$.
Let $H$ be the Hamiltonian defined by~\eqref{eqt:def_H}.
Assume that $H$ is in the form of $H_{x,t}(p_x,p_y) = H_{1}(x,t,p_x) + H_{2}(x,t,p_y)$ for any $x\in\Omega$, $t\in [0,T]$, $p_x,p_y\in\R$ for some functions $H_1$ and $H_2$. In this case, we choose $\hat L(x,t,\alpha_{11}, \alpha_{12},\alpha_{21}, \alpha_{22}) = L(x,t,\alpha_{11}) + L(x,t,\alpha_{12}) + L(x,t,\alpha_{21}) + L(x,t,\alpha_{22})$. Now, we want to show that such defined numerical Lagrangian gives a consistent numerical Hamiltonian.

Let $x\in \Omega$, $t\in [0,T]$, and $p_x, p_y\in \R$ be arbitrary vectors and numbers. For simplicity of the notations, we ignore the $(x,t)$ dependence in the functions $a_1$, $a_2$, and $L$. 
Define $h_{x,+}$, $h_{x,-}$, $h_{y,+}$, and $h_{y,-}$ by 
\begin{equation*}
\begin{split}
h_{x,+} = \sup_{\alpha: a_1^T\alpha \geq 0} \{-p_x a_1^T\alpha - L(\alpha)\}, \quad h_{x,-} = \sup_{\alpha: a_1^T\alpha \leq 0} \{-p_x a_1^T\alpha - L(\alpha)\},\\
h_{y,+} = \sup_{\alpha: a_2^T\alpha \geq 0} \{-p_y a_2^T\alpha - L(\alpha)\}, \quad h_{y,-} = \sup_{\alpha: a_2^T\alpha \leq 0} \{-p_y a_2^T\alpha - L(\alpha)\}. 
\end{split}
\end{equation*}
Since $L_{x,t}$ is a non-negative function with $L_{x,t}(0) = 0$, the numerical Hamiltonian defined in~\eqref{eqt:numericalH2d} satisfies $\hat H(x,t,p_x,p_x,p_y,p_y) = h_{x,+} + h_{x,-} + h_{y,+} + h_{y,-}$. Then, it suffices to prove that $h_{x,+} + h_{x,-} + h_{y,+} + h_{y,-} = H(x,t,p_x,p_y)$.
Denote the set of optimizers in~\eqref{eqt:def_H} by $\mathcal{A}(p)$, which, according to the assumptions on $f$ and $L$, is non-empty, closed, and convex for any $p\in\R^2$. With a similar argument as in the one-dimensional case, we have
\begin{equation*}
\begin{split}
h_{x,+} = 
\begin{dcases}
H(x,t,p_x,0) & \mathcal{A}(p_x,0)\cap \{\alpha\in \R^m: a_1^T\alpha \geq 0\} \neq \emptyset,\\
0 & \mathcal{A}(p_x,0)\cap \{\alpha\in \R^m: a_1^T\alpha \geq 0\} = \emptyset,
\end{dcases}\\
h_{x,-} = 
\begin{dcases}
H(x,t,p_x,0) & \mathcal{A}(p_x,0)\cap \{\alpha\in \R^m: a_1^T\alpha \leq 0\} \neq \emptyset,\\
0 & \mathcal{A}(p_x,0)\cap \{\alpha\in \R^m: a_1^T\alpha \leq 0\} = \emptyset,
\end{dcases}\\
h_{y,+} = 
\begin{dcases}
H(x,t,0,p_y) & \mathcal{A}(0,p_y)\cap \{\alpha\in \R^m: a_2^T\alpha \geq 0\} \neq \emptyset,\\
0 & \mathcal{A}(0,p_y)\cap \{\alpha\in \R^m: a_2^T\alpha \geq 0\} = \emptyset,
\end{dcases}\\
h_{y,-} = 
\begin{dcases}
H(x,t,0,p_y) & \mathcal{A}(0,p_y)\cap \{\alpha\in \R^m: a_2^T\alpha \leq 0\} \neq \emptyset,\\
0 & \mathcal{A}(0,p_y)\cap \{\alpha\in \R^m: a_2^T\alpha \leq 0\} = \emptyset.
\end{dcases}
\end{split}
\end{equation*}

If either $\mathcal{A}(p_x,0)\cap \{\alpha\in \R^m: a_1^T\alpha \geq 0\}$ or $\mathcal{A}(p_x,0)\cap \{\alpha\in \R^m: a_1^T\alpha \leq 0\}$ is empty, then we have $h_{x,+} + h_{x,-} = H(x,t,p_x,0)$. Otherwise, similarly as in the one-dimensional case, there exists $\bar\alpha\in \mathcal{A}(p_x,0)$ satisfying $a_1^T\bar\alpha = 0$, which implies $H(x,t,p_x,0) = 0$, and hence $h_{x,+} + h_{x,-} = 2H(x,t,p_x,0) = H(x,t,p_x,0)$ holds. Similarly, we have $h_{y,+} + h_{y,-} = H(x,t,0,p_y)$. Therefore, we get
\begin{equation*}
\begin{adjustbox}{width=0.99\textwidth}$
\begin{split}
h_{x,+} + h_{x,-} + h_{y,+} + h_{y,-} &= H(x,t,p_x,0) + H(x,t,0,p_y)\\
&= H_1(x,t,p_x)+ H_2(x,t,0) + H_1(x,t,0)+ H_2(x,t,p_y) \\
&= H(x,t,p_x,p_y) + H(x,t,0,0) \\
&= H(x,t,p_x,p_y) - \min_{\alpha\in\R^m} L_{x,t}(\alpha) = H(x,t,p_x,p_y),
\end{split}
$\end{adjustbox}
\end{equation*}
which ends the proof.

\begin{rem}
Similarly as in the one-dimensional case, the assumptions can be further relaxed. Assume that $f$ and $L$ are both continuous, $f_{x,t}$ is an affine function, and $L_{x,t}$ is convex and $1$-coercive for any $x\in \Omega$ and $t\in [0,T]$. Note that the $1$-coercivity assumption of $L_{x,t}$ is for the existence of the optimizer in~\eqref{eqt:def_H}, and the assumption that $f_{x,t}$ is affine makes the optimization problem in~\eqref{eqt:def_H} a convex problem. Under this general setup, we modify the saddle point problem~\eqref{eqt:saddle_point_det_2d_semi} to
\begin{equation*}
\begin{adjustbox}{width=0.99\textwidth}$
\begin{split}
\min_{\substack{\varphi\\ \varphi_{i,j}(0)=g(x_{i,j})}}\max_{\substack{\rho\geq 0, \alpha\\f_{1,i,j,t}(\alpha_{11,i,j}(t))\geq 0\forall i,j,t\\f_{1,i,j,t}(\alpha_{12,i,j}(t))\leq 0\forall i,j,t\\f_{2,i,j,t}(\alpha_{21,i,j}(t))\geq 0\forall i,j,t\\f_{2,i,j,t}(\alpha_{22,i,j}(t))\leq 0\forall i,j,t}}  \int_0^T \sum_{i=1}^{n_x} \sum_{j=1}^{n_y} \rho_{i,j}(t)\Bigg(\dot\varphi_{i,j}(t) - f_{1,i,j,t}(\alpha_{11,i,j}(t)) (D_x^+\varphi)_i(t) 
- f_{1,i,j,t}(\alpha_{12,i,j}(t)) (D_x^-\varphi)_i(t) \\
- f_{2,i,j,t}(\alpha_{21,i,j}(t)) (D_y^+\varphi)_i(t) - f_{2,i,j,t}(\alpha_{22,i,j}(t)) (D_y^-\varphi)_i(t) \\
- \hat L_{i,j,t}\left(\alpha_{11,i,j}(t), \alpha_{12,i,j}(t), \alpha_{21,i,j}(t), \alpha_{22,i,j}(t)\right) \Bigg)dt - c\sum_{i=1}^{n_x}\sum_{j=1}^{n_y} \varphi_{i,j}(T).
\end{split} 
$\end{adjustbox}
\end{equation*}
We define the numerical Lagrangian by
\begin{equation*}
\begin{adjustbox}{width=0.99\textwidth}$
\begin{split}
\hat L_{x,t}(\alpha_{11}, \alpha_{12}, \alpha_{21}, \alpha_{22}) = L_{x,t}(\alpha_{11}) + L_{x,t}(\alpha_{12}) + L_{x,t}(\alpha_{21}) + L_{x,t}(\alpha_{22}) - \min_{\alpha\in \R^m} L_{x,t}(\alpha)\\ - \min_{\substack{\alpha\in \R^m\\f_1(x,t,\alpha) = 0}} L_{x,t}(\alpha) - \min_{\substack{\alpha\in \R^m\\f_2(x,t,\alpha) = 0}} L_{x,t}(\alpha).
\end{split}
$\end{adjustbox}
\end{equation*}
The numerical Hamiltonian is then given by
\begin{equation*}
\begin{adjustbox}{width=0.99\textwidth}$
\begin{split}
\hat H(x,t,p_{x}^+, p_{x}^-, p_{y}^+, p_{y}^-) = \sup_{\substack{\alpha_{11}, \alpha_{12},\alpha_{21}, \alpha_{22}\in \R^m\\ f_{1,x,t}(\alpha_{11})\geq 0\\
f_{1,x,t}(\alpha_{12})\leq 0\\f_{2,x,t}(\alpha_{21})\geq 0\\
f_{2,x,t}(\alpha_{22})\leq 0}} \{-f_{1,x,t}(\alpha_{11}) p_x^+ - f_{1,x,t}(\alpha_{12}) p_x^--f_{2,x,t}(\alpha_{21}) p_y^+ \\
- f_{2,x,t}(\alpha_{22}) p_y^- - \hat L_{x,t}(\alpha_{11}, \alpha_{12},\alpha_{21}, \alpha_{22})\}.
\end{split}
$\end{adjustbox}
\end{equation*}
The consistency of this numerical Hamiltonian follows from a similar argument. Then, the implicit time discretization and PDHG updates can be similarly applied to this modified saddle point problem.
\end{rem}

\begin{rem}
Note that the assumption in the two-dimensional cases is more restricted than the one-dimensional cases. To be more specific, we assume that there exists $H_1$ and $H_2$ such that $H$ is in the form of $H_{x,t}(p_x,p_y) = H_{1}(x,t,p_x) + H_{2}(x,t,p_y)$ for any $x\in\Omega$, $t\in [0,T]$, $p_x,p_y\in\R$. If we consider the numerical Lagrangian $\hat L$ which can be written as $\hat L_{x,t}(\alpha_{11}, \alpha_{12}, \alpha_{21}, \alpha_{22}) = L_{11,x,t}(\alpha_{11}) + L_{12,x,t}(\alpha_{12}) + L_{21,x,t}(\alpha_{21}) + L_{22,x,t}(\alpha_{22})$ for some functions $L_{11}$, $L_{12}$, $L_{21}$, and $L_{22}$, then the corresponding numerical Hamiltonian is
\begin{equation*}
\begin{adjustbox}{width=0.99\textwidth}$
\begin{split}
\hat H(x,t,p_{x}^+, p_{x}^-, p_{y}^+, p_{y}^-) &= \sup_{\substack{\alpha_{11}, \alpha_{12},\alpha_{21}, \alpha_{22}\in \R^m\\ f_{1,x,t}(\alpha_{11})\geq 0\\
f_{1,x,t}(\alpha_{12})\leq 0\\f_{2,x,t}(\alpha_{21})\geq 0\\
f_{2,x,t}(\alpha_{22})\leq 0}} \{-f_{1,x,t}(\alpha_{11}) p_x^+ - f_{1,x,t}(\alpha_{12}) p_x^--f_{2,x,t}(\alpha_{21}) p_y^+ - f_{2,x,t}(\alpha_{22}) p_y^- \\
&\quad\quad\quad\quad\quad\quad\quad\quad - (L_{11,x,t}(\alpha_{11}) + L_{12,x,t}(\alpha_{12}) + L_{21,x,t}(\alpha_{21}) + L_{22,x,t}(\alpha_{22}))\}\\
&= \sup_{\substack{\alpha\in \R^m\\ f_{1,x,t}(\alpha)\geq 0}} \{-f_{1,x,t}(\alpha) p_x^+ - L_{11,x,t}(\alpha)\} + \sup_{\substack{\alpha\in \R^m\\ f_{1,x,t}(\alpha)\leq 0}} \{-f_{1,x,t}(\alpha) p_x^- - L_{12,x,t}(\alpha)\}\\
&\quad\quad 
+ \sup_{\substack{\alpha\in \R^m\\ f_{2,x,t}(\alpha)\geq 0}} \{-f_{2,x,t}(\alpha) p_y^+ - L_{21,x,t}(\alpha)\} + \sup_{\substack{\alpha\in \R^m\\ f_{2,x,t}(\alpha)\leq 0}} \{-f_{2,x,t}(\alpha) p_y^- - L_{22,x,t}(\alpha)\}.
\end{split}
$\end{adjustbox}
\end{equation*}
Then, if $\hat H$ is consistent, $H$ must be in the form of $H_{x,t}(p_x,p_y) = H_{1}(x,t,p_x) + H_{2}(x,t,p_y)$ with
\begin{equation*}
\begin{adjustbox}{width=0.99\textwidth}$
\begin{split}
H_1(x,t,p_x) &= \sup_{\substack{\alpha\in \R^m\\ f_{1,x,t}(\alpha)\geq 0}} \{-f_{1,x,t}(\alpha) p_x - L_{11,x,t}(\alpha)\} + \sup_{\substack{\alpha\in \R^m\\ f_{1,x,t}(\alpha)\leq 0}} \{-f_{1,x,t}(\alpha) p_x - L_{12,x,t}(\alpha)\},\\
H_2(x,t,p_y) &= \sup_{\substack{\alpha\in \R^m\\ f_{2,x,t}(\alpha)\geq 0}} \{-f_{2,x,t}(\alpha) p_y - L_{21,x,t}(\alpha)\} + \sup_{\substack{\alpha\in \R^m\\ f_{2,x,t}(\alpha)\leq 0}} \{-f_{2,x,t}(\alpha) p_y - L_{22,x,t}(\alpha)\}.
\end{split}
$\end{adjustbox}
\end{equation*}

For more general $H$, we need to design the numerical Lagrangian $\hat L$ in other forms. We also considered $\hat L_{x,t}(\alpha_{11}, \alpha_{12}, \alpha_{21}, \alpha_{22}) = L_{x,t}(\alpha_{11}+\alpha_{12}+\alpha_{21}+\alpha_{22})$, but it has too many degrees of freedom on $\alpha$. To illustrate this intuition, we consider the case when we have $m=n=2$ and $f(x,t,\alpha) = \alpha$. Denote $\alpha_{IJ} = (\alpha_{IJ,1}, \alpha_{IJ,2})\in\R^2$ for any $I,J=1,2$. Then, the numerical Hamiltonian is
\begin{equation*}
\begin{adjustbox}{width=0.99\textwidth}$
\begin{split}
\hat H(x,t,p_{x}, p_{x}, p_{y}, p_{y}) &= \sup_{\substack{\alpha_{11}, \alpha_{12},\alpha_{21}, \alpha_{22}\in \R^2\\ \alpha_{11,1}\geq 0\\
\alpha_{12,1}\leq 0\\
\alpha_{21,2}\geq 0\\
\alpha_{22,2}\leq 0}} \{-\alpha_{11,1} p_x - \alpha_{12,1} p_x-\alpha_{21,2} p_y - \alpha_{22,2} p_y - L_{x,t}(\alpha_{11}+\alpha_{12}+\alpha_{21}+\alpha_{22})\}.
\end{split}
$\end{adjustbox}
\end{equation*}
For any $\alpha_{11,1}, \alpha_{21,2}\geq 0$ and $\alpha_{12,1}, \alpha_{22,2}\leq 0$, we can always choose $\alpha_{11,2}$, $\alpha_{12,2}$, $\alpha_{21,1}$, $\alpha_{22,1}$ such that $L_{x,t}(\alpha_{11}+\alpha_{12}+\alpha_{21}+\alpha_{22})$ achieves the minimal value of $L_{x,t}$. Therefore, the numerical Hamiltonian $\hat H$ is infinity if $p_x$ is non-zero or $p_y$ is non-zero, which is not consistent. How to choose a numerical Lagrangian for more general cases is a possible interesting future direction.
\end{rem}

\section{Discretization for the viscous HJ PDEs}\label{sec:discretization_soc}
In this section, we provide details of the spatial and temporal discretization of the saddle point problem~\eqref{eqt:saddle_soc} for solving stochastic optimal control problems and the corresponding viscous HJ PDEs. We apply similar notation as in Section~\ref{sec:discretization_det}. Note that the difference between this part with Section~\ref{sec:discretization_det} is merely on the diffusion term.
\subsection{One-dimensional problems}
After upwind spatial discretization, the one-dimensional saddle point problem (divided by $\Delta x$) becomes
\begin{equation*}
\begin{adjustbox}{width=0.99\textwidth}$
\begin{split}
\min_{\substack{\varphi\\ \varphi_i(0)=g(x_i)}}\max_{\rho\geq 0, \alpha}  \int_0^T \sum_{i=1}^{n_x} \rho_i(t)\Bigg(\dot\varphi_i(t) - f_{i,t}(\alpha_{1,i}(t))_+ (D_x^+\varphi)_i(t) - f_{i,t}(\alpha_{2,i}(t))_- (D_x^-\varphi)_i(t) \\
- \hat L_{i,t}\left(\alpha_{1,i}(t), \alpha_{2,i}(t)\right) -\epsilon (D_{xx}\varphi)_i(t) \Bigg)dt - c\sum_{i=1}^{n_x} \varphi_i(T).
\end{split}
$\end{adjustbox}
\end{equation*}
Consider a stationary point $(\varphi, \rho, \alpha)$ in this saddle point problem. If we further assume $\rho_i(t) > 0$ for any $t\in [0,T]$, then the first order optimality condition is
\begin{equation*}
\begin{adjustbox}{width=0.99\textwidth}$
\begin{dcases}
\dot{\varphi}_i(t) - f_{i,t}(\alpha_{1,i}(t))_+ (D_x^+ \varphi)_i(t) - f_{i,t}(\alpha_{2,i}(t))_- (D_x^- \varphi)_i(t) - \hat L_{i,t}(\alpha_{1,i}(t), \alpha_{2,i}(t)) -\epsilon (D_{xx}\varphi)_i(t) = 0,\\
(\alpha_{1,i}(t), \alpha_{2,i}(t)) = \argmin_{\alpha_1, \alpha_2\in\R^m}\{f_{i,t}(\alpha_1)_+ (D_x^+\varphi)_i(t) + f_{i,t}(\alpha_{2})_- (D_x^- \varphi)_i(t) + \hat L_{i,t}(\alpha_1, \alpha_2)\}, \\
\dot{\rho}_i(t) - D_x^-(f_{i,t}(\alpha_{1,i}(t))_+\rho_i(t)) - D_x^+(f_{i,t}(\alpha_{2,i}(t))_-\rho_i(t)) + \epsilon (D_{xx}\rho)_i(t) = 0.
\end{dcases}
$\end{adjustbox}
\end{equation*}
Combining the first two lines in this optimality condition, we get
\begin{equation*}
\begin{adjustbox}{width=0.99\textwidth}$
\dot{\varphi}_i(t) +\sup_{\alpha_1,\alpha_2\in \R^m} \{-f_{i,t}(\alpha_{1})_+ (D_x^+ \varphi)_i(t) - f_{i,t}(\alpha_{2})_- (D_x^- \varphi)_i(t) - \hat L_{i,t}(\alpha_{1}, \alpha_{2})\} -\epsilon (D_{xx}\varphi)_i(t) = 0,
$\end{adjustbox}
\end{equation*}
which gives a semi-discrete scheme for the HJ PDE~\eqref{eqt:visc_HJ_initial} 
where the numerical Hamiltonian is defined by~\eqref{eqt:numericalH1d}.

With this discretization, the $\ell$-th update becomes
\begin{equation*}
{\scriptsize
\begin{dcases}
\rho_i^{\ell+1}(t) = \Big(\rho_i^{\ell}(t) + \tau_\rho\big(\dot{\tilde\varphi}^\ell_i(t) - f_{i,t}(\alpha^{\ell}_{1,i}(t))_+ (D_x^+ \tilde\varphi^\ell)_i(t) - f_{i,t}(\alpha^\ell_{2,i}(t))_- (D_x^- \tilde\varphi^\ell)_i(t) \\
\quad\quad\quad\quad\quad\quad\quad\quad
- \hat L_{i,t}(\alpha^\ell_{1,i}(t), \alpha^\ell_{2,i}(t)) -\epsilon (D_{xx}\tilde\varphi^\ell)_i(t)\big)\Big)_+.\\
(\alpha^{\ell+1}_{1,i}(t), \alpha^{\ell+1}_{2,i}(t)) = \argmin_{\alpha_1, \alpha_2\in\R^m}\Big\{f_{i,t}(\alpha_1)_+ (D_x^+\tilde\varphi^\ell)_i(t) + f_{i,t}(\alpha_2)_- (D_x^-\tilde\varphi^\ell)_i(t) + \hat L_{i,t}(\alpha_1, \alpha_2)\\
\quad\quad\quad\quad\quad\quad\quad\quad + \frac{\rho_i^{\ell+1}(t)}{2\tau_\alpha} \left(\|\alpha_1 - \alpha^\ell_{1,i}(t)\|^2 +  \|\alpha_2 - \alpha^\ell_{2,i}(t)\|^2\right)\Big\}.\\
\varphi^{\ell+1}_i(t) = \varphi^{\ell}_i(t) + \tau_{\varphi}(I-\partial_{tt}-D_{xx})^{-1}\Big(\dot{\rho}^{\ell+1}_i(t) - D_x^-(f_{i,t}(\alpha^{\ell+1}_{1,i}(t))_+\rho^{\ell+1}_i(t)) \\
\quad\quad\quad\quad\quad\quad\quad\quad - D_x^+(f_{i,t}(\alpha^{\ell+1}_{2,i}(t))_-\rho^{\ell+1}_i(t)) + \epsilon (D_{xx}\rho^{\ell+1})_i(t)\Big).\\
\tilde\varphi^{\ell+1} = 2\varphi^{\ell+1} - \varphi^\ell.
\end{dcases}
}
\end{equation*}

For the time discretization, we apply implicit Euler scheme for $\varphi$, and then the saddle point problem (divided by $\Delta t$) becomes
\begin{equation*}
{\scriptsize
\begin{split}
\min_{\substack{\varphi\\ \varphi_{i,1}=g(x_i)}}\max_{\rho\geq 0, \alpha}  \sum_{k=2}^{n_t} \sum_{i=1}^{n_x} \rho_{i,k}\Bigg((D_t^-\varphi)_{i,k} - f_{i,k}(\alpha_{1,i,k})_+ (D_x^+\varphi)_{i,k} - f_{i,k}(\alpha_{2,i,k})_- (D_x^-\varphi)_{i,k} \\
- \hat L_{i,k}\left(\alpha_{1,i,k}, \alpha_{2,i,k}\right) - \epsilon (D_{xx}\varphi)_{i,k} \Bigg) - \frac{c}{\Delta t}\sum_{i=1}^{n_x} \varphi_{i,n_t}.
\end{split}
}
\end{equation*}
The corresponding algorithm in the $\ell$-th iteration becomes
\begin{equation*}
{\scriptsize
\begin{dcases}
\rho_{i,k}^{\ell+1} = \Big(\rho_{i,k}^{\ell} + \tau_\rho\big((D_t^-{\tilde\varphi^\ell})_{i,k} - f_{i,k}(\alpha^{\ell}_{1,i,k})_+ (D_x^+ \tilde\varphi^\ell)_{i,k} - f_{i,k}(\alpha^\ell_{2,i,k})_- (D_x^- \tilde\varphi^\ell)_{i,k} \\
\quad\quad\quad\quad\quad\quad\quad\quad
- \hat L_{i,k}(\alpha^\ell_{1,i,k}, \alpha^\ell_{2,i,k}) - \epsilon (D_{xx}\tilde\varphi^\ell)_{i,k} \big)\Big)_+.\\
(\alpha^{\ell+1}_{1,i,k}, \alpha^{\ell+1}_{2,i,k}) = \argmin_{\alpha_1, \alpha_2\in\R^m}\Big\{f_{i,k}(\alpha_1)_+ (D_x^+\tilde\varphi^\ell)_{i,k} + f_{i,k}(\alpha_2)_- (D_x^-\tilde\varphi^\ell)_{i,k} + \hat L_{i,k}(\alpha_1, \alpha_2)\\
\quad\quad\quad\quad\quad\quad\quad\quad + \frac{\rho_{i,k}^{\ell+1}}{2\tau_\alpha} \left(\|\alpha_1 - \alpha^\ell_{1,i,k}\|^2 + \|\alpha_2 - \alpha^\ell_{2,i,k}\|^2\right)\Big\}.\\
\varphi^{\ell+1}_{i,k} = \varphi^{\ell}_{i,k} + \tau_{\varphi}(I-D_{tt}-D_{xx})^{-1}\Big((D_t^+{\rho}^{\ell+1})_{i,k} - D_x^-(f_{i,k}(\alpha^{\ell+1}_{1,i,k})_+\rho^{\ell+1}_{i,k}) \\
\quad\quad\quad\quad\quad\quad\quad\quad - D_x^+(f_{i,k}(\alpha^{\ell+1}_{2,i,k})_-\rho^{\ell+1}_{i,k}) + \epsilon (D_{xx}\rho^{\ell+1})_{i,k} \Big).\\
\tilde\varphi^{\ell+1} = 2\varphi^{\ell+1} - \varphi^\ell.
\end{dcases}
}
\end{equation*}
It's worth mentioning that selecting the numerical Lagrangian as $L_{x,t}(\alpha_{1}) + L_{x,t}(\alpha_{2})$ allows for the update process of $\alpha$ to be executed in parallel across the variables $\alpha_{1}$ and $\alpha_{2}$.

\subsection{Two-dimensional problems}
For two dimensional cases, with upwind spatial discretization, the saddle point problem (divided by $\Delta x \Delta y$) becomes
\begin{equation*}
\begin{adjustbox}{width=0.99\textwidth}$
\begin{split}
\min_{\substack{\varphi\\ \varphi_{i,j}(0)=g(x_{i,j})}}\max_{\rho\geq 0, \alpha}  \int_0^T \sum_{i=1}^{n_x} \sum_{j=1}^{n_y} \rho_{i,j}(t)\Bigg(\dot\varphi_{i,j}(t) - f_{1,i,j,t}(\alpha_{11,i,j}(t))_+ (D_x^+\varphi)_i(t) 
- f_{1,i,j,t}(\alpha_{12,i,j}(t))_- (D_x^-\varphi)_i(t) \\
- f_{2,i,j,t}(\alpha_{21,i,j}(t))_+ (D_y^+\varphi)_i(t) - f_{2,i,j,t}(\alpha_{22,i,j}(t))_- (D_y^-\varphi)_i(t) -\epsilon (D_{xx} \varphi + D_{yy} \varphi)_{i,j}(t) \\
- \hat L_{i,j,t}\left(\alpha_{11,i,j}(t), \alpha_{12,i,j}(t), \alpha_{21,i,j}(t), \alpha_{22,i,j}(t)\right) \Bigg)dt - c\sum_{i=1}^{n_x}\sum_{j=1}^{n_y} \varphi_{i,j}(T).
\end{split} 
$\end{adjustbox}
\end{equation*}
Consider a stationary point $(\varphi, \rho, \alpha)$ in this saddle point problem. If we further assume $\rho_{i,j}(t) > 0$ for any $t\in [0,T]$, then the first order optimality condition is
\begin{equation*}
\begin{adjustbox}{width=0.99\textwidth}$
\begin{dcases}
\dot{\varphi}_{i,j}(t) - f_{1,i,j,t}(\alpha_{11,i,j}(t))_+ (D_x^+ \varphi)_{i,j}(t) - f_{1,i,j,t}(\alpha_{12,i,j}(t))_- (D_x^- \varphi)_{i,j}(t) \\
\quad\quad - f_{2,i,j,t}(\alpha_{21,i,j}(t))_+ (D_y^+ \varphi)_{i,j}(t) - f_{2,i,j,t}(\alpha_{22,i,j}(t))_- (D_y^- \varphi)_{i,j}(t) \\
\quad\quad- \hat L_{i,j,t}(\alpha_{11,i,j}(t), \alpha_{12,i,j}(t), \alpha_{21,i,j}(t), \alpha_{22,i,j}(t)) -\epsilon (D_{xx} \varphi + D_{yy} \varphi)_{i,j}(t) = 0,\\
(\alpha_{11,i,j}(t), \alpha_{12,i,j}(t), \alpha_{21,i,j}(t), \alpha_{22,i,j}(t)) = \argmin_{\alpha_{11}, \alpha_{12},\alpha_{21}, \alpha_{22}\in\R^m}\{f_{1,i,j,t}(\alpha_{11})_+ (D_x^+\varphi)_{i,j}(t) \\
\quad\quad+ f_{1,i,j,t}(\alpha_{12})_- (D_x^- \varphi)_{i,j}(t) + f_{2,i,j,t}(\alpha_{21})_+ (D_y^+\varphi)_{i,j}(t) \\
\quad\quad + f_{2,i,j,t}(\alpha_{22})_- (D_y^- \varphi)_{i,j}(t) + \hat L_{i,j,t}(\alpha_{11}, \alpha_{12},\alpha_{21}, \alpha_{22})\}, \\
\dot{\rho}_{i,j}(t) - D_x^-(f_{1,i,j,t}(\alpha_{11,i,j}(t))_+\rho_{i,j}(t)) - D_x^+(f_{1,i,j,t}(\alpha_{12,i,j}(t))_-\rho_{i,j}(t))\\
\quad\quad - D_y^-(f_{2,i,j,t}(\alpha_{21,i,j}(t))_+\rho_{i,j}(t)) - D_y^+(f_{2,i,j,t}(\alpha_{22,i,j}(t))_-\rho_{i,j}(t)) + \epsilon (D_{xx} \rho + D_{yy} \rho)_{i,j}(t)= 0.
\end{dcases}
$\end{adjustbox}
\end{equation*}
Combining the first two lines in this optimality condition, we get
\begin{equation*}
\begin{adjustbox}{width=0.99\textwidth}$
\begin{split}
\dot{\varphi}_{i,j}(t) +\sup_{\alpha_{11},\alpha_{12},\alpha_{21},\alpha_{22}\in\R^m} \{-f_{1,i,j,t}(\alpha_{11})_+ (D_x^+ \varphi)_{i,j}(t) - f_{1,i,j,t}(\alpha_{12})_- (D_x^- \varphi)_{i,j}(t)
-f_{2,i,j,t}(\alpha_{21})_+ (D_y^+ \varphi)_{i,j}(t) \\
- f_{2,i,j,t}(\alpha_{22})_- (D_y^- \varphi)_{i,j}(t) 
- \hat L_{i,j,t}(\alpha_{11},\alpha_{12},\alpha_{21},\alpha_{22})\} - \epsilon (D_{xx} \varphi + D_{yy} \varphi)_{i,j}(t) = 0,
\end{split}
$\end{adjustbox}
\end{equation*}
which gives a semi-discrete scheme for the HJ PDE~\eqref{eqt:visc_HJ_initial} whose numerical Hamiltonian is defined in~\eqref{eqt:numericalH2d}.

With this discretization, the $\ell$-th update becomes
\begin{equation*}
\begin{adjustbox}{width=0.99\textwidth}$
\begin{dcases}
\rho_{i,j}^{\ell+1}(t) = \Big(\rho_{i,j}^{\ell}(t) + \tau_\rho\big(\dot{\tilde\varphi}^\ell_{i,j}(t) - f_{1,i,j,t}(\alpha^{\ell}_{11,i,j}(t))_+ (D_x^+ \tilde\varphi^\ell)_{i,j}(t) \\
\quad\quad\quad\quad
- f_{1,i,j,t}(\alpha^\ell_{12,i,j}(t))_- (D_x^- \tilde\varphi^\ell)_{i,j}(t) -f_{2,i,j,t}(\alpha^{\ell}_{21,i,j}(t))_+ (D_y^+ \tilde\varphi^\ell)_{i,j}(t) \\
\quad\quad\quad\quad
- f_{2,i,j,t}(\alpha^\ell_{22,i,j}(t))_- (D_y^- \tilde\varphi^\ell)_{i,j}(t) - \hat L_{i,j,t}(\alpha^\ell_{1,i}(t), \alpha^\ell_{2,i}(t)) - \epsilon(D_{xx}\tilde\varphi^\ell + D_{yy}\tilde\varphi^\ell)_{i,j}(t)\big)\Big)_+.\\
(\alpha^{\ell+1}_{11,i,j}(t), \alpha^{\ell+1}_{12,i,j}(t), \alpha^{\ell+1}_{21,i,j}(t), \alpha^{\ell+1}_{22,i,j}(t)) = \argmin_{\alpha_{11}, \alpha_{12}, \alpha_{21}, \alpha_{22}\in\R^m}\{f_{1,i,j,t}(\alpha_{11})_+ (D_x^+\tilde\varphi^\ell)_{i,j}(t) \\
\quad\quad\quad\quad
+ f_{1,i,j,t}(\alpha_{12})_- (D_x^-\tilde\varphi^\ell)_{i,j}(t) + f_{2,i,j,t}(\alpha_{21})_+ (D_y^+\tilde\varphi^\ell)_{i,j}(t) \\
\quad\quad\quad\quad
+ f_{2,i,j,t}(\alpha_{22})_- (D_y^-\tilde\varphi^\ell)_{i,j}(t) + \hat L_{i,j,t}(\alpha_{11}, \alpha_{12},\alpha_{21}, \alpha_{22})\\
\quad\quad\quad\quad
+ \frac{\rho_{i,j}^{\ell+1}(t)}{2\tau_\alpha} \left( \|\alpha_{11} - \alpha^\ell_{11,i,j}(t)\|^2 + \|\alpha_{12} - \alpha^\ell_{12,i,j}(t)\|^2 +  \|\alpha_{21} - \alpha^\ell_{21,i,j}(t)\|^2 +  \|\alpha_{22} - \alpha^\ell_{22,i,j}(t)\|^2\right)\}.\\
\varphi^{\ell+1}_{i,j}(t) = \varphi^{\ell}_{i,j}(t) + \tau_{\varphi}(I-\partial_{tt}-D_{xx} - D_{yy})^{-1}\Big(\dot{\rho}^{\ell+1}_{i,j}(t) - D_x^-\left(f_{1,i,j,t}(\alpha^{k+1}_{11,i,j}(t))_+\rho^{\ell+1}_{i,j}(t)\right) \\
\quad\quad\quad\quad - D_x^+\left(f_{1,i,j,t}(\alpha^{k+1}_{12,i,j}(t))_-\rho^{\ell+1}_{i,j}(t)\right) - D_y^-\left(f_{2,i,j,t}(\alpha^{k+1}_{21,i,j}(t))_+\rho^{\ell+1}_{i,j}(t)\right) \\
\quad\quad\quad\quad - D_y^+\left(f_{2,i,j,t}(\alpha^{k+1}_{22,i,j}(t))_-\rho^{\ell+1}_{i,j}(t)\right) + \epsilon(D_{xx}\rho^{\ell+1} + D_{yy}\rho^{\ell+1})_{i,j}(t) \Big).\\
\tilde\varphi^{\ell+1} = 2\varphi^{\ell+1} - \varphi^\ell.
\end{dcases}
$\end{adjustbox}
\end{equation*}

Then, we apply implicit Euler scheme for the time derivative of $\varphi$, and then the saddle point problem (divided by $\Delta t$) becomes
\begin{equation*}
\begin{adjustbox}{width=0.99\textwidth}$
\begin{split}
\min_{\substack{\varphi\\ \varphi_{i,j,1}=g(x_{i,j})}}\max_{\rho\geq 0, \alpha}  \sum_{k=2}^{n_t} \sum_{i=1}^{n_x} \sum_{j=1}^{n_y} \rho_{i,j,k}\Bigg((D_t^-\varphi)_{i,j,k} - f_{1,i,j,k}(\alpha_{11,i,j,k})_+ (D_x^+\varphi)_{i,j,k} 
- f_{1,i,j,k}(\alpha_{12,i,j,k})_- (D_x^-\varphi)_{i,j,k} \\
- f_{2,i,j,k}(\alpha_{21,i,j,k})_+ (D_y^+\varphi)_{i,j,k} - f_{2,i,j,k}(\alpha_{22,i,j,k})_- (D_y^-\varphi)_{i,j,k} - \epsilon(D_{xx}\varphi + D_{yy}\varphi)_{i,j,k} \\
- \hat L_{i,j,k}\left(\alpha_{11,i,j,k}, \alpha_{12,i,j,k}, \alpha_{21,i,j,k}, \alpha_{22,i,j,k}\right) \Bigg) - \frac{c}{\Delta t}\sum_{i=1}^{n_x} \sum_{j=1}^{n_y} \varphi_{i,j,n_t}.
\end{split}
$\end{adjustbox}
\end{equation*}
The corresponding algorithm in the $\ell$-th iteration becomes
\begin{equation*}
\begin{adjustbox}{width=0.99\textwidth}$
\begin{dcases}
\rho_{i,j,k}^{\ell+1} = \Big(\rho_{i,j,k}^{\ell} + \tau_\rho\big((D_t^-{\tilde\varphi^\ell})_{i,j,k} - f_{1,i,j,k}(\alpha^{\ell}_{11,i,j,k})_+ (D_x^+ \tilde\varphi^\ell)_{i,j,k} - f_{1,i,j,k}(\alpha^\ell_{12,i,j,k})_- (D_x^- \tilde\varphi^\ell)_{i,j,k} \\
\quad\quad\quad\quad - f_{2,i,j,k}(\alpha^{\ell}_{21,i,j,k})_+ (D_y^+ \tilde\varphi^\ell)_{i,j,k} - f_{2,i,j,k}(\alpha^\ell_{22,i,j,k})_- (D_y^- \tilde\varphi^\ell)_{i,j,k}\\
\quad\quad\quad\quad
- \hat L_{i,j,k}(\alpha^\ell_{11,i,j,k}, \alpha^\ell_{12,i,j,k},\alpha^\ell_{21,i,j,k}, \alpha^\ell_{22,i,j,k}) - \epsilon(D_{xx}\tilde\varphi^\ell + D_{yy}\tilde\varphi^\ell)_{i,j,k}\big)\Big)_+.\\
(\alpha^{\ell+1}_{11,i,j,k}, \alpha^{\ell+1}_{12,i,j,k}, \alpha^{\ell+1}_{21,i,j,k}, \alpha^{\ell+1}_{22,i,j,k}) = \argmin_{\alpha_{11}, \alpha_{12}, \alpha_{21}, \alpha_{22}\in\R^m}\{f_{1,i,j,k}(\alpha_{11})_+ (D_x^+\tilde\varphi^\ell)_{i,j,k} \\
\quad\quad\quad\quad
+ f_{1,i,j,k}(\alpha_{12})_- (D_x^-\tilde\varphi^\ell)_{i,j,k} + f_{2,i,j,k}(\alpha_{21})_+ (D_y^+\tilde\varphi^\ell)_{i,j,k} \\
\quad\quad\quad\quad
+ f_{2,i,j,k}(\alpha_{22})_- (D_y^-\tilde\varphi^\ell)_{i,j,k} + \hat L_{i,j,k}(\alpha_{11}, \alpha_{12}, \alpha_{21}, \alpha_{22})\\
\quad\quad\quad\quad
+ \frac{\rho_{i,j,k}^{\ell+1}}{2\tau_\alpha} \left( \|\alpha_{11} - \alpha^\ell_{11,i,j,k}\|^2 + \|\alpha_{12} - \alpha^\ell_{12,i,j,k}\|^2 +  \|\alpha_{21} - \alpha^\ell_{21,i,j,k}\|^2 +  \|\alpha_{22} - \alpha^\ell_{22,i,j,k}\|^2\right)\}.\\
\varphi^{\ell+1}_{i,j,k} = \varphi^{\ell}_{i,j,k} + \tau_{\varphi}(I-D_{tt}-D_{xx} - D_{yy})^{-1}\Big((D_t^+{\rho}^{\ell+1})_{i,j,k} - D_x^-\left(f_{1,i,j,k}(\alpha^{k+1}_{11,i,j,k})_+\rho^{\ell+1}_{i,j,k}\right) \\
\quad\quad\quad\quad - D_x^+\left(f_{1,i,j,k}(\alpha^{k+1}_{12,i,j,k})_-\rho^{\ell+1}_{i,j,k}\right) - D_y^-\left(f_{2,i,j,k}(\alpha^{k+1}_{21,i,j,k})_+\rho^{\ell+1}_{i,j,k}\right) \\
\quad\quad\quad\quad - D_y^+\left(f_{2,i,j,k}(\alpha^{k+1}_{22,i,j,k})_-\rho^{\ell+1}_{i,j,k}\right) + \epsilon(D_{xx}\rho^{\ell+1} + D_{yy}\rho^{\ell+1})_{i,j,k}\Big).\\
\tilde\varphi^{\ell+1} = 2\varphi^{\ell+1} - \varphi^\ell.
\end{dcases}
$\end{adjustbox}
\end{equation*}
Just as with the one-dimensional scenario, by selecting the numerical Lagrangian as $L_{x,t}(\alpha_{11}) + L_{x,t}(\alpha_{12}) + L_{x,t}(\alpha_{21}) + L_{x,t}(\alpha_{22})$, the update process for $\alpha$ can be carried out in parallel across the variables $\alpha_{11}$, $\alpha_{12}$, $\alpha_{21}$, and $\alpha_{22}$.

\end{document}